\documentclass[9pt]{article}
\usepackage{amssymb}
\usepackage{amsmath}
\usepackage{amsthm}

\input pictex.tex

\newcommand{\vsp}{\vspace{0.1cm}}

\newcommand{\esp}{\hspace{0.05cm}}
\theoremstyle{definition}

\newtheorem{thm}{Theorem}[section]

\newtheorem{prop}[thm]{Proposition}

\newtheorem{lem}[thm]{Lemma}

\newtheorem{rem}[thm]{Remark}

\newcommand{\dac}{\mathrm{Diff}_+^{1+\alpha}([0,1])}

\title{Sharp regularity for certain nilpotent group actions on the interval}
\author{G.~Castro, \esp E.~Jorquera \esp \& \esp A. Navas}
\date{}

\begin{document}

\maketitle

\vspace{-0.8cm}

\begin{abstract}

\vspace{0.1cm}

According to the classical Plante-Thurston Theorem, all nilpotent groups of
$C^2$-diffeomorphisms of the closed interval are Abelian. Using techniques
coming from the works of Denjoy and Pixton, Farb and Franks constructed a
faithful action by $C^1$-diffeomorphisms of $[0,1]$ for every finitely-generated,
torsion-free, non-Abelian nilpotent group. In this work, we give
a version of this construction that is sharp in what concerns the H\"older
regularity of the derivatives. Half of the proof relies on results on
random paths on Heisenberg-like groups that are interesting by themselves.
\end{abstract}

\tableofcontents


\section{Introduction}

\hspace{0.45cm} Much work has been done on centralizers of $C^2$-diffeomorphisms of the
interval \cite{elena,Kop,serg,sz}. This theory has been extensively used for studying the
algebraic constraints  of finitely-generated subgroups of $\mathrm{Diff}_+^2([0,1])$.
For example, using the famous Kopell lemma \cite{Kop}, Plante and Thurston showed
that nilpotent groups of $C^2$-diffeomorphisms of $[0,1[$ (resp. $]0,1[$) are
Abelian (resp. metabelian); see \cite{PT}.

As is well known, most of the rigidity properties are lost when we consider
centralizers of $C^1$-diffeomor- phisms. In relation to Plante-Thurston's theorem,
this fact is corroborated by the work of Farb and Franks. In \cite{FF}, they construct
an embedding $\phi_{_{FF}}$ of $N_d$ into $\mathrm{Diff}^1_+([0,1])$, where $N_d$ denotes
the (nilpotent) group of $(d+1) \times (d+1)$ lower-triangular matrices whose entries are
integers which equal 1 on the diagonal (see \S \ref{accion} for the details). Since every
finitely-generated, torsion-free, nilpotent group embeds into $N_d$ for some $d \geq 1$
(see \cite{Ra}), one concludes that all these groups can be realized as groups of
$C^1$-diffeomorphisms of the (closed) interval (compare \cite{jorquera}).

Major progress has been recently made in the understanding of the loss of rigidity
for centralizers in intermediate differentiability classes, that is, between $C^1$
and $C^2$ (see \cite{DKN,KN,growth}). Recall that, for $0 < \alpha < 1$, a
diffeomorphism $f$ is said to be of class $C^{1+\alpha}$ if its derivative is
$\alpha$-H\"older continuous. In other words, there exists a constant $M$ such that
for all $x,y$,
\begin{equation}\label{def-holder}
|f'(x) - f'(y)| \leq M |x - y|^{\alpha}.
\end{equation}
We denote the group of $C^{1+\alpha}$-diffeomorphisms of $[0,1]$ by
$\mathrm{Diff}^{1+\alpha}_+([0,1])$. In the first part of this work
we show the following result. (Notice that for $d = 2$, the theorem
below still holds and follows from Plante-Thurston's theorem.)

\vspace{0.5cm}

\noindent{\bf Theorem A.} {\em If $d \geq 3$ and \esp $\alpha > \frac{2}{d(d-1)}$,
\esp then the action $\phi_{_{FF}}$ is not topologically conjugated to an action
by $C^{1+\alpha}$-diffeomorphisms of \esp $[0,1]$.}

\vspace{0.5cm}

This theorem should be considered as a partial complement to \cite[Theorem B]{growth}
which establishes that, for all $0 < \alpha < 1$, every subgroup $\Gamma$ of
$\mathrm{Diff}_+^{1+\alpha}([0,1])$ without free subsemigroups is virtually nilpotent.
(Although the last result still holds for the open interval $]0,1[$, Theorem A above
fails to be true in this context, but it extends --with the very same proof-- to
the case of the half-closed interval). For the proof of our theorem, the main
technical achievement consists in controlling the distortion of suitable compositions
of elements in any regularity larger than the critical one. To do this, we develop a
nontrivial modification of the probabilistic techniques of \cite{DKN,KN}. Recall that
\cite[Theorem B]{DKN} deals with Abelian group actions that are dynamically very
similar to $\phi_{_{FF}}$, and a direct application of it shows that $\phi_{_{FF}}$
is not conjugated to an action by $C^{1+\alpha}$-diffeomorphisms of $[0,1[$ for any
$\alpha > \frac{1}{d-1}$. The fact that our critical regularity here is actually smaller
relies on that compared to the Abelian actions of \cite{DKN}, the action $\phi_{_{FF}}$
has a more complicated combinatorial dynamics in that the growth of certain orbits is
polynomial with degree precisely equal to \esp $\frac{d(d-1)}{2}$. We should point
out that similar combinatorial dynamics appear for the actions of the natural quotients
of the Grigorchuk-Machi's group \cite{grig} for which the method of this article should
also provide the best possible regularity (compare \cite[Theorem A]{growth}). Moreover,
it is worth mentioning that the very same arguments show that Theorem A above still
applies to topological {\em semiconjugacies.}

Although not directly related, all the results described above should be compared to
(and have potential relations with) Borichev's extension \cite{borichev} to intermediate
regularity of Polterovich-Sodin's theorem \cite{PS} concerning distortion of interval
diffeomorphisms.

\vsp\vsp

The second part of this work is devoted to a converse of Theorem A.
The next theorem improves the main result of \cite{FF}.

\vspace{0.5cm}

\noindent{\bf Theorem B.} {\em For each $d \geq 2$ and  \esp
$\alpha < \frac{2}{d(d-1)}$, \esp the action $\phi_{_{FF}}$ is topologically
conjugated to an action by $C^{1+\alpha}$-diffeomorphisms of \esp $[0,1]$.}

\vspace{0.5cm}

The proof of this theorem is based on classical constructions of Denjoy and Pixton
(a clever exposition of these techniques appears in \cite{TsP}; see also
\cite{navas-book}). Nevertheless, putting these methods in practice in the
present case is far from being straightforward. The computations are quite
involved, and in this part of the  work some of them are just sketched.

As in \cite{DKN,KN}, here we were unable to settle the
$C^{1 + \frac{2}{d(d-1)}}$ case, though we conjecture that the
rigidity ({\em i.e.} Theorem A) still holds for this critical regularity.

Theorems A and B suggest that, attached to each finitely-generated, torsion-free
nilpotent group $\Gamma$, there should be a positive exponent $\alpha(\Gamma) \leq 1$
that is critical for embedding $\Gamma$ into $\mathrm{Diff}_+^{1+\alpha}([0,1])$. However, 
it is still unclear to us what should be the value of $\alpha(\Gamma)$. Indeed, Theorem B
only deals with very particular actions, and many nilpotent groups admit actions that
are fairly different from these. In order to corroborate this point, in the last
part of this work we improve another construction of \cite{FF}, thus proving
the next

\vspace{0.35cm}

\noindent{\bf Theorem C.} {\em For every $\alpha < 1$ and each $d \geq 1$, 
the group $\dac$ contains a metabelian subgroup of nilpotence degree $d$.}


\vspace{0.25cm}

\noindent{\bf Acknowledgments.} We are indebted to Yves de Cornulier and Romain Tessera
for a crucial remark concerning the growth of certain orbits of the action of $N_d$, to
Bertrand Deroin and Victor Kleptsyn for many valuable remarks to the second part of 
this article, and to Takashi Tsuboi for a clever hint for the constructions therein. The 
second-named author was funded by the Fondecyt Iniciaci\'on Project 11121316, 
while third-named author was partially funded by the Fondecyt Project 1100536.


\section{Non-existence of smoothing for $\alpha > \frac{2}{d(d-1)}$}

\subsection{A remind on Farb-Franks' action $\phi_{_{FF}}$}
\label{accion}

\hspace{0.45cm} We deal with the group $N_d$ of \esp $(d+1) \times (d+1)$ \esp
lower-triangular matrices with integer entries, all of which are equal to 1 on the
diagonal. Notice that $N_2$ corresponds to the Heisenberg group. In general, $N_d$ is
a nilpotent group of nilpotence degree $d$. A nice system of generators of $N_{d}$
is $\{f_{2,1},\ldots,f_{d+1,d}\}$, where $f_{i,j}$ is the elementary matrix whose
unique nonzero entry outside the diagonal is the $(i,j)$-entry (with $i > j$).

The group $N_d$ acts linearly on $\mathbb{Z}^{d+1}$ with the affine hyperplane
$1 \times \mathbb{Z}^{d}$ remaining invariant. The thus-induced action on $\mathbb{Z}^{d}$
allows producing an action on the interval as follows. Let $\big\{I_{i_1,\ldots,i_{d}}\!\!:
\esp\esp (i_1,\ldots,i_{d}) \!\in\! \mathbb{Z}^d \big\}$ be a family of intervals such
that the sum $\sum_{i_1,\ldots,i_{d}} \vert I_{i_1,\ldots,i_{d}}\vert$ is finite,
say equal to 1 after normalization. We join these intervals lexicographically on the
closed interval $[0,1]$, and we identify $f_{j+1,j}$ to a certain homeomorphism
sending each interval $I = I_{i_1,\ldots,i_{d}}$ into the interval $J$ given by:

\vspace{0.15cm}

\noindent $\bullet$ \esp $J := I_{i_1+1, i_{2},\ldots,i_{d-1},i_{d}}$, \esp for $j = 1$,

\vspace{0.15cm}

\noindent $\bullet$ \esp $J := I_{i_1,\ldots,i_{j-1}, i_j + i_{j-1},i_{j+1},\ldots,i_{d}}$,
for $2 \!\leq\! j \!\leq\! d$.

\vspace{0.15cm}

\noindent It is not hard to perform this procedure in a equivariant way (for instance, using
piecewise-affine maps), thus preserving the group structure and hence obtaining an embedding
of $N_d$ into $\mathrm{Homeo}_+([0,1])$. (Much harder is to obtain an embedding into the group
of diffeomorphisms.) For this action, an interval of the form $I_{i_1,\ldots,i_{d}}$ is sent
by $f \in N_d$ into $I_{j_1,\ldots,j_{d}}$, where \esp $f \big( (1,i_1,\ldots,i_{d})^T \big)
= (1, j_1, \ldots, j_{d})^T$. \esp
Notice that up to topological conjugacy, all the actions obtained by this
procedure are equivalent. This includes Farb-Franks' action $\phi_{_{FF}}$, which is obtained
via this method for a well-chosen family of diffeomorphisms between the intervals
of type $I,J$ above so that the resulting $f_{i,j}$'s are $C^1$-diffeomorphisms.


\subsection{From control of distortion to the proof of Theorem A}
\label{main-idea}

\hspace{0.45cm} Let us begin by stating a general principle from \cite{DKN} in the
form of the following

\vspace{0.1cm}

\begin{prop} {\em Let $f_1,\ldots,f_k$ be $C^1$-diffeomorphisms of the interval
$[0,1]$ that commute with a $C^1$-diffeomorphism $g$. Assume that $g$ fixes a
subinterval $I$ of $[0,1]$ and its restriction to $I$ is nontrivial. Assume
moreover that for a certain $0 < \alpha < 1$ and a sequence of indexes
$i_j \in \{ 1,\ldots,k \}$, the sum}
\begin{equation}\label{clave}
L_{\alpha} := \sum_{j \geq 0} \big| f_{i_j} \cdots f_{i_1} (I) \big|^{\alpha}
\end{equation}
{\em is finite. Then $f_1,\ldots,f_k$ cannot be all of class $C^{1+\alpha}$.}
\end{prop}

\noindent{\bf Proof.} Let $x_0 \in I$ be such that $g(x_0) \neq x_0$. Denote
by $[a,b]$ the shortest interval containing $x_0$ that is fixed by $g$.
For each $j \geq 1, n \geq 1$ and $z \in [a,b]$, the equality \esp
$g^n = (f_{i_j} \cdots f_{i_1})^{-1} \circ g^n \circ (f_{i_j} \cdots f_{i_1})$
\esp yields
$$\log Dg^n (z) = \log D(f_{i_j} \cdots f_{i_1})(z) +
\log Dg^n (f_{i_j} \cdots f_{i_1} (z)) - \log D(f_{i_j} \cdots f_{i_1}) (g^n (z)).$$
Fix a constant $M$ such that (\ref{def-holder}) holds for all $f \in \{ f_1,\ldots,f_k \}$
and all $x,y$ in $[0,1]$. Letting $z_{n} := g^n(z)$ and noticing that
$z_n$ belongs to $[a,b] \subset I$ for all $n \geq 1$, we obtain
\begin{eqnarray*}
|\log Dg^n(z)|
&\leq& |\log Dg^n(f_{i_j} \cdots f_{i_1} (z))|
     + \sum_{m=1}^{j} \big| \log Df_{i_m} (f_{i_{m-1}} \cdots f_{i_1}(z))
     - \log Df_{i_m} (f_{i_{m-1}} \cdots f_{i_1}(z_n)) \big|\\
&\leq& |\log Dg^n (f_{i_j} \cdots f_{i_1} (z))| +
\sum_{m=1}^{j} M \big| f_{i_{m-1}} \cdots f_{i_1} (z)
    - f_{i_{m-1}} \cdots f_{i_1} (z_n) \big|^{\alpha}\\
&\leq& |\log Dg^n (f_{i_j} \cdots f_{i_1} (z))| +
     M \sum_{m=1}^{j} | f_{i_{m-1}} \cdots f_{i_1} (I)|^{\alpha}\\
&\leq& |\log Dg^n (f_{i_j} \cdots f_{i_1} (z))| + M L_{\alpha}.
\end{eqnarray*}
The length of the intervals \esp $f_{i_j} \cdots f_{i_1} (I)$ \esp
must necessarily converge to zero as $j$ goes to infinite. Moreover, since $g^n$ fixes $I$
and commutes with $f_1,\ldots,f_k$, on each of these intervals there must be a point at which
its derivative equals 1. By the continuity of $D g^n$, we conclude that the value of \esp
$Dg^n (f_{i_j} \cdots f_{i_1} (z))$ \esp converges to 1 as $j$ goes to infinite. Hence we obtain
\esp $Dg^n (z) \leq e^{M L_{\alpha}}$ \esp for all \esp $n \geq 1$ \esp and all \esp $z \in [a,b]$, \esp
which certainly contradicts the fact that the restriction of $g$ to $[a,b]$ is nontrivial. $\hfill\square$

\vspace{0.43cm}

Let us come back to the action $\phi_{_{FF}}$. Notice that the group $N_{d-1}$ can be naturally
viewed as the subgroup of $N_d$ formed by the elements whose last row coincide with that
of the identity. We will denote by $N_{d-1}^*$ the copy of $N_{d-1}$ inside $N_d$.

Notice that the element $g := f_{d+1,1} \in N_{d}$ is centralized by $N_{d-1}^*$.
Under the action $\phi_{_{FF}}$, this element fixes the interval
\begin{equation}\label{def-interval}
I^* := \bigcup_{j \in \mathbb{Z}} I_{0,\ldots,0,j}.
\end{equation}
Moreover, this interval is sent into a disjoint one by any nontrivial
element of $N_{d-1}^*$. We are hence in a situation close to that of the preceding
proposition. Thus, we need to ensure the existence of a systems of generators for
$N_{d-1}^*$ and a sequence of compositions for which the associated sum (\ref{clave})
is finite provided that
$\alpha \!> \!\frac{2}{d(d-1)}$. To do this, we will use the system of
generators $\{f_{2,1}, f_{3,1}, \ldots , f_{d,1} \} \cup \{f_{2,1}, f_{3,2}, \ldots,
f_{d,d-1} \}$.

\vspace{0.1cm}

It is worth mentioning that this is an analogous problem to that of the $\mathbb{Z}^d$-actions
on the interval considered in \cite[Th\'eor\`eme B]{DKN}. However, the $\mathbb{Z}^d$-case
is easier in that the generators of the dynamics commute, hence the orbit graph of the
associated interval $I^*$ has a simpler structure. Indeed, the space of infinite paths of this
graph can be endowed of a natural probability measure such that for appropriately large values
of $\alpha$ (namely, for $\alpha > 1/d$), almost every path has a finite $L_{\alpha}$-series.
In order to establish this, besides the restriction on the exponent $\alpha$, the main
property of the underlying process is that the arrival probabilities up to time $k$
are equidistributed along the sphere of radius $k$
(centered at the origin) for every $k \geq 1$. Although in \cite{DKN} this is modeled
via a Polya urn like model that charges only the positive powers of the generators, an
alternative model sharing this property that charges both positive and negative powers
of the generators is the Markov process depicted in Figure~1 below for the case $d=2$
(the reader will easily check the equidistribution property along spheres as well as
the general rule for the transition probabilities; the generalization for higher
values of $d$ is not very hard).

\begin{rem} It seems to be an interesting and nontrivial problem to determine general
conditions for an infinite graph ensuring the existence of a Markov process satisfying
the equidistribution property above.
\end{rem}

\vspace{0.5cm}


\beginpicture

\setcoordinatesystem units <1.2cm,1.2cm>

\putrule from -3 -4 to -3 4
\putrule from -2 -4 to -2 4
\putrule from -1 -4 to -1 4
\putrule from 0 -4 to 0 4
\putrule from 1 -4 to 1 4
\putrule from 2 -4 to 2 4
\putrule from 3 -4 to 3 4

\putrule from -4 -3 to 4 -3
\putrule from -4 -2 to 4 -2
\putrule from -4 -1 to 4 -1
\putrule from -4 0 to 4 0
\putrule from -4 1 to 4 1
\putrule from -4 2 to 4 2
\putrule from -4 3 to 4 3


\plot 0.5 3
0.35 3.05 /
\plot 0.5 3
0.35 2.95 /

\plot 0.5 2
0.35 2.05 /
\plot 0.5 2
0.35 1.95 /

\plot 0.5 1
0.35 1.05 /
\plot 0.5 1
0.35 0.95 /

\plot 0.5 0
0.35 0.05 /
\plot 0.5 0
0.35 -0.05 /

\plot 0.5 -1
0.35 -0.95 /
\plot 0.5 -1
0.35 -1.05 /

\plot 0.5 -2
0.35 -1.95 /
\plot 0.5 -2
0.35 -2.05 /

\plot 0.5 -3
0.35 -2.95 /
\plot 0.5 -3
0.35 -3.05 /


\plot 1.5 3
1.35 3.05 /
\plot 1.5 3
1.35 2.95 /

\plot 1.5 2
1.35 2.05 /
\plot 1.5 2
1.35 1.95 /

\plot 1.5 1
1.35 1.05 /
\plot 1.5 1
1.35 0.95 /

\plot 1.5 0
1.35 0.05 /
\plot 1.5 0
1.35 -0.05 /

\plot 1.5 -1
1.35 -0.95 /
\plot 1.5 -1
1.35 -1.05 /

\plot 1.5 -2
1.35 -1.95 /
\plot 1.5 -2
1.35 -2.05 /

\plot 1.5 -3
1.35 -2.95 /
\plot 1.5 -3
1.35 -3.05 /


\plot 2.5 3
2.35 3.05 /
\plot 2.5 3
2.35 2.95 /

\plot 2.5 2
2.35 2.05 /
\plot 2.5 2
2.35 1.95 /

\plot 2.5 1
2.35 1.05 /
\plot 2.5 1
2.35 0.95 /

\plot 2.5 0
2.35 0.05 /
\plot 2.5 0
2.35 -0.05 /

\plot 2.5 -1
2.35 -0.95 /
\plot 2.5 -1
2.35 -1.05 /

\plot 2.5 -2
2.35 -1.95 /
\plot 2.5 -2
2.35 -2.05 /

\plot 2.5 -3
2.35 -2.95 /
\plot 2.5 -3
2.35 -3.05 /


\plot 3.5 3
3.35 3.05 /
\plot 3.5 3
3.35 2.95 /

\plot 3.5 2
3.35 2.05 /
\plot 3.5 2
3.35 1.95 /

\plot 3.5 1
3.35 1.05 /
\plot 3.5 1
3.35 0.95 /

\plot 3.5 0
3.35 0.05 /
\plot 3.5 0
3.35 -0.05 /

\plot 3.5 -1
3.35 -0.95 /
\plot 3.5 -1
3.35 -1.05 /

\plot 3.5 -2
3.35 -1.95 /
\plot 3.5 -2
3.35 -2.05 /

\plot 3.5 -3
3.35 -2.95 /
\plot 3.5 -3
3.35 -3.05 /


\plot -0.5 3
-0.35 3.05 /
\plot -0.5 3
-0.35 2.95 /

\plot -0.5 2
-0.35 2.05 /
\plot -0.5 2
-0.35 1.95 /

\plot -0.5 1
-0.35 1.05 /
\plot -0.5 1
-0.35 0.95 /

\plot -0.5 0
-0.35 0.05 /
\plot -0.5 0
-0.35 -0.05 /

\plot -0.5 -1
-0.35 -0.95 /
\plot -0.5 -1
-0.35 -1.05 /

\plot -0.5 -2
-0.35 -1.95 /
\plot -0.5 -2
-0.35 -2.05 /

\plot -0.5 -3
-0.35 -2.95 /
\plot -0.5 -3
-0.35 -3.05 /


\plot -1.5 3
-1.35 3.05 /
\plot -1.5 3
-1.35 2.95 /

\plot -1.5 2
-1.35 2.05 /
\plot -1.5 2
-1.35 1.95 /

\plot -1.5 1
-1.35 1.05 /
\plot -1.5 1
-1.35 0.95 /

\plot -1.5 0
-1.35 0.05 /
\plot -1.5 0
-1.35 -0.05 /

\plot -1.5 -1
-1.35 -0.95 /
\plot -1.5 -1
-1.35 -1.05 /

\plot -1.5 -2
-1.35 -1.95 /
\plot -1.5 -2
-1.35 -2.05 /

\plot -1.5 -3
-1.35 -2.95 /
\plot -1.5 -3
-1.35 -3.05 /


\plot -2.5 3
-2.35 3.05 /
\plot -2.5 3
-2.35 2.95 /

\plot -2.5 2
-2.35 2.05 /
\plot -2.5 2
-2.35 1.95 /

\plot -2.5 1
-2.35 1.05 /
\plot -2.5 1
-2.35 0.95 /

\plot -2.5 0
-2.35 0.05 /
\plot -2.5 0
-2.35 -0.05 /

\plot -2.5 -1
-2.35 -0.95 /
\plot -2.5 -1
-2.35 -1.05 /

\plot -2.5 -2
-2.35 -1.95 /
\plot -2.5 -2
-2.35 -2.05 /

\plot -2.5 -3
-2.35 -2.95 /
\plot -2.5 -3
-2.35 -3.05 /


\plot -3.5 3
-3.35 3.05 /
\plot -3.5 3
-3.35 2.95 /

\plot -3.5 2
-3.35 2.05 /
\plot -3.5 2
-3.35 1.95 /

\plot -3.5 1
-3.35 1.05 /
\plot -3.5 1
-3.35 0.95 /

\plot -3.5 0
-3.35 0.05 /
\plot -3.5 0
-3.35 -0.05 /

\plot -3.5 -1
-3.35 -0.95 /
\plot -3.5 -1
-3.35 -1.05 /

\plot -3.5 -2
-3.35 -1.95 /
\plot -3.5 -2
-3.35 -2.05 /

\plot -3.5 -3
-3.35 -2.95 /
\plot -3.5 -3
-3.35 -3.05 /


\plot 3 0.5
3.05 0.35 /
\plot 3 0.5
2.95 0.35 /

\plot 2 0.5
2.05 0.35 /
\plot 2 0.5
1.95 0.35 /

\plot 1 0.5
1.05 0.35 /
\plot 1 0.5
0.95 0.35 /

\plot 0 0.5
0.05 0.35 /
\plot 0 0.5
-0.05 0.35 /

\plot -1 0.5
-0.95 0.35 /
\plot -1 0.5
-1.05 0.35 /

\plot -2 0.5
-1.95 0.35 /
\plot -2 0.5
-2.05 0.35 /

\plot -3 0.5
-2.95 0.35 /
\plot -3 0.5
-3.05 0.35 /


\plot 3 1.5
3.05 1.35 /
\plot 3 1.5
2.95 1.35 /

\plot 2 1.5
2.05 1.35 /
\plot 2 1.5
1.95 1.35 /

\plot 1 1.5
1.05 1.35 /
\plot 1 1.5
0.95 1.35 /

\plot 0 1.5
0.05 1.35 /
\plot 0 1.5
-0.05 1.35 /

\plot -1 1.5
-0.95 1.35 /
\plot -1 1.5
-1.05 1.35 /

\plot -2 1.5
-1.95 1.35 /
\plot -2 1.5
-2.05 1.35 /

\plot -3 1.5
-2.95 1.35 /
\plot -3 1.5
-3.05 1.35 /


\plot 3 2.5
3.05 2.35 /
\plot 3 2.5
2.95 2.35 /

\plot 2 2.5
2.05 2.35 /
\plot 2 2.5
1.95 2.35 /

\plot 1 2.5
1.05 2.35 /
\plot 1 2.5
0.95 2.35 /

\plot 0 2.5
0.05 2.35 /
\plot 0 2.5
-0.05 2.35 /

\plot -1 2.5
-0.95 2.35 /
\plot -1 2.5
-1.05 2.35 /

\plot -2 2.5
-1.95 2.35 /
\plot -2 2.5
-2.05 2.35 /

\plot -3 2.5
-2.95 2.35 /
\plot -3 2.5
-3.05 2.35 /


\plot 3 3.5
3.05 3.35 /
\plot 3 3.5
2.95 3.35 /

\plot 2 3.5
2.05 3.35 /
\plot 2 3.5
1.95 3.35 /

\plot 1 3.5
1.05 3.35 /
\plot 1 3.5
0.95 3.35 /

\plot 0 3.5
0.05 3.35 /
\plot 0 3.5
-0.05 3.35 /

\plot -1 3.5
-0.95 3.35 /
\plot -1 3.5
-1.05 3.35 /

\plot -2 3.5
-1.95 3.35 /
\plot -2 3.5
-2.05 3.35 /

\plot -3 3.5
-2.95 3.35 /
\plot -3 3.5
-3.05 3.35 /


\plot 3 -0.5
3.05 -0.35 /
\plot 3 -0.5
2.95 -0.35 /

\plot 2 -0.5
2.05 -0.35 /
\plot 2 -0.5
1.95 -0.35 /

\plot 1 -0.5
1.05 -0.35 /
\plot 1 -0.5
0.95 -0.35 /

\plot 0 -0.5
0.05 -0.35 /
\plot 0 -0.5
-0.05 -0.35 /

\plot -1 -0.5
-0.95 -0.35 /
\plot -1 -0.5
-1.05 -0.35 /

\plot -2 -0.5
-1.95 -0.35 /
\plot -2 -0.5
-2.05 -0.35 /

\plot -3 -0.5
-2.95 -0.35 /
\plot -3 -0.5
-3.05 -0.35 /


\plot 3 -1.5
3.05 -1.35 /
\plot 3 -1.5
2.95 -1.35 /

\plot 2 -1.5
2.05 -1.35 /
\plot 2 -1.5
1.95 -1.35 /

\plot 1 -1.5
1.05 -1.35 /
\plot 1 -1.5
0.95 -1.35 /

\plot 0 -1.5
0.05 -1.35 /
\plot 0 -1.5
-0.05 -1.35 /

\plot -1 -1.5
-0.95 -1.35 /
\plot -1 -1.5
-1.05 -1.35 /

\plot -2 -1.5
-1.95 -1.35 /
\plot -2 -1.5
-2.05 -1.35 /

\plot -3 -1.5
-2.95 -1.35 /
\plot -3 -1.5
-3.05 -1.35 /


\plot 3 -2.5
3.05 -2.35 /
\plot 3 -2.5
2.95 -2.35 /

\plot 2 -2.5
2.05 -2.35 /
\plot 2 -2.5
1.95 -2.35 /

\plot 1 -2.5
1.05 -2.35 /
\plot 1 -2.5
0.95 -2.35 /

\plot 0 -2.5
0.05 -2.35 /
\plot 0 -2.5
-0.05 -2.35 /

\plot -1 -2.5
-0.95 -2.35 /
\plot -1 -2.5
-1.05 -2.35 /

\plot -2 -2.5
-1.95 -2.35 /
\plot -2 -2.5
-2.05 -2.35 /

\plot -3 -2.5
-2.95 -2.35 /
\plot -3 -2.5
-3.05 -2.35 /


\plot 3 -3.5
3.05 -3.35 /
\plot 3 -3.5
2.95 -3.35 /

\plot 2 -3.5
2.05 -3.35 /
\plot 2 -3.5
1.95 -3.35 /

\plot 1 -3.5
1.05 -3.35 /
\plot 1 -3.5
0.95 -3.35 /

\plot 0 -3.5
0.05 -3.35 /
\plot 0 -3.5
-0.05 -3.35 /

\plot -1 -3.5
-0.95 -3.35 /
\plot -1 -3.5
-1.05 -3.35 /

\plot -2 -3.5
-1.95 -3.35 /
\plot -2 -3.5
-2.05 -3.35 /

\plot -3 -3.5
-2.95 -3.35 /
\plot -3 -3.5
-3.05 -3.35 /

\put{$\bullet$} at 0 0

\small

\put{Figure 1} at -0.04 -4.25

\put{} at -6.9 0

\tiny

\put{$\bullet$} at -3 -3
\put{$\bullet$} at -2 -3
\put{$\bullet$} at -1 -3
\put{$\bullet$} at 0 -3
\put{$\bullet$} at 1 -3
\put{$\bullet$} at 2 -3
\put{$\bullet$} at 3 -3

\put{$\bullet$} at -3 -2
\put{$\bullet$} at -2 -2
\put{$\bullet$} at -1 -2
\put{$\bullet$} at 0 -2
\put{$\bullet$} at 1 -2
\put{$\bullet$} at 2 -2
\put{$\bullet$} at 3 -2

\put{$\bullet$} at -3 -1
\put{$\bullet$} at -2 -1
\put{$\bullet$} at -1 -1
\put{$\bullet$} at 0 -1
\put{$\bullet$} at 1 -1
\put{$\bullet$} at 2 -1
\put{$\bullet$} at 3 -1

\put{$\bullet$} at -3 0
\put{$\bullet$} at -2 0
\put{$\bullet$} at -1 0
\put{$\bullet$} at 1 0
\put{$\bullet$} at 2 0
\put{$\bullet$} at 3 0

\put{$\bullet$} at -3 1
\put{$\bullet$} at -2 1
\put{$\bullet$} at -1 1
\put{$\bullet$} at 0 1
\put{$\bullet$} at 1 1
\put{$\bullet$} at 2 1
\put{$\bullet$} at 3 1

\put{$\bullet$} at -3 2
\put{$\bullet$} at -2 2
\put{$\bullet$} at -1 2
\put{$\bullet$} at 0 2
\put{$\bullet$} at 1 2
\put{$\bullet$} at 2 2
\put{$\bullet$} at 3 2

\put{$\bullet$} at -3 3
\put{$\bullet$} at -2 3
\put{$\bullet$} at -1 3
\put{$\bullet$} at 0 3
\put{$\bullet$} at 1 3
\put{$\bullet$} at 2 3
\put{$\bullet$} at 3 3

\put{$\frac{1}{4}$} at 0.5 0.2
\put{$\frac{1}{4}$} at -0.5 0.2
\put{$\frac{1}{4}$} at 0.2 0.5
\put{$\frac{1}{4}$} at 0.2 -0.5


\put{$\frac{1}{2}$} at 1.5 0.2
\put{$\frac{1}{4}$} at 1.2 0.5
\put{$\frac{1}{4}$} at 1.2 -0.5

\put{$\frac{1}{2}$} at -1.5 0.2
\put{$\frac{1}{4}$} at -0.8 0.5
\put{$\frac{1}{4}$} at -0.8 -0.5

\put{$\frac{1}{2}$} at 0.2 1.5
\put{$\frac{1}{4}$} at 0.5 1.2
\put{$\frac{1}{4}$} at -0.5 1.2

\put{$\frac{1}{2}$} at 0.2 -1.5
\put{$\frac{1}{4}$} at 0.5 -0.8
\put{$\frac{1}{4}$} at -0.5 -0.8


\put{$\frac{2}{3}$} at 2.5 0.2
\put{$\frac{1}{6}$} at 2.2 0.5
\put{$\frac{1}{6}$} at 2.2 -0.5

\put{$\frac{2}{3}$} at -2.5 0.2
\put{$\frac{1}{6}$} at -1.8 0.5
\put{$\frac{1}{6}$} at -1.8 -0.5

\put{$\frac{2}{3}$} at 0.2 2.5
\put{$\frac{1}{6}$} at 0.5 2.2
\put{$\frac{1}{6}$} at -0.5 2.2

\put{$\frac{2}{3}$} at 0.2 -2.5
\put{$\frac{1}{6}$} at 0.5 -1.8
\put{$\frac{1}{6}$} at -0.5 -1.8


\put{$\frac{1}{2}$} at 1.5 1.2
\put{$\frac{1}{2}$} at 1.2 1.5

\put{$\frac{1}{2}$} at 1.5 -0.8
\put{$\frac{1}{2}$} at 1.2 -1.5

\put{$\frac{1}{2}$} at -0.8 1.5
\put{$\frac{1}{2}$} at -1.5 1.2

\put{$\frac{1}{2}$} at -0.8 -1.5
\put{$\frac{1}{2}$} at -1.5 -0.8


\put{$\frac{3}{4}$} at 3.5 0.2
\put{$\frac{1}{8}$} at 3.2 0.5
\put{$\frac{1}{8}$} at 3.2 -0.5

\put{$\frac{3}{4}$} at -3.5 0.2
\put{$\frac{1}{8}$} at -2.8 0.5
\put{$\frac{1}{8}$} at -2.8 -0.5

\put{$\frac{3}{4}$} at 0.2 3.5
\put{$\frac{1}{8}$} at 0.5 3.2
\put{$\frac{1}{8}$} at -0.5 3.2

\put{$\frac{3}{4}$} at 0.2 -3.5
\put{$\frac{1}{8}$} at 0.5 -2.8
\put{$\frac{1}{8}$} at -0.5 -2.8


\put{$\frac{5}{8}$} at 2.5 1.2
\put{$\frac{3}{8}$} at 1.5 2.2

\put{$\frac{5}{8}$} at -2.5 1.2
\put{$\frac{3}{8}$} at -1.5 2.2

\put{$\frac{5}{8}$} at 2.5 -0.8
\put{$\frac{3}{8}$} at 1.5 -1.8

\put{$\frac{5}{8}$} at -2.5 -0.8
\put{$\frac{3}{8}$} at -1.5 -1.8


\put{$\frac{5}{8}$} at 1.2 2.5
\put{$\frac{3}{8}$} at 2.2 1.5

\put{$\frac{5}{8}$} at -0.8 2.5
\put{$\frac{3}{8}$} at -1.8 1.5

\put{$\frac{5}{8}$} at 1.2 -2.5
\put{$\frac{3}{8}$} at 2.2 -1.5

\put{$\frac{5}{8}$} at -0.8 -2.5
\put{$\frac{3}{8}$} at -1.8 -1.5

\endpicture


\vspace{0.5cm}

Let us now consider the orbit of the interval $I^*$ defined by (\ref{def-interval}) for the
action of $N_{d-1}^*$. For simplicity, let us first deal with the case $d = 3$. With respect to
the generators $f_{2,1},f_{3,1},f_{3,2}$ of $N_2^*$, the orbit graph is depicted in Figure 2 below.
Here, $f_{2,1}$ corresponds to the generator whose action on the the graph is moving to the
right, whereas the action of both $f_{3,1}$ and $f_{3,2}$ consists in moving up, the former
by one unit and the latter with an amplitude that depends on the position. (Notice that the
directions of the arrows mean that we are only considering positive powers of the
generators.)

Now, the difficulty comes from that, as the reader may easily check, it is impossible to put
probability distributions on this graph yielding the equidistribution property along the
spheres centered at the origin. (This is already impossible
for the sphere of radius 4.) To overcome this problem, we will use the counting argument of
(the first part of) \cite{KN}, which actually corresponds to a deterministic counterpart of
the random walk argument above. Indeed, this argument is more robust in that it does not need
any equidistribution property, though it requires a certain extra argument to obtain our desired
infinite path as a concatenation of finite paths that behave nicely for certain finite processes.

\vspace{0.5cm}


\beginpicture

\setcoordinatesystem units <0.9cm,0.9cm>

\putrule from -4 -4 to -4 3
\putrule from -3 -4 to -3 3
\putrule from -2 -4 to -2 3
\putrule from -1 -4 to -1 3
\putrule from 0 -4 to 0 3
\putrule from 1 -4 to 1 3
\putrule from 2 -4 to 2 3
\putrule from 3 -4 to 3 3

\putrule from -4 -4 to 3 -4
\putrule from -4 -3 to 3 -3
\putrule from -4 -2 to 3 -2
\putrule from -4 -1 to 3 -1
\putrule from -4 0 to 3 0
\putrule from -4 1 to 3 1
\putrule from -4 2 to 3 2
\putrule from -4 3 to 3 3

\circulararc -52 degrees from -2 -4 center at 0 -3
\circulararc -52 degrees from -2 -3 center at 0 -2
\circulararc -52 degrees from -2 -2 center at 0 -1
\circulararc -52 degrees from -2 -1 center at 0 0
\circulararc -52 degrees from -2 0 center at 0 1
\circulararc -52 degrees from -2 1 center at 0 2

\circulararc -53 degrees from -1 -4 center at 2 -2.5
\circulararc -53 degrees from -1 -3 center at 2 -1.5
\circulararc -53 degrees from -1 -2 center at 2 -0.5
\circulararc -53 degrees from -1 -1 center at 2 0.5
\circulararc -53 degrees from -1 0 center at 2 1.5

\circulararc -48 degrees from 0 -4 center at 4.5 -2
\circulararc -48 degrees from 0 -3 center at 4.5 -1
\circulararc -48 degrees from 0 -2 center at 4.5 0
\circulararc -48 degrees from 0 -1 center at 4.5 1

\circulararc -53 degrees from 1 -4 center at 6 -1.5
\circulararc -53 degrees from 1 -3 center at 6 -0.5
\circulararc -53 degrees from 1 -2 center at 6 0.5

\circulararc -52 degrees from 2 -4 center at 8.2 -1
\circulararc -52 degrees from 2 -3 center at 8.2 0

\circulararc -48 degrees from 3 -4 center at 10.8 -0.5

\put{$\bullet$} at -4 -4
\put{} at -9.75 0

\small

\put{Figure 2} at -0.5 -4.35

\tiny

\plot -4 2.5 -4.05 2.4 /
\plot -4 2.5 -3.95 2.4 /
\plot -4 1.5 -4.05 1.4 /
\plot -4 1.5 -3.95 1.4 /
\plot -4 0.5 -4.05 0.4 /
\plot -4 0.5 -3.95 0.4 /
\plot -4 -0.5 -4.05 -0.6 /
\plot -4 -0.5 -3.95 -0.6 /
\plot -4 -1.5 -4.05 -1.6 /
\plot -4 -1.5 -3.95 -1.6 /
\plot -4 -2.5 -4.05 -2.6 /
\plot -4 -2.5 -3.95 -2.6 /
\plot -4 -3.5 -4.05 -3.6 /
\plot -4 -3.5 -3.95 -3.6 /

\plot -3 2.5 -3.05 2.4 /
\plot -3 2.5 -2.95 2.4 /
\plot -3 1.5 -3.05 1.4 /
\plot -3 1.5 -2.95 1.4 /
\plot -3 0.5 -3.05 0.4 /
\plot -3 0.5 -2.95 0.4 /
\plot -3 -0.5 -3.05 -0.6 /
\plot -3 -0.5 -2.95 -0.6 /
\plot -3 -1.5 -3.05 -1.6 /
\plot -3 -1.5 -2.95 -1.6 /
\plot -3 -2.5 -3.05 -2.6 /
\plot -3 -2.5 -2.95 -2.6 /
\plot -3 -3.5 -3.05 -3.6 /
\plot -3 -3.5 -2.95 -3.6 /

\plot -2 2.5 -2.05 2.4 /
\plot -2 2.5 -1.95 2.4 /
\plot -2 1.5 -2.05 1.4 /
\plot -2 1.5 -1.95 1.4 /
\plot -2 0.5 -2.05 0.4 /
\plot -2 0.5 -1.95 0.4 /
\plot -2 -0.5 -2.05 -0.6 /
\plot -2 -0.5 -1.95 -0.6 /
\plot -2 -1.5 -2.05 -1.6 /
\plot -2 -1.5 -1.95 -1.6 /
\plot -2 -2.5 -2.05 -2.6 /
\plot -2 -2.5 -1.95 -2.6 /
\plot -2 -3.5 -2.05 -3.6 /
\plot -2 -3.5 -1.95 -3.6 /

\plot -1 2.5 -1.05 2.4 /
\plot -1 2.5 -0.95 2.4 /
\plot -1 1.5 -1.05 1.4 /
\plot -1 1.5 -0.95 1.4 /
\plot -1 0.5 -1.05 0.4 /
\plot -1 0.5 -0.95 0.4 /
\plot -1 -0.5 -1.05 -0.6 /
\plot -1 -0.5 -0.95 -0.6 /
\plot -1 -1.5 -1.05 -1.6 /
\plot -1 -1.5 -0.95 -1.6 /
\plot -1 -2.5 -1.05 -2.6 /
\plot -1 -2.5 -0.95 -2.6 /
\plot -1 -3.5 -1.05 -3.6 /
\plot -1 -3.5 -0.95 -3.6 /

\plot 0 2.5 0.05 2.4 /
\plot 0 2.5 -0.05 2.4 /
\plot 0 1.5 0.05 1.4 /
\plot 0 1.5 -0.05 1.4 /
\plot 0 0.5 0.05 0.4 /
\plot 0 0.5 -0.05 0.4 /
\plot 0 -0.5 0.05 -0.6 /
\plot 0 -0.5 -0.05 -0.6 /
\plot 0 -1.5 0.05 -1.6 /
\plot 0 -1.5 -0.05 -1.6 /
\plot 0 -2.5 0.05 -2.6 /
\plot 0 -2.5 -0.05 -2.6 /
\plot 0 -3.5 0.05 -3.6 /
\plot 0 -3.5 -0.05 -3.6 /

\plot 0 2.5 0.05 2.4 /
\plot 0 2.5 -0.05 2.4 /
\plot 0 1.5 0.05 1.4 /
\plot 0 1.5 -0.05 1.4 /
\plot 0 0.5 0.05 0.4 /
\plot 0 0.5 -0.05 0.4 /
\plot 0 -0.5 0.05 -0.6 /
\plot 0 -0.5 -0.05 -0.6 /
\plot 0 -1.5 0.05 -1.6 /
\plot 0 -1.5 -0.05 -1.6 /
\plot 0 -2.5 0.05 -2.6 /
\plot 0 -2.5 -0.05 -2.6 /
\plot 0 -3.5 0.05 -3.6 /
\plot 0 -3.5 -0.05 -3.6 /

\plot 1 2.5 0.95 2.4 /
\plot 1 2.5 1.05 2.4 /
\plot 1 1.5 0.95 1.4 /
\plot 1 1.5 1.05 1.4 /
\plot 1 0.5 0.95 0.4 /
\plot 1 0.5 1.05 0.4 /
\plot 1 -0.5 0.95 -0.6 /
\plot 1 -0.5 1.05 -0.6 /
\plot 1 -1.5 0.95 -1.6 /
\plot 1 -1.5 1.05 -1.6 /
\plot 1 -2.5 0.95 -2.6 /
\plot 1 -2.5 1.05 -2.6 /
\plot 1 -3.5 0.95 -3.6 /
\plot 1 -3.5 1.05 -3.6 /

\plot 2 2.5 1.95 2.4 /
\plot 2 2.5 2.05 2.4 /
\plot 2 1.5 1.95 1.4 /
\plot 2 1.5 2.05 1.4 /
\plot 2 0.5 1.95 0.4 /
\plot 2 0.5 2.05 0.4 /
\plot 2 -0.5 1.95 -0.6 /
\plot 2 -0.5 2.05 -0.6 /
\plot 2 -1.5 1.95 -1.6 /
\plot 2 -1.5 2.05 -1.6 /
\plot 2 -2.5 1.95 -2.6 /
\plot 2 -2.5 2.05 -2.6 /
\plot 2 -3.5 1.95 -3.6 /
\plot 2 -3.5 2.05 -3.6 /

\plot 3 2.5 2.95 2.4 /
\plot 3 2.5 3.05 2.4 /
\plot 3 1.5 2.95 1.4 /
\plot 3 1.5 3.05 1.4 /
\plot 3 0.5 2.95 0.4 /
\plot 3 0.5 3.05 0.4 /
\plot 3 -0.5 2.95 -0.6 /
\plot 3 -0.5 3.05 -0.6 /
\plot 3 -1.5 2.95 -1.6 /
\plot 3 -1.5 3.05 -1.6 /
\plot 3 -2.5 2.95 -2.6 /
\plot 3 -2.5 3.05 -2.6 /
\plot 3 -3.5 2.95 -3.6 /
\plot 3 -3.5 3.05 -3.6 /


\plot 2.2 -4 2.1 -4.05 /
\plot 2.2 -4 2.1 -3.95 /
\plot 2.2 -3 2.1 -3.05 /
\plot 2.2 -3 2.1 -2.95 /
\plot 2.2 -2 2.1 -2.05 /
\plot 2.2 -2 2.1 -1.95 /
\plot 2.2 -1 2.1 -1.05 /
\plot 2.2 -1 2.1 -0.95 /
\plot 2.2 0 2.1 -0.05 /
\plot 2.2 0 2.1 0.05 /
\plot 2.2 1 2.1 1.05 /
\plot 2.2 1 2.1 0.95 /
\plot 2.2 2 2.1 2.05 /
\plot 2.2 2 2.1 1.95 /
\plot 2.2 3 2.1 3.05 /
\plot 2.2 3 2.1 2.95 /

\plot 1.2 -4 1.1 -4.05 /
\plot 1.2 -4 1.1 -3.95 /
\plot 1.2 -3 1.1 -3.05 /
\plot 1.2 -3 1.1 -2.95 /
\plot 1.2 -2 1.1 -2.05 /
\plot 1.2 -2 1.1 -1.95 /
\plot 1.2 -1 1.1 -1.05 /
\plot 1.2 -1 1.1 -0.95 /
\plot 1.2 0 1.1 -0.05 /
\plot 1.2 0 1.1 0.05 /
\plot 1.2 1 1.1 1.05 /
\plot 1.2 1 1.1 0.95 /
\plot 1.2 2 1.1 2.05 /
\plot 1.2 2 1.1 1.95 /
\plot 1.2 3 1.1 3.05 /
\plot 1.2 3 1.1 2.95 /

\plot 0.2 -4 0.1 -4.05 /
\plot 0.2 -4 0.1 -3.95 /
\plot 0.2 -3 0.1 -3.05 /
\plot 0.2 -3 0.1 -2.95 /
\plot 0.2 -2 0.1 -2.05 /
\plot 0.2 -2 0.1 -1.95 /
\plot 0.2 -1 0.1 -1.05 /
\plot 0.2 -1 0.1 -0.95 /
\plot 0.2 0 0.1 -0.05 /
\plot 0.2 0 0.1 0.05 /
\plot 0.2 1 0.1 1.05 /
\plot 0.2 1 0.1 0.95 /
\plot 0.2 2 0.1 2.05 /
\plot 0.2 2 0.1 1.95 /
\plot 0.2 3 0.1 3.05 /
\plot 0.2 3 0.1 2.95 /

\plot -0.8 -4 -0.9 -4.05 /
\plot -0.8 -4 -0.9 -3.95 /
\plot -0.8 -3 -0.9 -3.05 /
\plot -0.8 -3 -0.9 -2.95 /
\plot -0.8 -2 -0.9 -2.05 /
\plot -0.8 -2 -0.9 -1.95 /
\plot -0.8 -1 -0.9 -1.05 /
\plot -0.8 -1 -0.9 -0.95 /
\plot -0.8 0 -0.9 -0.05 /
\plot -0.8 0 -0.9 0.05 /
\plot -0.8 1 -0.9 1.05 /
\plot -0.8 1 -0.9 0.95 /
\plot -0.8 2 -0.9 2.05 /
\plot -0.8 2 -0.9 1.95 /
\plot -0.8 3 -0.9 3.05 /
\plot -0.8 3 -0.9 2.95 /

\plot -1.8 -4 -1.9 -4.05 /
\plot -1.8 -4 -1.9 -3.95 /
\plot -1.8 -3 -1.9 -3.05 /
\plot -1.8 -3 -1.9 -2.95 /
\plot -1.8 -2 -1.9 -2.05 /
\plot -1.8 -2 -1.9 -1.95 /
\plot -1.8 -1 -1.9 -1.05 /
\plot -1.8 -1 -1.9 -0.95 /
\plot -1.8 0 -1.9 -0.05 /
\plot -1.8 0 -1.9 0.05 /
\plot -1.8 1 -1.9 1.05 /
\plot -1.8 1 -1.9 0.95 /
\plot -1.8 2 -1.9 2.05 /
\plot -1.8 2 -1.9 1.95 /
\plot -1.8 3 -1.9 3.05 /
\plot -1.8 3 -1.9 2.95 /

\plot -2.8 -4 -2.9 -4.05 /
\plot -2.8 -4 -2.9 -3.95 /
\plot -2.8 -3 -2.9 -3.05 /
\plot -2.8 -3 -2.9 -2.95 /
\plot -2.8 -2 -2.9 -2.05 /
\plot -2.8 -2 -2.9 -1.95 /
\plot -2.8 -1 -2.9 -1.05 /
\plot -2.8 -1 -2.9 -0.95 /
\plot -2.8 0 -2.9 -0.05 /
\plot -2.8 0 -2.9 0.05 /
\plot -2.8 1 -2.9 1.05 /
\plot -2.8 1 -2.9 0.95 /
\plot -2.8 2 -2.9 2.05 /
\plot -2.8 2 -2.9 1.95 /
\plot -2.8 3 -2.9 3.05 /
\plot -2.8 3 -2.9 2.95 /

\plot -3.8 -4 -3.9 -4.05 /
\plot -3.8 -4 -3.9 -3.95 /
\plot -3.8 -3 -3.9 -3.05 /
\plot -3.8 -3 -3.9 -2.95 /
\plot -3.8 -2 -3.9 -2.05 /
\plot -3.8 -2 -3.9 -1.95 /
\plot -3.8 -1 -3.9 -1.05 /
\plot -3.8 -1 -3.9 -0.95 /
\plot -3.8 0 -3.9 -0.05 /
\plot -3.8 0 -3.9 0.05 /
\plot -3.8 1 -3.9 1.05 /
\plot -3.8 1 -3.9 0.95 /
\plot -3.8 2 -3.9 2.05 /
\plot -3.8 2 -3.9 1.95 /
\plot -3.8 3 -3.9 3.05 /
\plot -3.8 3 -3.9 2.95 /

\plot -2.085 -3.8 -1.985 -3.9 /
\plot -2.085 -3.8 -2.1 -3.93 /
\plot -2.085 -2.8 -1.985 -2.9 /
\plot -2.085 -2.8 -2.1 -2.93 /
\plot -2.085 -1.8 -1.985 -1.9 /
\plot -2.085 -1.8 -2.1 -1.93 /
\plot -2.085 -0.8 -1.985 -0.9 /
\plot -2.085 -0.8 -2.1 -0.93 /
\plot -2.085 0.2 -1.985 0.1 /
\plot -2.085 0.2 -2.1 0.07 /
\plot -2.085 1.2 -1.985 1.1 /
\plot -2.085 1.2 -2.1 1.07 /

\plot -1.085 -3.8 -0.985 -3.9 /
\plot -1.085 -3.8 -1.1 -3.93 /
\plot -1.085 -2.8 -0.985 -2.9 /
\plot -1.085 -2.8 -1.1 -2.93 /
\plot -1.085 -1.8 -0.985 -1.9 /
\plot -1.085 -1.8 -1.1 -1.93 /
\plot -1.085 -0.8 -0.985 -0.9 /
\plot -1.085 -0.8 -1.1 -0.93 /
\plot -1.085 0.2 -0.985 0.1 /
\plot -1.085 0.2 -1.1 0.07 /

\plot -0.085 -3.8 0.015 -3.9 /
\plot -0.085 -3.8 -0.1 -3.93 /
\plot -0.085 -2.8 0.015 -2.9 /
\plot -0.085 -2.8 -0.1 -2.93 /
\plot -0.085 -1.8 0.015 -1.9 /
\plot -0.085 -1.8 -0.1 -1.93 /
\plot -0.085 -0.8 0.015 -0.9 /
\plot -0.085 -0.8 -0.1 -0.93 /

\plot 0.91 -3.8 0.9 -3.91 /
\plot 0.91 -3.8 1 -3.88 /
\plot 0.91 -2.8 0.9 -2.91 /
\plot 0.91 -2.8 1 -2.88 /
\plot 0.91 -1.8 0.9 -1.91 /
\plot 0.91 -1.8 1 -1.88 /

\plot 1.91 -3.8 1.9 -3.91 /
\plot 1.91 -3.8 2 -3.88 /
\plot 1.91 -2.8 1.9 -2.91 /
\plot 1.91 -2.8 2 -2.88 /

\plot 2.91 -3.8 2.9 -3.91 /
\plot 2.91 -3.8 3 -3.88 /

\put{$\bullet$} at -4 -3
\put{$\bullet$} at -4 -2
\put{$\bullet$} at -4 -1
\put{$\bullet$} at -4 0
\put{$\bullet$} at -4 1
\put{$\bullet$} at -4 2
\put{$\bullet$} at -4 3

\put{$\bullet$} at -3 -4
\put{$\bullet$} at -2 -4
\put{$\bullet$} at -1 -4
\put{$\bullet$} at 0 -4
\put{$\bullet$} at 1 -4
\put{$\bullet$} at 2 -4
\put{$\bullet$} at 3 -4

\put{$\bullet$} at -3 -3
\put{$\bullet$} at -2 -3
\put{$\bullet$} at -1 -3
\put{$\bullet$} at 0 -3
\put{$\bullet$} at 1 -3
\put{$\bullet$} at 2 -3
\put{$\bullet$} at 3 -3

\put{$\bullet$} at -3 -2
\put{$\bullet$} at -2 -2
\put{$\bullet$} at -1 -2
\put{$\bullet$} at 0 -2
\put{$\bullet$} at 1 -2
\put{$\bullet$} at 2 -2
\put{$\bullet$} at 3 -2

\put{$\bullet$} at -3 -1
\put{$\bullet$} at -2 -1
\put{$\bullet$} at -1 -1
\put{$\bullet$} at 0 -1
\put{$\bullet$} at 1 -1
\put{$\bullet$} at 2 -1
\put{$\bullet$} at 3 -1

\put{$\bullet$} at -3 0
\put{$\bullet$} at -2 0
\put{$\bullet$} at -1 0
\put{$\bullet$} at 0 0
\put{$\bullet$} at 1 0
\put{$\bullet$} at 2 0
\put{$\bullet$} at 3 0

\put{$\bullet$} at -3 1
\put{$\bullet$} at -2 1
\put{$\bullet$} at -1 1
\put{$\bullet$} at 0 1
\put{$\bullet$} at 1 1
\put{$\bullet$} at 2 1
\put{$\bullet$} at 3 1

\put{$\bullet$} at -3 2
\put{$\bullet$} at -2 2
\put{$\bullet$} at -1 2
\put{$\bullet$} at 0 2
\put{$\bullet$} at 1 2
\put{$\bullet$} at 2 2
\put{$\bullet$} at 3 2

\put{$\bullet$} at -3 3
\put{$\bullet$} at -2 3
\put{$\bullet$} at -1 3
\put{$\bullet$} at 0 3
\put{$\bullet$} at 1 3
\put{$\bullet$} at 2 3
\put{$\bullet$} at 3 3

\endpicture


\vspace{0.5cm}

To close this section, let us finally explain why the exponent $\frac{2}{d(d-1)}$ is critical for
the action $\phi_{_{FF}}$. For simplicity, let us first consider the case $d=3$. Looking at the graph
of Figure 2 above, one easily computes the growth of the balls. This appears to be cubic, in the
sense that the number of points at distance \esp $\leq \! n$ \esp from the origin is \esp
$\frac{n^3 + 11n + 6}{6} \sim n^3$. These points correspond to intervals in the orbit of
$I^*$ obtained up to \esp $\leq \! n$ \esp compositions of the generators. Since these
intervals are disjoint, the length of a typical one should be of order \esp $\sim \! 1/n^3$.
\esp Hence, along a generic sequence of compositions, the value of the corresponding sum
$L_{\alpha}$ should be of order
$$\sum_{n \geq 1} \left( \frac{1}{n^{3}} \right)^{\alpha},$$
which is finite for $\alpha > \frac{1}{3} = \frac{2}{3(3-1)}$, as expected.

\vspace{0.1cm}

For the case of a general $d \geq 3$, it is very tempting trying to argue as before 
for any $\alpha$ larger than the inverse of the degree of growth of the graph of the orbit 
of $I$. Now, recall that according to the Bass-Guivarch formula (see \cite[Appendix]{gromov}), 
the growth of $N_{d-1}^*$ is polynomial of degree \esp $\sum_{i=1}^{d-1} i(d-i).$ \esp 
Moreover, the stabilizer of $I$ under the action of $N_{d-1}^*$ is the subgroup of 
$N_{d-1}^*$ made of the matrices whose first column is \esp $(1,0,0,\ldots,0)^T$. 
This subgroup naturally identifies with $N_{d-2}$, whose growth is polynomial of 
degree \esp $\sum_{i=1}^{d-2} i (d-i-1)$. \esp The difference of these degrees 
equals
\begin{equation}\label{growth-degree}
\sum_{i=1}^{d-1} i(d-i) - \sum_{i=1}^{d-1} i(d-i-1) = \sum_{i=1}^{d-1} i 
= \frac{d(d-1)}{2}.
\end{equation}
Since the graph of the orbit of $I$ identifies with the space of cosets 
$N_{d-1} / N_{d-2}$, one should expect that its growth is polynomial 
of degree given by (\ref{growth-degree}), and this is actually the case.


\subsection{Proof of Theorem A: the case $d = 3$}
\label{caso d=4}

\hspace{0.45cm} The proof of Theorem A is somewhat technical and requires hard notation.
This is the reason why we have chosen to first give the proof for the case $d = 3$, where
most of the ideas become more transparent and an important technical problem is overcomed
by a trick consisting in the introduction of a small parameter $\varepsilon > 0$. For the
general case, we use a slightly modified construction keeping essentially the same 
arguments. We begin with a lemma in the spirit of \cite[Lemma 2.2]{KN}.

\vspace{0.15cm}

\begin{lem} \label{densidades}
{\em Let $n \geq 1$ be an integer and let $C_1,C_2, \varepsilon$
be positive constants. Let $P$ be a set of \esp $\leq C_1 n^{3+\varepsilon}$ pairs
of non-negative integers $(i,j)$ associated to which there is a number $\ell_{i,j} > 0$
such that $\sum_{(i,j) \in P} \ell_{i,j} \leq 1$. Suppose that $P$ partitioned into
$n' \geq n^2 / C_2$ (resp. $n' \geq n^{2+\varepsilon} / C_2$) disjoint subsets
$P_1,\ldots,P_{n'}$. Then, given $A > 1$ and $1 > \alpha > 0$, the proportion
of indexes $m \in \{1,\ldots,n' \}$ for which}
$$\sum_{(i,j) \in P_m} \ell_{i,j}^\alpha \leq
\frac{A C_1^{1-\alpha} C_2}{n^{3\alpha - 1 - \varepsilon (1-\alpha)}} \qquad
\Big(\mbox{resp.} \!\!\!\! \quad \sum_{(i,j) \in P_m} \ell_{i,j}^\alpha \leq
\frac{A C_1^{1-\alpha} C_2}{n^{3\alpha - 1 + \varepsilon \alpha}} \esp \Big)$$
{\em is at least $1 - 1/A$.}
\end{lem}

\noindent{\bf Proof.} Since $\sum_{(i,j) \in P} \ell_{i,j} \leq 1$ and $P$ consists of at
most $C_1 n^{3+\varepsilon}$ pairs, a direct application of H\"older's inequality yields
$$\sum_{(i,j) \in P} \ell_{i,j}^{\alpha} \leq (C_1 n^{3+\varepsilon})^{1-\alpha}.$$
Hence,
$$\frac{1}{n'} \sum_{m=1}^{n'} \sum_{(i,j) \in P_m} \ell_{i,j}^{\alpha} 
\esp \leq \esp 
\frac{C_1^{1-\alpha} n^{(3+\varepsilon)(1-\alpha)}}{n'},$$
and the latter expression is less than or equal to \esp
$C_1^{1- \alpha} C_2 n^{1 - 3 \alpha + \varepsilon (1-\alpha)}$ \esp
(resp. \esp $C_1^{1- \alpha} C_2 n^{1 - 3\alpha - \varepsilon \alpha}$).
The lemma then follows as a direct application of Chebyshev's inequality. $\hfill\square$

\vspace{0.38cm}

Let us now come back to the graph associated to the action $\phi_{_{FF}}$ depicted in Figure 2,
and let us set $\ell_{i,j} := | f_{2,1}^i f_{3,1}^j (I^*)|$. Fix positive constants
$\alpha,\varepsilon$ such that
\begin{equation}\label{conditions}
\alpha > \frac{1}{3} = \frac{2}{(3-1)(3-2)}, \qquad
\varepsilon < \max \left\{ \frac{3 \alpha - 1}{1-\alpha}, 1 \right\}.
\end{equation}

For any real numbers $M \leq N$, we let $[[ M,N ]] := [M,N] \cap \mathbb{Z}$.
Given an integer $n \geq 2$, we consider the set $P(n) := [[ n,8n-1 ]] \times
[[ 0, n^{2+\varepsilon} ]]$. This set $P(n)$ consists of
$7n([n^{2+\varepsilon}] + 1) \leq 10 \esp n^{3+\varepsilon}$ points (with
$[\esp \cdot \esp]$ standing for the integer part), and is partitioned into
the $n' = [n^{2+\varepsilon}] + 1 \geq n^{2+\varepsilon}$ disjoint sets
(horizontal paths) \esp $P(n,1), P(n,2), \ldots, P(n,n')$ \esp given by
$$P(n,m) := \big\{ (n,m), (n+1,m), \ldots, (8n-1,m) \big\}.$$
By the preceding lemma, for each $0 < A_n < 1$,
the proportion of indexes $m \in \{1,\ldots,n'\}$ for which
\begin{equation}\label{first}
\sum_{i=n}^{8n -1} \ell_{i,m}^{\alpha} \esp\esp\esp = \!\!\sum_{(i,j) \in P(n,m)} \ell_{i,j}^{\alpha}
\esp\esp\esp \leq \esp\esp\esp \frac{A_n 10^{1-\alpha}}{n^{3 \alpha - 1 + \varepsilon \alpha}}
\end{equation}
is at least $1 - 1/A_n$. Notice that each path $P(n,m)$
comes from the action of the generator $f_{2,1}$.

\vspace{0.08cm}

Similarly, for each integer $n \geq 2$, let us consider the set
$Q(n) := [[ n, 2n -1 ]] \times [[ 0, n^{2+\varepsilon}]] $ consisting of \esp
$n ([n^{2+\varepsilon}] + 1) \leq 2 n^{3 + \varepsilon}$ \esp points. Although in general
there is no partition of $Q(n)$ into paths induced by the action of $f_{3,1}, f_{3,2}$
all of them having the same number of points, a partition that almost satisfies this
property (and that will be sufficient for our purposes) can be defined as follows.
For each \esp $n \leq m \leq 2n-1$ \esp we divide the set
$\{(m,0), (m,1), \ldots \}$ into $n$ paths via the following rules:

\vspace{0.15cm}

\noindent -- For each $0 \leq j \leq n-2$, there is a path starting at $(m,j)$ jumping upwards
of $m$ units;

\vspace{0.15cm}

\noindent -- The path starting at $(m,n-1)$ makes $m-n$ jumps upwards of 1 unit
and then makes a jump of $m$ units;

\vspace{0.15cm}

\noindent -- The picture repeats ``periodically'', so that each infinite path is
made of $n-1$ consecutive jumps of $m$ units followed by $m-n$ jumps of 1 unit.

\vspace{0.15cm}

Figure 3 illustrates the case where $n = 3$ and $m = 5$ though the resulting paths are disposed
horizontally instead of vertically by obvious reasons. Although one may give precise formulas for
the points in each of these paths, this is not completely necessary. The main property that we
will retain is the obvious fact that the number of points of each of them inside any rectangle
$[[ n, 2n-1 ]] \times [[ 0, K - 1 ]]$ lies between \esp $\frac{K}{n} - 2n$ \esp and \esp
$\frac{K}{n} + 2n$. \esp (An alternative construction leading to a much better -logarithmic-
control of the deviation
will be given in \S \ref{section-general}.) In particular, we have an induced partition
of $Q(n)$ into $n'' = n^2$ paths \esp $Q(n,1), Q(n,2), \ldots, Q(n,n'')$ \esp for
which the preceding lemma yields that for each $A_n > 0$, the
proportion of indexes $m \in \{1,\ldots,n''\}$ satisfying
\begin{equation}\label{second}
\sum_{(i,j) \in Q(n,m)} \ell_{i,j}^{\alpha} \leq
\frac{A_n 2^{1-\alpha}}{n^{3 \alpha - 1 - \varepsilon (1-\alpha)}}
\end{equation}
is at least $1 - 1/A_n$. Notice again that each of these paths comes from the action of
the generators $f_{3,1}$ and $f_{3,2}$ according to the amplitude of the jump.

\vspace{0.84cm}


\beginpicture

\setcoordinatesystem units <0.418cm,0.418cm>

\circulararc 80 degrees from 1 1 center at 3.5 4
\circulararc 80 degrees from 2 1 center at 4.5 4
\circulararc 80 degrees from 5 1 center at 7.5 4
\circulararc 80 degrees from 6 1 center at 8.5 4
\circulararc 80 degrees from 9 1 center at 11.5 4
\circulararc 80 degrees from 10 1 center at 12.5 4
\circulararc 80 degrees from 13 1 center at 15.5 4
\circulararc 80 degrees from 14 1 center at 16.5 4
\circulararc 80 degrees from 17 1 center at 19.5 4
\circulararc 80 degrees from 18 1 center at 20.5 4
\circulararc 80 degrees from 21 1 center at 23.5 4
\circulararc 80 degrees from 22 1 center at 24.5 4
\circulararc 80 degrees from 25 1 center at 27.5 4
\circulararc 80 degrees from 26 1 center at 28.5 4
\circulararc 80 degrees from 29 1 center at 31.5 4
\circulararc 80 degrees from 30 1 center at 32.5 4
\circulararc 80 degrees from 33 1 center at 35.5 4
\circulararc 80 degrees from 34 1 center at 36.5 4

\circulararc 180 degrees from 3 1 center at 3.5 1
\circulararc 180 degrees from 4 1 center at 4.5 1
\circulararc 180 degrees from 7 1 center at 7.5 1
\circulararc 180 degrees from 8 1 center at 8.5 1
\circulararc 180 degrees from 11 1 center at 11.5 1
\circulararc 180 degrees from 12 1 center at 12.5 1
\circulararc 180 degrees from 15 1 center at 15.5 1
\circulararc 180 degrees from 16 1 center at 16.5 1
\circulararc 180 degrees from 19 1 center at 19.5 1
\circulararc 180 degrees from 20 1 center at 20.5 1
\circulararc 180 degrees from 23 1 center at 23.5 1
\circulararc 180 degrees from 24 1 center at 24.5 1
\circulararc 180 degrees from 27 1 center at 27.5 1
\circulararc 180 degrees from 28 1 center at 28.5 1
\circulararc 180 degrees from 31 1 center at 31.5 1
\circulararc 180 degrees from 32 1 center at 32.5 1
\circulararc 180 degrees from 35 1 center at 35.5 1
\circulararc 180 degrees from 36 1 center at 36.5 1
\circulararc 180 degrees from 39 1 center at 39.5 1

\put{} at 0.75 0

\small

\put{Figure 3} at 20.6 -0.8

\tiny

\put{$\bullet$} at 40 1
\put{$\bullet$} at 39 1
\put{$\bullet$} at 38 1
\put{$\bullet$} at 37 1
\put{$\bullet$} at 36 1
\put{$\bullet$} at 35 1
\put{$\bullet$} at 34 1
\put{$\bullet$} at 33 1
\put{$\bullet$} at 32 1
\put{$\bullet$} at 31 1
\put{$\bullet$} at 30 1
\put{$\bullet$} at 29 1
\put{$\bullet$} at 28 1
\put{$\bullet$} at 27 1
\put{$\bullet$} at 26 1
\put{$\bullet$} at 25 1
\put{$\bullet$} at 24 1
\put{$\bullet$} at 23 1
\put{$\bullet$} at 22 1
\put{$\bullet$} at 21 1
\put{$\bullet$} at 20 1
\put{$\bullet$} at 19 1
\put{$\bullet$} at 18 1
\put{$\bullet$} at 17 1
\put{$\bullet$} at 16 1
\put{$\bullet$} at 15 1
\put{$\bullet$} at 14 1
\put{$\bullet$} at 13 1
\put{$\bullet$} at 12 1
\put{$\bullet$} at 11 1
\put{$\bullet$} at 1 1
\put{$\bullet$} at 2 1
\put{$\bullet$} at 3 1
\put{$\bullet$} at 4 1
\put{$\bullet$} at 5 1
\put{$\bullet$} at 6 1
\put{$\bullet$} at 7 1
\put{$\bullet$} at 8 1
\put{$\bullet$} at 9 1
\put{$\bullet$} at 10 1

\endpicture


\vspace{0.75cm}

We will apply the preceding construction for each integer \esp $n = n_k := 4^k$, \esp
where $k \geq 1$. The choice of the constants $A_{n_k}$ is as follows. First, we let
$r_k$ (resp. $s_k$) be the minimum (resp. maximum) number of points of a path of the
form $Q(n_k,m)$ inside $Q(n_k)$. Similarly, we let $r_k'$ (resp. $s_k'$) be the
minimum (resp. maximum) number of points in a path of the form $Q(n_k,m)$ inside
$P(n_{k-1}) \cap Q(n_k)$. Finally, we let
\begin{equation}\label{cte-B}
B := \prod_{k \geq 2} \frac{s_k}{r_k} \frac{s_k'}{r_k'}.
\end{equation}
Notice that the value of $B$ is finite. Indeed, by the discussion above, we have
$$4^{k(1+\varepsilon)} - 2 \cdot 4^{k} =
n_k^{1+\varepsilon} - 2n_k \leq r_k \leq s_k \leq n_k^{1+\varepsilon} + 2n_k
=4^{k(1+\varepsilon)} + 2 \cdot 4^{k}$$
and
$$4^{k + k \varepsilon -1} - 2 \cdot 4^{k+1}
= \frac{n_k^{2+\varepsilon}}{n_{k+1}} - 2n_{k+1} \leq r_k'
\leq s_k' \leq \frac{n_k^{2+\varepsilon}}{n_{k+1}} + 2n_{k+1}
= 4^{k + k \varepsilon -1} + 2 \cdot 4^{k+1},$$
which easily yield the convergence of the infinite product in the definition of $B$.
We will also use the constant
\begin{equation}\label{ctes-C-D}
C := 4 \sum_{k \geq 1} \frac{1}{2^{k(3\alpha - 1 -\varepsilon(1-\alpha))}}.
\end{equation}
Notice again that since (\ref{conditions}) implies that \esp
$3 \alpha - 1 - \varepsilon (1-\alpha) > 0$, we have $C < \infty$.

\vspace{0.1cm}

We now fix \esp $A_{n_1} \geq 2^{2+k(3\alpha - 1 - \varepsilon(1-\alpha))} BC$ \esp such
that (\ref{first}) holds for $n = n_1$ and {\em every} $m$ in the corresponding range.
Finally, for $k \geq 2$, we set
$$A_{n_k} := B \esp C \esp 2^{k(3\alpha - 1 - \varepsilon(1-\alpha))}.$$

We next state a key lemma whose proof is postponed in order to proceed
immediately to the proof of Theorem~A in the case $d = 3$.

\vspace{0.25cm}

\begin{lem} \label{fin}
{\em There are two infinite sequences of paths $P(n_k,m'_k)$ and $Q(n_k,m_k'')$ such that
{\em (\ref{first})} (resp. {\em (\ref{second}))} holds for $n = n_k$ and $m = m'_k$
(resp. $m = m''_k$) and such that $P(n_k,m'_k)$ intersects both $Q(n_k,m''_k)$ and
$Q(n_{k+1},m''_{k+1})$ for all $k \geq 1$.}
\end{lem}

\vspace{0.4cm}


\beginpicture

\setcoordinatesystem units <0.77cm,0.77cm>

\putrule from -3 -4 to -3 -2
\putrule from -2 -4 to -2 -2
\putrule from 0 -4 to 0 3
\putrule from 3 -4 to 3 3
\putrule from 9 -4 to 9 4
\putrule from -3 -4 to 4 -4
\putrule from -3 -2 to 4 -2
\putrule from 0 3 to 4 3
\putrule from 3 -4 to 15 -4
\putrule from 3 3 to 15 3
\putrule from 15 -4 to 15 4

\putrule from -2.6 -3.6 to -2.6 -2.4
\putrule from -2.612 -3.6 to -2.612 -2.4
\putrule from -2.582 -3.6 to -2.582 -2.4

\putrule from -2.6 -2.4 to 2 -2.4
\putrule from -2.6 -2.412 to 2 -2.412
\putrule from -2.6 -2.388 to 2 -2.388

\putrule from 2 -2.4 to 2 -3.6
\putrule from 2.012 -2.4 to 2.012 -3.6
\putrule from 1.988 -2.4 to 1.988 -3.6

\putrule from 2 -3.6 to 13 -3.6
\putrule from 2 -3.612 to 13 -3.612
\putrule from 2 -3.588 to 13 -3.588

\putrule from 13 -3.6 to 13 4
\putrule from 13.012 -3.6 to 13.012 4
\putrule from 12.988 -3.6 to 12.988 4

\putrule from -5.6 -3.6 to -2.6 -3.6
\putrule from -5.6 -3.612 to -2.6 -3.612
\putrule from -5.6 -3.588 to -2.6 -3.588

\put{} at -4 0

\small

\put{Figure 4} at 5.05 -5

\put{0} at -4.8 -4
\put{$[(4^{k})^{2+\varepsilon}]$} at -4.8 -2
\put{$[(4^{k+1})^{2+\varepsilon}]$} at -4.8 3

\tiny

\put{$Q(n_k,m_k'')$} at -2.5 -1.7
\put{$P(n_k,m_k')$} at -1 -2.75
\put{$Q(n_{k+1},m_{k+1}'')$} at 2 3.3
\put{$P(n_{k+1},m_{k+1}')$} at 6.5 -3.3
\put{$Q(n_{k+2},m''_{k+2})$} at 13 4.32
\put{$P(n_{k-1},m_{k-1}')$} at -4.3 -3.3

\put{$\bullet$} at -2 -3.6
\put{$\bullet$} at -2.2 -3.6
\put{$\bullet$} at -2.4 -3.6
\put{$\bullet$} at -2.8 -3.6
\put{$\bullet$} at -3 -3.6
\put{$\bullet$} at -3.2 -3.6
\put{$\bullet$} at -3.4 -3.6
\put{$\bullet$} at -3.6 -3.6
\put{$\bullet$} at -3.8 -3.6
\put{$\bullet$} at -4 -3.6
\put{$\bullet$} at -4.2 -3.6
\put{$\bullet$} at -4.4 -3.6
\put{$\bullet$} at -4.6 -3.6
\put{$\bullet$} at -4.8 -3.6
\put{$\bullet$} at -5 -3.6
\put{$\bullet$} at -5.2 -3.6
\put{$\bullet$} at -5.4 -3.6
\put{$\bullet$} at -5.6 -3.6

\put{$\bullet$} at -3 -2.4
\put{$\bullet$} at -2.8 -2.4
\put{$\bullet$} at -2.6 -2.4
\put{$\bullet$} at -2.4 -2.4
\put{$\bullet$} at -2.2 -2.4
\put{$\bullet$} at -2 -2.4
\put{$\bullet$} at -1.8 -2.4
\put{$\bullet$} at -1.6 -2.4
\put{$\bullet$} at -1.4 -2.4
\put{$\bullet$} at -1.2 -2.4
\put{$\bullet$} at -1 -2.4
\put{$\bullet$} at -0.8 -2.4
\put{$\bullet$} at -0.6 -2.4
\put{$\bullet$} at -0.4 -2.4
\put{$\bullet$} at -0.2 -2.4
\put{$\bullet$} at 0 -2.4
\put{$\bullet$} at 0.2 -2.4
\put{$\bullet$} at 0.4 -2.4
\put{$\bullet$} at 0.6 -2.4
\put{$\bullet$} at 0.8 -2.4
\put{$\bullet$} at 1 -2.4
\put{$\bullet$} at 1.2 -2.4
\put{$\bullet$} at 1.4 -2.4
\put{$\bullet$} at 1.6 -2.4
\put{$\bullet$} at 1.8 -2.4
\put{$\bullet$} at 2 -2.4
\put{$\bullet$} at 2.2 -2.4
\put{$\bullet$} at 2.4 -2.4
\put{$\bullet$} at 2.6 -2.4
\put{$\bullet$} at 2.8 -2.4
\put{$\bullet$} at 3 -2.4

\put{$\bullet$} at -2.6 -2
\put{$\bullet$} at -2.6 -2.4
\put{$\bullet$} at -2.6 -2.8
\put{$\bullet$} at -2.6 -3.2
\put{$\bullet$} at -2.6 -3.6
\put{$\bullet$} at -2.6 -4
\put{$\bullet$} at -2.6 -3.4

\put{$\bullet$} at 2 3
\put{$\bullet$} at 2 2.4
\put{$\bullet$} at 2 1.8
\put{$\bullet$} at 2 1.2
\put{$\bullet$} at 2 0.6
\put{$\bullet$} at 2 0
\put{$\bullet$} at 2 -0.6
\put{$\bullet$} at 2 -1.2
\put{$\bullet$} at 2 -1.8
\put{$\bullet$} at 2 -2.4
\put{$\bullet$} at 2 -3
\put{$\bullet$} at 2 -3.6

\put{$\bullet$} at 2 -2.6
\put{$\bullet$} at 2 -2.8
\put{$\bullet$} at 2 0.8
\put{$\bullet$} at 2 1

\put{$\bullet$} at 0 -3.6
\put{$\bullet$} at 0.2 -3.6
\put{$\bullet$} at 0.4 -3.6
\put{$\bullet$} at 0.6 -3.6
\put{$\bullet$} at 0.8 -3.6
\put{$\bullet$} at 1 -3.6
\put{$\bullet$} at 1.2 -3.6
\put{$\bullet$} at 1.4 -3.6
\put{$\bullet$} at 1.6 -3.6
\put{$\bullet$} at 1.8 -3.6
\put{$\bullet$} at 2 -3.6
\put{$\bullet$} at 2.2 -3.6
\put{$\bullet$} at 2.4 -3.6
\put{$\bullet$} at 2.6 -3.6
\put{$\bullet$} at 2.8 -3.6
\put{$\bullet$} at 3 -3.6
\put{$\bullet$} at 3.2 -3.6
\put{$\bullet$} at 3.4 -3.6
\put{$\bullet$} at 3.6 -3.6
\put{$\bullet$} at 3.8 -3.6
\put{$\bullet$} at 4 -3.6
\put{$\bullet$} at 4.2 -3.6
\put{$\bullet$} at 4.4 -3.6
\put{$\bullet$} at 4.6 -3.6
\put{$\bullet$} at 2.8 -3.6
\put{$\bullet$} at 3 -3.6
\put{$\bullet$} at 3.2 -3.6
\put{$\bullet$} at 3.4 -3.6
\put{$\bullet$} at 3.6 -3.6
\put{$\bullet$} at 3.8 -3.6
\put{$\bullet$} at 4 -3.6
\put{$\bullet$} at 4.2 -3.6
\put{$\bullet$} at 4.4 -3.6
\put{$\bullet$} at 4.6 -3.6
\put{$\bullet$} at 2.8 -3.6
\put{$\bullet$} at 3 -3.6
\put{$\bullet$} at 3.2 -3.6
\put{$\bullet$} at 3.4 -3.6
\put{$\bullet$} at 3.6 -3.6
\put{$\bullet$} at 3.8 -3.6
\put{$\bullet$} at 4 -3.6
\put{$\bullet$} at 4.2 -3.6
\put{$\bullet$} at 4.4 -3.6
\put{$\bullet$} at 4.6 -3.6
\put{$\bullet$} at 3.8 -3.6
\put{$\bullet$} at 4 -3.6
\put{$\bullet$} at 4.2 -3.6
\put{$\bullet$} at 4.4 -3.6
\put{$\bullet$} at 4.6 -3.6
\put{$\bullet$} at 2.8 -3.6
\put{$\bullet$} at 3 -3.6
\put{$\bullet$} at 3.2 -3.6
\put{$\bullet$} at 3.4 -3.6
\put{$\bullet$} at 3.6 -3.6
\put{$\bullet$} at 3.8 -3.6
\put{$\bullet$} at 4 -3.6
\put{$\bullet$} at 4.2 -3.6
\put{$\bullet$} at 4.4 -3.6
\put{$\bullet$} at 4.6 -3.6
\put{$\bullet$} at 4.8 -3.6
\put{$\bullet$} at 5 -3.6
\put{$\bullet$} at 5.2 -3.6
\put{$\bullet$} at 5.4 -3.6
\put{$\bullet$} at 5.6 -3.6
\put{$\bullet$} at 5.8 -3.6
\put{$\bullet$} at 6 -3.6
\put{$\bullet$} at 6.2 -3.6
\put{$\bullet$} at 6.4 -3.6
\put{$\bullet$} at 6.6 -3.6
\put{$\bullet$} at 6.8 -3.6
\put{$\bullet$} at 7 -3.6
\put{$\bullet$} at 7.2 -3.6
\put{$\bullet$} at 7.4 -3.6
\put{$\bullet$} at 7.6 -3.6
\put{$\bullet$} at 7.8 -3.6
\put{$\bullet$} at 8 -3.6
\put{$\bullet$} at 8.2 -3.6
\put{$\bullet$} at 8.4 -3.6
\put{$\bullet$} at 8.6 -3.6
\put{$\bullet$} at 8.8 -3.6
\put{$\bullet$} at 9 -3.6
\put{$\bullet$} at 9.2 -3.6
\put{$\bullet$} at 9.4 -3.6
\put{$\bullet$} at 9.6 -3.6
\put{$\bullet$} at 9.8 -3.6
\put{$\bullet$} at 10 -3.6
\put{$\bullet$} at 10.2 -3.6
\put{$\bullet$} at 10.4 -3.6
\put{$\bullet$} at 10.6 -3.6
\put{$\bullet$} at 10.8 -3.6
\put{$\bullet$} at 11 -3.6
\put{$\bullet$} at 11.2 -3.6
\put{$\bullet$} at 11.4 -3.6
\put{$\bullet$} at 11.6 -3.6
\put{$\bullet$} at 11.8 -3.6
\put{$\bullet$} at 12 -3.6
\put{$\bullet$} at 12.2 -3.6
\put{$\bullet$} at 12.4 -3.6
\put{$\bullet$} at 12.6 -3.6
\put{$\bullet$} at 12.8 -3.6
\put{$\bullet$} at 13 -3.6
\put{$\bullet$} at 13.2 -3.6
\put{$\bullet$} at 13.4 -3.6
\put{$\bullet$} at 13.6 -3.6
\put{$\bullet$} at 13.8 -3.6
\put{$\bullet$} at 14 -3.6
\put{$\bullet$} at 14.2 -3.6
\put{$\bullet$} at 14.4 -3.6
\put{$\bullet$} at 14.6 -3.6
\put{$\bullet$} at 14.8 -3.6
\put{$\bullet$} at 15 -3.6

\put{$\bullet$} at 13 -3.4
\put{$\bullet$} at 13 -2.6
\put{$\bullet$} at 13 -1.8
\put{$\bullet$} at 13 -1
\put{$\bullet$} at 13 -0.2
\put{$\bullet$} at 13 0.6
\put{$\bullet$} at 13 1.4
\put{$\bullet$} at 13 2.2
\put{$\bullet$} at 13 3
\put{$\bullet$} at 13 3.8

\put{$\bullet$} at 13 3.2
\put{$\bullet$} at 13 3.4
\put{$\bullet$} at 13 3.6

\put{$\bullet$} at 13 -2
\put{$\bullet$} at 13 -2.2
\put{$\bullet$} at 13 -2.4

\setdots

\putrule from -4.5 -4 to -3 -4
\putrule from -3.7 -2 to -3 -2
\putrule from -3.5 3 to 0 3

\put{$4^{k}$} at -3 -4.25
\put{$4^{k+2}$} at 9 -4.25
\put{$2 \! \cdot \! 4^{k} - 1$} at -2 -4.5
\put{$4^{k+1}$} at 0 -4.25
\put{$2 \! \cdot \! 4^{k+1} - 1$} at 3 -4.5
\put{$2 \! \cdot \! 4^{k+2} - 1$} at 14.8 -4.5

\endpicture


\vspace{0.7cm}

Assuming this lemma, the proof of Theorem~A in the case $d = 3$ is at hand. Indeed, the concatenation
of the sequence of finite paths provided by the lemma naturally yields an infinite path without
loops which is in correspondence with a sequence of compositions by
$f_{2,1},f_{3,1},f_{3,1}^{-1},f_{3,2},f_{3,2}^{-1}$ (see
Figure 4). By construction, for this sequence of iterations, the value of the corresponding
$L_{\alpha}$-sum (\ref{clave}) for the interval $I^{**}$ corresponding to the initial point
of $Q(n_1,m_1')$ is less than or equal to
\begin{eqnarray*}
10^{1 - \alpha} \sum_{k \geq 1} \frac{A_{n_k}}{n_{k}^{3\alpha - 1 + \varepsilon \alpha}} +
2^{1-\alpha} \sum_{k \geq 1} \frac{A_{n_k}}{n_k^{3 \alpha - 1 - \varepsilon (1-\alpha)}}
&\leq&
\frac{20 A_{n_1}}{4^{3\alpha - 1 - \varepsilon(1-\alpha)}} +
\sum_{k \geq 2} \frac{40 BC}{2^{k(3 \alpha - 1 - \varepsilon (1-\alpha))}}\\
&\leq& 80 A_{n_1} 4^{\varepsilon(1-\alpha)} + 40 B C^2.
\end{eqnarray*}
This interval $I^{**}$ is in the orbit of $I^*$, from which it can
be reached in no more than \esp $(2 \cdot 4^{1} - 1) + 4 = 11$
\esp iterations of the generator $f_{2,1}$. By concatenating this
finite path to the previous one, we obtain an infinite path
associated to which the $L_{\alpha}$-sum corresponding to $I^*$ is
finite, which allows to conclude the proof by the arguments
developed in \S \ref{main-idea}.

\vspace{0.1cm}

All that remains for completing the proof of Theorem A in the case $d = 3$ is the

\vspace{0.35cm}

\noindent{\bf Proof of Lemma \ref{fin}.} The argument is similar to that of \cite[Lemma 2.3]{KN},
but it needs a slight modification. Namely, for each $k \geq 1$, we let $D_k'$ be the density of
indexes $m' \in \{1, \ldots, [n_k^{2+\varepsilon}] + 1\}$ such that $P(n_k,m')$ is ``reached''
by a sequence of paths $Q(n_1,m_1''), P(n_1,m_1'), \ldots, Q(n_k,m_k'')$ satisfying:

\vspace{0.1cm}

\noindent -- $P(n_i,m_i')$ intersects both $Q(n_i,m_i'')$ and $Q(n_{i+1},m_{i+1}'')$ for all
$1 \leq i \leq k-1$, whereas $P(n_k,m')$ intersects $Q(n_k,m_k'')$;

\vspace{0.1cm}

\noindent -- Inequality (\ref{first}) (resp. (\ref{second})) holds for $(n,m) = (n_i,m_i')$
whenever $1 \leq i \leq k-1$ as well as for $(n,m) = (n_k,m')$ (resp. for
$(n,m) = (n_i,m_i'')$ whenever $1 \leq i \leq k$).

\vspace{0.15cm}

\noindent Similarly, we denote by $D_k''$ the density of indexes
$m'' \in \{1, \ldots, n_k^2\}$ such that $Q(n_k,m'')$ is reached by a sequence
of paths $Q(n_1,m_1''), P(n_1,m_1'), \ldots, P(n_{k-1},m_{k-1}')$ satisfying:

\vspace{0.1cm}

\noindent -- $P(n_i,m_i')$ intersects both $Q(n_i,m_i'')$ and $Q(n_{i+1},m_{i+1}'')$ for all
$1 \leq i \leq k-1$;

\vspace{0.1cm}

\noindent -- Inequality (\ref{first}) (resp. (\ref{second})) holds for $(n,m) = (n_i,m_i')$
(resp. for $(n,m) = (n_i,m_i'')$) whenever $1 \leq i \leq k-1$ as well as for $(n,m) = (n_k,m'')$.

\vspace{0.1cm}

We claim that the following relations hold:
\begin{equation}\label{key}
1- D_k' \leq (1-D_k'') \frac{s_k}{r_k} + \frac{1}{A_{n_k}}, \esp \qquad \esp
1 - D_{k+1}'' \leq (1-D_k') \frac{s_{k+1}'}{r_{k+1}'} + \frac{1}{A_{n_{k+1}}}.
\end{equation}
Assuming this for a while, we obtain for each $k \geq 1$,
$$1 - D_{k}' \leq (1-D_{k-1}') \frac{s_k}{r_k}
\frac{s_k'}{r_k'} + \frac{1}{A_{n_k}} \frac{s_k}{r_k} + \frac{1}{A_{n_k}}.$$
Using induction, this easily yields
$$1-D_k' \leq (1-D_1') \prod_{i=2}^{k} \frac{s_i}{r_i} \frac{s_i'}{r_i'}
+ 2 \sum_{i=2}^{k} \frac{1}{A_{n_i}} \prod_{j=2}^{i} \frac{s_j}{r_j}.$$
From the definition $n_i := 4^{i}$ and that of the constant
$B$ in (\ref{cte-B}), one concludes that for each $k \geq 1$,
$$1 - D_k' \leq (1-D_1') B + 2B \sum_{i = 1}^{k} \frac{1}{A_{n_i}}.$$
Now, the choice of $A_{n_1}$ was made so that $D_1' = 1$, hence
$$1 - D_k' \leq 2B \sum_{i \geq 1} \frac{1}{A_{n_i}} \leq \frac{1}{2}.$$
Thus, $D_k' \geq 1/2$ holds for all $k \geq 1$, which provides finite paths
satisfying the desired properties of length as large as we want. The infinite path
claimed to exist is obtained easily from this by means of a Cantor diagonal type argument.

Finally, it remains to show (\ref{key}). The proof follows the same principle of that
of \cite[Lemma 2.3]{KN} but requires a little adjustment. First, we denote by
$\hat{D}_k''$ the density of points in $Q(n_k)$ that are ``well-attainable''
in the sense that they belong to the last of a sequence of consecutively
intersecting paths $Q(n_1,m_1''), P(n_1,m_1'), \ldots, P(n_{k-1},m_{k-1}')$,
$\!Q(n_k,m_k'')$ for which inequalities of type (\ref{first}) or
(\ref{second}) hold according to the case. We have
\begin{equation}\label{small-key}
(1 - D_k') \leq (1 - \hat{D}_k'') + \frac{1}{A_{n_k}}.
\end{equation}
Indeed, the term $1/A_{n_k}$ corresponds to the density of horizontal paths in $P(n_k)$ that are
``bad by themselves'' in the sense that the corresponding type (\ref{first}) inequality does
not hold for them. The term $(1-\hat{D}_k'')$ corresponds to the density of paths in $P(n_k)$
that may be good by themselves but intersect $Q(n_k)$ at a set formed only by non-well-attainable
points. (Notice that we are using the fact that all horizontal paths in $P(n_k)$ have the same
number of points in $Q(n_k)$.) The left-side inequality in (\ref{key}) then follows as a
combination of (\ref{small-key}) and the inequality
$$1 - \hat{D}_k'' \leq
(1 - D_k'') \frac{s_k}{r_k},$$
where the correction factor comes from the fact that although the number of points in
each path of the form $Q(n_k,m)$ is not constant, it varies between $r_k$ and $s_k$.

Similarly, in the right-side inequality, the term $1/A_{n_{k+1}}$ corresponds to the density
of bad-by-themselves  paths of the form $Q(n_{k+1},m)$ in $Q(n_{k+1})$. The term $(1-D_k')$
corresponds to the ``accumulated density of bad paths'' up to $P(n_k)$, and equals the density
of ``non-well-attainable'' points in $P(n_k) \cap Q (n_{k+1})$. Finally, the correction factor
comes from the fact that the number of points in $P(n_k) \cap Q(n_{k+1})$ contained in each
path of the form $Q(n_{k+1},m)$ lies between $r_{k+1}'$ and $s_{k+1}'$.


\subsection{Proof of Theorem A: the general case}
\label{section-general}

\hspace{0.45cm} To deal with the general case we will follow a similar strategy, though most of
the computations become more involved. We will now consider paths inside parallelepipeds of
dimension $d-1$ having sides of length of (relative) order $k,k^2,\ldots,k^{d-1}$. This will
make naturally appear the exponent $\frac{d(d-1)}{2}$ in relation to the total number of
points in the parallelepiped. The most relevant difficulty will be related to the decomposition
into paths. Indeed, the construction of the preceding section illustrated by Figure 3 is no
longer satisfactory, and we will need to carry out a nontrivial modification of it. Since
this is of independent interest and has potential applications in other contexts, the discussion
of the new construction will be the subject of \S \ref{section-comb}. Here we content ourselves
in stating what we need for our purposes, which is summarized in the next

\vspace{0.15cm}

\begin{lem} \label{decompos} {\em Let $M > N$ be positive integer numbers, with $N$ of the form
$1+2^k$. There exists a decomposition of $\mathbb{N}_0 := \{0,1,\ldots\}$ into $N$
subsets (paths) satisfying:}

\vspace{0.1cm}

\noindent (i) {\em The distance (jump) between two 
consecutive points of each path is either $M$ or 1;}

\vspace{0.1cm}

\noindent (ii) {\em For all $0 \leq K_1 < K_2$, the maximal number
of points of a path contained in $[[ K_1,K_2 ]]$ differs from the
minimal one by at most \esp $4 + 2\frac{M-1}{N-1} + 4\log_2(N-1)$.}
\end{lem}

\vspace{0.15cm}

We now proceed to the proof of Theorem A. Recall that the graph of the $N_{d-1}^*$-orbit of the
interval $I^*$ defined by (\ref{def-interval}) has $\mathbb{Z}^{d-1}$ as its
set of vertices. We will hence inductively define parallelepipeds
$Q(n) \!\subset\! \mathbb{Z}^{d-1}$. We start with $Q(0) := [[ 1,1+4^{d+1} ]]^{d-1}$. Assuming that
\esp $Q(n)~:=~[[x_{1,n},y_{1,n}]] \times \cdots \times [[x_{d-1,n},y_{d-1,n}]]$ \esp has been already
defined, we let $i(n) \!\in\! \{1,\ldots,d-1\}$ be the residue class (mod. $d-1$) of $n$, and we set
\esp $Q(n+1) := \cdots \times [[1+4^{i(n)} (x_{i(n),n}-1), y_{i(n),n}]] \times 
[[x_{i(n)+1,n},1 + 4^{i(n)+1} (y_{i(n)+1,n}-1)]] \times \cdots,$
\esp where the dots mean that the corresponding factors remain untouched. (See Figure 5 for an
illustration of the case $d=4$, with $n \equiv 1$ (mod. $3$).)

Notice that $x_{i,n}, y_{i,n}$ are of the form $1+2^k$ for all $i,n$.
Although one may give precise formulas for $x_{i,n},y_{i,n}$, we will only need to record the
(easy to check) fact that for some universal constants $C_1,C_2,C_3,C_4$, we have the estimates
\begin{equation}\label{diferencia}
C_1 4^{\frac{in}{d-1}} \leq y_{i,n} - x_{i,n} \leq C_2 4^{\frac{in}{d-1}}
\end{equation}
and
\begin{equation}\label{crecimiento}
C_3 4^{\frac{in}{d-1}} \leq x_{i,n} \leq C_4 4^{\frac{in}{d-1}}.
\end{equation}
In particular, the number of points in $Q(n)$ is
\begin{equation}\label{number-of-points}
|Q(n)| = \prod_{j=1}^{d-1} (y_{j,n} - x_{j,n}) \leq \prod_{j=1}^{d-1} C_2 4^{\frac{jn}{d-1}}
= C_2^{d-1} 4^{\frac{n}{d-1} \sum_{j=1}^{d-1} j} = C_2^{d-1} 4^{\frac{nd}{2}}.
\end{equation}
Each $Q(n)$ is decomposed into paths pointing in the $i(n)^{th}$-direction as follows.
If $i(n) = 1$, then we decompose $Q(n)$ into ``horizontal'' paths of jump 1 at
each step, so that the number of paths is
$$\prod_{j \neq 1} (y_{j,n} - x_{j,n})
\geq \prod_{j \neq 1} C_1 4^{\frac{jn}{d-1}}
= C_1^{d-2} 4^{\frac{n}{d-1} \sum_{j \neq 1} j}
= C_1^{d-2} 4^{n \big[\frac{d}{2} - \frac{1}{d-1} \big]}.$$
If $i(n) \neq 1$, then for each fixed coordinates
$z_j \in [[x_{j,n},y_{j,n}]]$, with $j \neq i(n)$, we identify
$$\{z_1\} \times \cdots \times \{z_{i(n)-1}\} \times [[x_{i(n),n},y_{i(n),n}]]
\times \{z_{i(n)+1}\} \times \cdots \times \{z_{d-1}\} \esp \esp \sim
\esp\esp [[x_{i(n),n},y_{i(n),n}]] \subset \mathbb{N}$$
and we decompose this set into $N := x_{i(n)-1,n}$ paths making jumps (in the $i(n)^{th}$-direction)
of either 1 or $M := z_{i(n)-1,n}$ steps following the strategy of Lemma \ref{decompos}.
The corresponding number of paths now equals
\begin{small}
$$x_{i(n)-1,n} \prod_{j \neq i(n)} (y_{j,n} - x_{j,n})
\geq C_3 4^{\frac{(i(n)-1)n}{d-1}} \prod_{j \neq i(n)} C_1 4^{\frac{jn}{d-1}}
= C_3 4^{\frac{(i(n)-1) n}{d-1}} C_1^{d-2} 4^{\frac{n}{d-1} \sum_{j \neq i(n)} j}
= C_3 C_1^{d-2} 4^{n \big[ \frac{d}{2} - \frac{1}{d-1} \big]}.$$
\end{small}In either case, we denote by $Q(n,1), \ldots, Q(n,m_n)$ these paths, and
we let \esp $C_5 := \min \{ C_3 C_1^{d-2}, C_1^{d-2} \}$, \esp so that \esp $m_n \geq C_5
4^{n \big[ \frac{d}{2} - \frac{n}{d-1} \big] }$. \esp What is important in the construction
above is that each of these paths has a concrete dynamical meaning for the action of
$N_{d-1}^{*} \subset N_d$. Namely, if $i(n) =1$, they are induced by the action of the generator
$f_{2,1}$, whereas for $i(n) \neq 1$, they are induced by the action of $f_{i(n),1}$ and
$f_{i(n)-1,i(n)}$, where the first generator appears for 1-step jumps and the second one
for jumps of amplitude $z_{i(n)-1,n}$.

Associated to each point $(i_1,\ldots,i_{d-1}) \in \mathbb{Z}^{d-1}$ there is a
positive number $\ell_{i_1,\ldots,i_{d-1}}$, namely the length of the interval
$$I_{i_1,\ldots,i_{d-1}}^* := \bigcup_{j \in \mathbb{Z}} I_{i_1,\ldots,i_{d-1},j}.$$
Notice that the total sum of the $\ell_{i_1,\ldots,i_{d-1}}$'s equals 1. Moreover,
all the intervals $I_{i_1,\ldots,i_{d-1}}^*$ are in the $N_{d-1}^*$-orbit of
$I^* = I_{0,\ldots,0}^*$; see (\ref{def-interval}). Hence, as in the case $d = 3$,
what we need to do is to ensure the existence of an infinite sequence of
intersecting paths in $Q(1), Q(2), \ldots$ along which the total $L_{\alpha}$-sum
is finite provided that $\alpha > \frac{2}{d(d-1)}$.
To do this, we start with the next

\vspace{0.15cm}

\begin{lem} {\em Given \esp $0 < \alpha < 1$, there exists a constant
$C_6 > 0$ such that for all $A > 0$ and all $n \geq 1$, the subset
of indexes $m \in \{1, \ldots, m_n \}$ satisfying}
\begin{equation}\label{estimate-general}
\sum_{(i_1,\ldots,i_{d-1}) \in Q(n,m)} \ell_{i_1,\ldots,i_{d-1}}^{\alpha} \leq
\frac{A \esp C_6} {4^{n \big[ \frac{d \alpha}{2} - \frac{1}{d-1} \big] }}
\end{equation}
{\em has density at least $1 - 1/A$.}
\end{lem}

\noindent{\bf Proof.} As in the case $d = 3$, by H\"older's inequality we have
$$\sum_{(i_1,\ldots,i_{d-1}) \in Q(n)} \ell_{i_1,\ldots,i_{d-1}}^{\alpha} \leq
|Q(n)|^{1 - \alpha} \leq C_2^{(d-1)(1 - \alpha)} 4^{\frac{n d (1 - \alpha)}{2}}.$$
Hence,
$$\frac{1}{m_n} \sum_{m=1}^{m_n} \sum_{(i_1,\ldots,i_{d-1}) \in Q(n,m)}
\ell_{i_1,\ldots,i_{d-1}}^{\alpha}
\leq \frac{C_2^{(d-1)(1 - \alpha)} 4^{\frac{n d (1 - \alpha)}{2}}}
{C_5 4^{\frac{n d}{2} - \frac{n}{d-1}}}
= \frac{C_2^{(d-1)(1 - \alpha)}} {C_5 4^{\frac{n d \alpha}{2} - \frac{n}{d-1}}},$$
and the claim of the lemma follows from Chebyshev's inequality for
\esp $C_6 := C_2^{(d-1)(1-\alpha)} / C_5$. $\hfill\square$

\vspace{0.55cm}


\beginpicture

\setcoordinatesystem units <0.825cm,0.825cm>

\putrule from -4 0 to 4 0
\putrule from 0 0 to 0 1.5
\putrule from -4 1.5 to 4 1.5
\putrule from -4 0 to -4 1.5
\putrule from 4 0 to 4 1.5
\putrule from 5 1 to 5 2.5
\putrule from 7 3 to 7 5.3
\putrule from -3 2.5 to 5 2.5
\putrule from 1.8 3.3 to 1.8 7.8
\plot 1.8 7.8 3 9 /
\putrule from 5.8 3.3 to 5.8 7.8
\plot 5.8 7.8 7 9 /

\putrule from 3 9 to 15 9
\putrule from 1.8 7.8 to 13.8 7.8
\plot 13.8 7.8 15 9 /

\putrule from 1.8 5.3 to 13.8 5.3
\plot 13.8 5.3 15 6.5 /

\putrule from 15 6.5 to 15 9
\putrule from 13.8 5.3 to 13.8 7.8

\putrule from 1.8 3.3 to 5.8 3.3
\putrule from 5.8 3.3 to 5.8 1.8

\plot 4 0 7 3 /
\plot -4 1.5 -3 2.5 /
\plot 4 1.5 7 4.5 /
\plot 0 1.5 1.8 3.3 /

\plot 2.835 1.25 5.085 3.5 /
\plot 2.85 1.25 5.1 3.5 /
\plot 2.855 1.254 5.102 3.504 /
\plot 2.865 1.265 5.1 3.5065 /

\putrule from -2.4 1.24 to 2.85 1.24
\putrule from -2.4 1.23 to 2.85 1.23
\putrule from -2.4 1.22 to 2.85 1.22

\putrule from 5.1 3.5 to 5.1 7.215
\putrule from 5.118 3.5 to 5.118 7.215
\putrule from 5.09 3.5 to 5.09 7.215

\putrule from 5.1 7.215 to 11.1 7.215
\putrule from 5.1 7.2 to 11.1 7.2
\putrule from 5.1 7.185 to 11.1 7.185

\plot 11.1 7.2 13.7 9.8 /
\plot 11.082 7.182 13.682 9.782 /
\plot 11.07 7.17 13.67 9.77 /

\putrule from -2.4 1.23 to -2.4 -0.5
\putrule from -2.41 1.23 to -2.41 -0.5
\putrule from -2.39 1.23 to -2.39 -0.5


\setdots
\plot 0 0 3 3 /
\putrule from 3 3 to 3 4.5
\putrule from 3 3 to 7 3

\putrule from -3 1 to 5 1
\putrule from -3 1 to -3 2.5
\plot -4 0 -3 1 /
\putrule from 1 1 to 1 2.5

\putrule from 1.8 3.3 to 1.8 2
\putrule from 1.8 1.8 to 5.8 1.8

\putrule from 3 6.5 to 15 6.5
\plot 2 5.5 3 6.5 /

\putrule from 3 4.5 to 3 9
\putrule from 7 5.3 to 7 9
\plot 5.8 5.3 7 6.5 /
\putrule from 3 4.5 to 7 4.5

\plot 1.8 3.3 3 4.5 /


\small

\put{$\bullet$} at -2.4 1
\put{$\bullet$} at -2.4 0.25
\put{$\bullet$} at -2.4 -0.5
\put{$\bullet$} at -2.4 -0.25

\put{$\bullet$} at 5.25 3.63
\put{$\bullet$} at 5.1 3.48
\put{$\bullet$} at 4.95 3.33
\put{$\bullet$} at 4.05 2.43
\put{$\bullet$} at 3.15 1.53
\put{$\bullet$} at 3 1.38

\put{$\bullet$} at -3.4 1.23
\put{$\bullet$} at -3.15 1.23
\put{$\bullet$} at -2.9 1.23
\put{$\bullet$} at -2.65 1.23
\put{$\bullet$} at -2.4 1.23
\put{$\bullet$} at -2.15 1.23
\put{$\bullet$} at -1.9 1.23
\put{$\bullet$} at -1.65 1.23
\put{$\bullet$} at -1.4 1.23
\put{$\bullet$} at -1.15 1.23
\put{$\bullet$} at -0.9 1.23
\put{$\bullet$} at -0.65 1.23
\put{$\bullet$} at -0.4 1.23
\put{$\bullet$} at -0.15 1.23
\put{$\bullet$} at 0.1 1.23
\put{$\bullet$} at 0.35 1.23
\put{$\bullet$} at 0.6 1.23
\put{$\bullet$} at 0.85 1.23
\put{$\bullet$} at 1.1 1.23
\put{$\bullet$} at 1.35 1.23
\put{$\bullet$} at 1.6 1.23
\put{$\bullet$} at 1.85 1.23
\put{$\bullet$} at 2.1 1.23
\put{$\bullet$} at 2.35 1.23
\put{$\bullet$} at 2.6 1.23
\put{$\bullet$} at 2.85 1.23
\put{$\bullet$} at 3.1 1.23
\put{$\bullet$} at 3.35 1.23
\put{$\bullet$} at 3.6 1.23
\put{$\bullet$} at 3.85 1.23
\put{$\bullet$} at 4.1 1.23
\put{$\bullet$} at 4.35 1.23

\put{$\bullet$} at 5.1 2.98
\put{$\bullet$} at 5.1 3.25
\put{$\bullet$} at 5.1 3.48
\put{$\bullet$} at 5.1 3.715
\put{$\bullet$} at 5.1 4.95
\put{$\bullet$} at 5.1 5.2
\put{$\bullet$} at 5.1 6.7
\put{$\bullet$} at 5.1 6.95
\put{$\bullet$} at 5.1 7.2
\put{$\bullet$} at 5.1 8.7

\put{$\bullet$} at 2.85 7.2
\put{$\bullet$} at 3.1 7.2
\put{$\bullet$} at 3.35 7.2
\put{$\bullet$} at 3.6 7.2
\put{$\bullet$} at 3.85 7.2
\put{$\bullet$} at 4.1 7.2
\put{$\bullet$} at 4.35 7.2
\put{$\bullet$} at 4.6 7.2
\put{$\bullet$} at 4.85 7.2
\put{$\bullet$} at 5.1 7.2
\put{$\bullet$} at 5.35 7.2
\put{$\bullet$} at 5.6 7.2
\put{$\bullet$} at 5.85 7.2
\put{$\bullet$} at 6.1 7.2
\put{$\bullet$} at 6.35 7.2
\put{$\bullet$} at 6.6 7.2
\put{$\bullet$} at 6.85 7.2
\put{$\bullet$} at 7.1 7.2
\put{$\bullet$} at 7.35 7.2
\put{$\bullet$} at 7.6 7.2
\put{$\bullet$} at 7.85 7.2
\put{$\bullet$} at 8.1 7.2
\put{$\bullet$} at 8.35 7.2
\put{$\bullet$} at 8.6 7.2
\put{$\bullet$} at 8.85 7.2
\put{$\bullet$} at 9.1 7.2
\put{$\bullet$} at 9.35 7.2
\put{$\bullet$} at 9.6 7.2
\put{$\bullet$} at 9.85 7.2
\put{$\bullet$} at 10.1 7.2
\put{$\bullet$} at 10.35 7.2
\put{$\bullet$} at 10.6 7.2
\put{$\bullet$} at 10.85 7.2
\put{$\bullet$} at 11.1 7.2
\put{$\bullet$} at 11.35 7.2
\put{$\bullet$} at 11.6 7.2
\put{$\bullet$} at 11.85 7.2
\put{$\bullet$} at 12.1 7.2
\put{$\bullet$} at 12.35 7.2
\put{$\bullet$} at 12.6 7.2
\put{$\bullet$} at 12.85 7.2
\put{$\bullet$} at 13.1 7.2
\put{$\bullet$} at 13.35 7.2
\put{$\bullet$} at 13.6 7.2
\put{$\bullet$} at 13.85 7.2
\put{$\bullet$} at 14.1 7.2
\put{$\bullet$} at 14.35 7.2
\put{$\bullet$} at 14.6 7.2

\put{$\bullet$} at 10.8 6.9
\put{$\bullet$} at 10.95 7.05
\put{$\bullet$} at 11.1 7.2
\put{$\bullet$} at 11.25 7.35
\put{$\bullet$} at 12.15 8.25
\put{$\bullet$} at 13 9.1
\put{$\bullet$} at 13.75 9.85

\put{Figure 5} at 5.16 -1.1
\put{$Q(n)$} at -1.5 2.8
\put{$Q(n+1)$} at 0.5 3
\put{$Q(n+2)$} at 0.94 5.2
\put{$Q(n+3)$} at 9.5 9.3

\put{$x_{1,n}$} at -4 -0.43
\put{$y_{1,n} = y_{1,n+1}$} at 4 -0.43
\put{$x_{1,n+1}$} at 0 -0.43

\put{$x_{2,n} = x_{2,n+1}$} at 5.5 0
\put{$y_{2,n}$} at 5.7 1
\put{$y_{2,n+1} = y_{2,n+2}$} at 8.6 3
\put{$x_{2,n+2}$} at 6.7 1.8

\put{$x_{3,n}$} at -4.5 0.1
\put{$y_{3,n}$} at -4.5 1.6

\put{} at -4.5 0

\endpicture


\vspace{0.5cm}

From now on, we fix \esp $\alpha > \frac{d(d-1)}{2}$. \esp
We start by letting $r_n$ (resp. $s_n$) be the minimum (resp. maximum) of points
in a path of the form $Q(n,m)$ inside $Q(n) \cap Q(n+1)$. Similarly, we denote by
$r_n'$ (resp. $s_n'$) the minimum (resp. maximum) number of points of a path of
the form $Q(n+1,m)$ inside $Q(n) \cap Q(n+1)$. Then we let
$$B:= \prod_{n \geq 1} \frac{s_n}{r_n} \frac{s_n'}{r_n'}.$$
We claim that the value of $B$ is finite. Indeed, we have $r_n = s_n$ whenever 
$i(n) = 1$, whereas $s_n' = r_n'$ whenever $i(n) = d-1$. For the other values 
of $i := i(n)$, the condition (ii) in Lemma \ref{decompos} together with the 
inequalities \esp $2 \frac{y_{i-1,n}-1}{x_{i-1,n}-1} \leq 4^{d+2}$ \esp and \esp 
$2 \frac{y_{i,n+1}-1}{x_{i,n+1}-1} \leq 4^{d+2}$ \esp yield the estimates
$$\frac{y_{i,n+1} - x_{i,n+1}}{x_{i-1,n}} - 4 - 4^{d+2} - 4\log_2 (x_{i-1,n}-1) 
\esp \leq \esp r_n \esp \leq s_n \esp \leq \esp 
\frac{y_{i,n+1} - x_{i,n+1}}{x_{i-1,n}} + 4 + 4^{d+2} + 4\log_2 (x_{i-1,n}-1)$$
and
$$\frac{y_{i+1,n}-x_{i+1,n}}{x_{i,n+1}} - 4 - 4^{d+2} 
- 4 \log_2 ( x_{i,n+1}-1 ) 
\esp \leq \esp r_n' \esp \leq \esp s_n' \esp \leq \esp 
\frac{y_{i+1,n}-x_{i+1,n}}{x_{i,n+1}} + 4 + 4^{d+2} + 4\log_2 ( x_{i,n+1}-1 ),$$
which together with (\ref{diferencia}) and (\ref{crecimiento}) easily imply the
finiteness of $B$.

We will also use the (finite) constant
$$C:= 2 \sum_{n \geq 1} \frac{1}{2^{n \big[ \frac{d \alpha}{2} - \frac{1}{d-1} \big] }}.$$
Now we fix \esp $A_1 \geq B \esp C \esp 2^{\frac{d \alpha}{2} - \frac{1}{d-1}}$ \esp
such (\ref{estimate-general}) holds for $n=1$ and {\em every} \esp
$m \in \{1,\ldots,m_1\}$ \esp when letting $A = A_1$.
Finally, for $n \geq 2$, we set
$$A_n := B \esp C \esp 2^{n \big[ \frac{d \alpha}{2} - \frac{1}{d-1} \big]}.$$

\vspace{0.1cm}

\begin{lem} \label{al-fin}
{\em There exists an infinite sequence of paths of the form $Q(n,m_n')$
in $Q(n)$ such that, for all $n \geq 1$, the path $Q(n+1,m_{n+1}')$ intersects
$Q(n,m_n')$ and} (\ref{estimate-general}) {\em holds for $m = m_n'$ and $A = A_n$.}
\end{lem}

\noindent{\bf Proof.} As in the case $d = 3$, for each $n \geq 1$ we let $D_n$
be the density of indexes $m \in \{ 1,\ldots,m_n \}$ such that there exists a finite
sequence of paths $Q(1,n_1'), \ldots, Q(n,m_n')$ satisfying:

\vspace{0.1cm}

\noindent -- For each $1 \leq k \leq n-1$, the path $Q(k+1,m_{k+1}')$ intersects $Q(k,m_k')$;

\vspace{0.1cm}

\noindent -- Estimate (\ref{estimate-general}) holds for each $m = m_k'$ and $A = A_k$.

\vspace{0.1cm}

\noindent Similar arguments to those leading to (\ref{key}) yield
$$(1-D_{n+1}) \leq (1-D_n) \frac{s_n}{r_n} \frac{s_n'}{r_n'} + \frac{1}{A_n}.$$
Indeed, the product $\frac{s_n s_n'}{r_n r_n'}$ acts as a correction factor for the
passage from $Q(n)$ to $Q(n+1)$ taking into account that the paths of the form
$Q(n,m)$ do not have the same number of points in $Q(n) \cap Q(n+1)$, and similarly
for those of the form $Q(n+1,m)$. By induction, the preceding inequality yields
$$1- D_n \leq (1 - D_1) \prod_{k=1}^{n-1} \frac{s_k}{r_k} \frac{s_k'}{r_k'} +
\sum_{k=1}^{n-1} \frac{1}{A_k} \prod_{j=1}^{k-1} \frac{s_j}{r_j} \frac{s_j'}{r_j'}
\leq (1-D_1) B + B \sum_{k \geq 1} \frac{1}{A_k} .$$
The choice of $A_1$ was made so that $D_1 = 1$, hence \esp
$$1-D_n \leq B \sum_{k \geq 1} \frac{1}{A_k} \leq \frac{1}{2}.$$
As a consequence, $D_n \geq 1/2$, which implies that for each $n$ we may obtain
a finite sequence of $n$ paths with the desired properties. The infinite sequence
is obtained via a Cantor diagonal type argument. $\hfill\square$

\vspace{0.45cm}

The proof of Theorem A is now at hand. Indeed, the concatenation of the paths provided
by the preceding lemma yields an infinite sequence of points in $\mathbb{Z}^{d-1}$
along which the value of the $L_{\alpha}$-sum is bounded from above by
$$\sum_{n \geq 1} \frac{A_n C_6}{4^{n \big[ \frac{d \alpha}{2} - \frac{1}{d-1} \big]}}
\leq \frac{A_1 C_6}{4^{\frac{d \alpha}{2} - \frac{1}{d-1} }} +
\sum_{n \geq 2} \frac{B \esp C \esp C_6}{2^{n \big[ \frac{d
\alpha}{2} - \frac{1}{d-1} \big]}} \leq 2 A_1 C_6 + 2 B C_6.$$
This is in correspondence to a sequence of intervals of the form
$I_{i_1,\ldots,i_{d-1}}$ each of which is obtained from the
preceding one by applying one of the generators in
$\{f_{2,1},f_{3,1}, \ldots, f_{d,1} \} \cup \{ f_{2,1}, f_{3,2}.
\ldots, f_{d,d-1} \}$. Joining this infinite sequence to a finite
one from the origin to a point in $Q(1,n_1')$, we obtain an
infinite sequence of intervals in the $N_{d-1}^*$-orbit of the
interval $I^*$ for which the $L_{\alpha}$-sum is finite, and hence
the arguments of \S \ref{main-idea} may be applied. This concludes
the proof.


\subsection{An independent combinatorial lemma}
\label{section-comb}

\hspace{0.45cm} The aim of this Section is to give the proof of Lemma \ref{decompos}. 
We first give the details of the construction of the partition of $\mathbb{N}_0$ into 
$N$ sets (paths) $P_1,....,P_N$, and latter we check the desired properties. The 
construction is made in two steps, the former of which applies to arbitrary values 
of $N$, whereas the latter is restricted to integers of the form $1 + 2^k$.

\vspace{0.3cm}

\noindent{\bf Step 1.} Let $M>N$ be positive integers. Assume that we are given a partition 
$$[[ 0, M-1 ]] = R_0 \bigcup R_1 \bigcup \ldots \bigcup R_{N-1}$$ 
into ``consecutive'' sets, that is, such that $1 + \max {R_i} = \min {R_{i+1}}$ holds 
for all $0 \leq i \leq N-2$. Then this induces a partition of $\mathbb{N}_0$ 
as follows. Denoting $R \oplus k := \{n+k \!: n \in R \}$, we define

\vspace{0.1cm}

\noindent $\bullet$ \esp $S_1 := \bigcup \limits_{j=1}^{N-1} R_j \oplus j (M-1)$,

\vspace{0.1cm}

\noindent $\bullet$ \esp $S_i := \bigcup \limits_{j=i-1}^{N-1} R_j \oplus (j-i+1) (M-1)  
\hspace{0.1cm} \bigcup \hspace{0.1cm} 
\bigcup\limits_{j=1}^{i-2} R_j \oplus (j-i+N)(M-1)$, \esp for $2 \leq i \leq N$.

\vspace{0.1cm}

\noindent (Notice that, by definition, the second term in the definition of $S_i$ 
above is empty for $i=2$.) Now, what defines our partition of $\mathbb{N}_0$ is 
the ``periodic repetition'' of the sets $S_1,...,S_N$. More precisely, we let

\vspace{0.1cm}

\noindent $\bullet$ \esp 
$P_1 := R_0 \bigcup \bigcup \limits_{k=0}^{\infty} S_1 \oplus kN(M-1) $,

\vspace{0.1cm}

\noindent $\bullet$ \esp 
$P_i := \bigcup\limits_{k=0}^{\infty} S_i \oplus kN(M-1)$, \esp for $2 \leq i \leq N$.

\vspace{0.1cm}

To have a clearer view of this construction, the reader may easily check that for the 
particular choice $R_0 := \{0\}, \esp R_1 := \{1\}, \ldots, \esp R_{N-2} := \{N-2\}$ 
and $R_{N-1} := \{N-1,N,N+1,\ldots,M-1\}$, it yields to the paths constructed in 
\S \ref{caso d=4} (see again Figure 3 for an illustration).

It is sometimes better to think on our paths as concatenations of ``patches''. 
In this view, for $2 \leq i \leq N$, the sequence representing $S_i$ is \esp 
$R_{i-1} R_{i} \ldots R_{N-1} R_1 R_2 \ldots R_{i-2}$, which in notation 
modulo $N-1$ may be rewritten as \esp $R_{i-1} R_i \ldots R_{i+N-2}$. \esp 
This means that $S_i$ is made of a copy of $R_{i-1}$ followed by a copy of $R_i$ 
translated by $M-1$ units, a copy of $R_{i+1}$ translated by another $M-1$ units, 
and so on. \esp Similarly, our paths $P_i$ may be seen as infinite sequences of 
patches. Thinking on each $S_i$ as a patch as well, for $2 \leq i \leq N$, the 
path $P_i$ is represented by \esp $S_i S_i S_i \ldots$. \esp The sequence 
representing $P_1$ corresponds to \esp $R_0 S_1 S_1 S_1 \ldots$.

\vspace{0.3cm}

\noindent{\bf Step 2.} Assuming that $N$ has the form $1+2^k$, we 
will associate to it a particular choice of sets $R_1,...,R_N$. 
Let $p \geq 1$ and $q \geq 0$ be the integers such that
$$M-1 = (N-1)p + q, \quad \mbox{with } q < N-1.$$
Let us consider the binary expansion of $q$:
$$q = 2^{r_1} + \ldots + 2^{r_l}, \quad \mbox{with } 
r_1 >....> r_l \geq 0.$$ 
(Notice that since $q < N-1= 2^k$, we have $k > r_1$.) 
Now, for $1 \leq i \leq N-1$, define $s_i$ as being the largest integer $s$ 
such that $2^{k-r_s}$ divides $i$ whenever there is such an index $s$, and as 
being equal to zero otherwise. We claim that the following relation holds:
\begin{equation}\label{igualdad-clave}
s_1 + s_2 + \ldots + s_{N-1} = q.
\end{equation}
Indeed, by definition, $s_i$ equals $s > 0$ if and only if $i$ is a multiple 
of $2^{k-r_s}$ but not a multiple of $2^{k-r_{s+1}}$. Now, in $\{1,2,\ldots,N-1\}$, 
there are exactly $2^{r_s}$ multiples of $2^{k-r_s}$, namely the products of 
$2^{k-r_s}$ with the integers in $\{1,2,3,\ldots, 2^{r_s}\}$. Hence, the 
left-side expression in (\ref{igualdad-clave}) equals 
\begin{equation}\label{gen}
\sum_{s=1}^l s \esp \big| \{ i \!: \esp s_i = s \} \big| 
\esp = \esp 
\sum_{s=1}^{l-1} s \big( 2^{r_s} - 2^{r_{s+1}} \big) \esp + \esp\esp l 2^{r_l} 
\esp = \esp 
\sum_{s=1}^l 2^{r_s} 
\esp = \esp q.
\end{equation}

Finally, let us inductively define:

\vspace{0.1cm}

\noindent $\bullet$ \esp $R_0 := \{ 0 \}$,

\vspace{0.1cm}

\noindent $\bullet$ \esp $R_i := \{1+ \max {R_{i-1}},...,p + s_i + \max {R_{i-1}} \}$, 
\esp where $1 \leq i \leq N-1$.

\vspace{0.1cm}

Notice that for $1 \leq i \leq N-1$, the number of points of $R_i$ equals
\begin{equation}\label{numero-de-ptos}
p + s_i \leq p + l \leq p + k = p + \log_2 (N-1).
\end{equation} 
Using (\ref{igualdad-clave}), we conclude that the number 
of points contained in the union of the $R_i$'s equals  
$$1 + p(N-1) + s_1 + \ldots +s_{N-1} \esp
= \esp 1 + p(N-1) + q \esp = \esp M.$$ 
Thus, the $R_i$'s \esp yield a partition of $[[0,M-1]]$ into consecutive sets. We claim that 
the corresponding partition of $\mathbb{N}_0$ into the paths $P_1,\ldots,P_N$ produced as in 
Step 1 satisfies the desired properties. 

\vspace{0.3cm}

\noindent{\bf Step 3.} We first notice that in order to prove property (ii) of Lemma 
\ref{decompos}, we may restrict ourselves to intervals of the form $[[0,K]]$ instead 
of general intervals $[[K_1,K_2]]$ provided we obtain the better bound \esp 
$2 + \frac{M-1}{N-1} + 2 \log_2 (N)$ \esp for the maximal difference of 
points in $[[0,K]]$ among our $N$ paths. This is what we now proceed to do. 

Let $a,b$ be non-negative integers such that 
$$K = a N(M-1) + b, \quad \mbox{with } \esp b < N(M-1),$$
Let us first consider a path $P_i$ such that $2 \leq i \leq N$. 
In terms of patch sequences, and using notation modulo $N-1$, 
the intersection of $P_i$ with $[[0,K]]$ has the form 
$$S_i \ldots S_i R_{i-1} R_i \ldots R_{i-1+t} T, \quad \mbox{with } \esp t \leq N-1.$$
Here, the patch $T$ is a starting part of the patch $R_{i+t}$. 
Moreover, the patch $S_i$ appears precisely $a$ times.

By construction, the number of points in the set represented above is $a$ times the number 
of points in $S_i$ plus the sum of the number of points in \esp $R_{i-1} \ldots R_{i-1+t}$ \esp 
plus the number of points in $T$. The former equals $a(M-1)$, hence it is independent of 
$i \in \{2,\ldots,N\}$, 
whereas the latter is smaller than or equal to \esp $p + s_{i+t} \leq p + \log_2 (N-1)$; \esp 
see (\ref{numero-de-ptos}). As a consequence, the difference with respect to the number of 
points in $[[0,K]] \cap P_j$ (with $2 \leq j \leq N$) is at most $p + \log_2 (N-1)$ plus 
the difference between the number of points in \esp $R_{i-1} \ldots R_{i-1+t}$ \esp 
and \esp $R_{j-1} \ldots R_{j-1+t}$. Since \esp $p \leq 1 + \frac{M-1}{N-1}$, \esp 
our task reduces to show that the last difference is at most \esp $\log_2 (N-1)$.

Now, the number of points in the first (resp. second) sequence above equals
$$( p + s_{i-1} ) + ( p + s_i ) + \ldots + ( p + s_{i-1+t} )  
\esp = \esp tp + s_{i-1} + \ldots + s_{i-1+t}$$
$$\big(\mbox{resp. } ( p + s_{j-1} ) + ( p + s_j ) + \ldots + 
( p + s_{j-1+t} ) \esp = \esp tp + s_{j-1} + \ldots + s_{j-1+t} \big).$$
Define $\rho_{s,i}$ (resp. $\rho_{s,j}$) as being the number of indexes in 
$\{i-1,\ldots,i-1+t\}$ (resp. $\{j-1,\ldots,j-1+t\}$) that are multiples of 
$2^{k-r_s}$. A similar argument to that leading to (\ref{gen}) yields 
$$s_{i-1} + \ldots + s_{i-1+t} = \rho_{1,i} + \rho_{2,i} + \ldots + \rho_{l,i} 
\quad \big( \mbox{resp. } 
s_{j-1} + \ldots + s_{j-1+t} = \rho_{1,j} + \rho_{2,j} + \ldots + \rho_{l,j} \big).$$
Since 
$$\frac{t}{2^{k-r_s}} \leq \rho_{s,i} \leq 1 + \frac{t}{2^{k-r_s}} \quad 
\big( \mbox{resp. } \frac{t}{2^{k-r_s}} \leq \rho_{s,j} \leq 1 + \frac{t}{2^{k-r_s}} \big),$$
the value of \esp $|\rho_{s,i} - \rho_{s,j}|$ \esp equals either zero or 1.
We thus conclude that 
$$\big| s_{i-1} + \ldots + s_{i-1+t} - s_{j-1} - \ldots - s_{j-1+t} \big| \esp 
\leq |\rho_{1,i} - \rho_{1,j}| + \ldots + |\rho_{l,i} - \rho_{l,j}| 
\esp \leq \esp l \esp \leq \esp k \esp = \esp \log_2 (N-1),$$
as we wanted to show. 

Actually, so far we have obtained the upper bound \esp $1 + \frac{M-1}{N-1} + 2\log_2 (N-1)$ \esp 
for the difference between the number of points in $P_i \cap [[0,K]]$ and $P_j \cap [[0,K]]$. 
The extra 1 which lacks appears when making comparisons with the path $P_1$, taking into 
account that $P_1$ starts with $R_0 = \{ 0 \}$. The proof of this follows the same 
ideas above. We leave the details to the reader.


\section{Construction of smoothings for $\alpha < \frac{2}{d(d-1)}$}

\subsection{A reminder on Denjoy-Pixton actions}
\label{denjoy-pixton}

\hspace{0.45cm} For the constructions leading to the proofs of Theorems 
B and C, we will use Pixton's technique \cite{Pix}. The main technical 
tool will be the following lemma from \cite{TsP}.

\vspace{0.1cm}

\begin{lem} {\em For a certain universal constant $M$ there exists a family of diffeomorphisms 
$\varphi_{I',I}^{J',J} \!: I \rightarrow J$ where $I,I',J,J'$ are non-degenerate intervals 
and $I'$ (resp. $J'$) is contiguous by the left to $I$ (resp. $J$), satisfying \esp  
$\varphi_{J',J}^{K',K} \circ \varphi_{I',I}^{J',J} = \varphi_{I',I}^{K',K}$ \esp and}
$$\frac{\big| \log \big( D \varphi_{I',I}^{J',J}(u) \big) - 
\log \big( D \varphi_{I',I}^{J',J}(v) \big) \big|}{|u - v|} 
\esp \leq \esp \frac{M}{|I|} \left| \frac{|I| |J'|}{|J| |I'|} - 1 \right| $$ 
for all $u,v$ in $I$ provided that \esp \esp \esp 
$\max \{ |I|,|I'|,|J|,|J'| \} \leq 2 \min \{ |I|,|I'|,|J|,|J'| \}$.
\label{cle}
\end{lem}

\vspace{0.1cm} 

The proof of this lemma proceeds as follows. Following \cite{TsP}, 
let $\xi(x)(\frac{\partial}{\partial x})$ be a $C^{\infty}$ vector 
field on $[0,1]$ such that \esp $\xi(x)=x$ \esp near $0$, and \esp 
$\xi(x)=0$ \esp on $[1/2,1]$. 
Moreover, assume that for all $x$,   
$$\left\lvert D \xi (x) \right\lvert \leq 1.$$
Let $\psi_{t}(x)$ be the solution of the differential equation
$$\frac{d\psi_{t}}{dt}(x)=\xi(\psi_{t}(x)), \quad \psi_{0}(x)=x.$$
Let us consider the diffeomorphism $x \mapsto b \esp \esp \psi_{t}(x/a)$ 
sending the interval $[0,a]$ onto the interval $[0,b]$. For any real 
numbers $a',a,b',b$ such that \esp $a'<0<a$ \esp and \esp $b'<0<b$, \esp let 
$\phi_{a',a}^{b',b}$ be the diffeomorphism from $[0,a]$ onto $[0,b]$ defined by 
$$\phi_{a',a}^{b',b}(x) = b \esp \psi_{\log (b'a/a'b)}(x/a).$$
Its is easy to check that for all positive $a,b,c$ and all negative $a',b',c'$, one has
$$\phi_{b',b}^{c',c} \circ \phi_{a',a}^{b',b} \esp = \esp \phi_{a',a}^{c',c}$$
Moreover, as is shown in \cite{TsP},
\begin{equation}
\log D\phi_{a',a}^{b',b} (x) = 
\log \frac{b}{a} + \log D\psi_{\log(b'a/a'b)} \left(\frac{x}{a}\right),
\label{uno-tsuboi}
\end{equation}
\begin{equation}
\left\lvert\log D\psi_{\log(b'a/a'b)}\right\lvert 
\leq \left\lvert\log\frac{b'a}{a'b}\right\lvert 
= \left\lvert\log\frac{b'}{a'}-\log\frac{b}{a}\right\lvert.
\label{dos-tsuboi}
\end{equation}
Furthermore, letting $M > 0$ be a constant such that 
$\lvert D^{2} \xi(x) \lvert\leq M$ for all $x$, we have 
\begin{equation}
\left\lvert D \log D \phi_{a',a}^{b',b} (x)\right\lvert 
\leq\frac{M}{a}\left\lvert\frac{b'a}{a'b}-1\right\lvert.
\label{tres-tsuboi}
\end{equation}
Starting with the maps $\phi_{a',a}^{b',b}$, we construct the desired family 
$\{\varphi_{I',I}^{J',J}\}$ 
as follows. Letting $I=[w,w+a]$, $I'=[w+a',w]$, $J=[w',w'+b]$, and $J'=[w'+b',w']$, 
where $a' < 0 < a$ and $b' < 0 < b$, we let 
$$\varphi_{I',I}^{J',J} = \phi_{a',a}^{b',b}(x-w) + w'.$$


\subsection{Sharp embeddings of $N_d$}

\hspace{0.45cm} In order to 
prove our Theorem B, we fix once and for all an arbitrary positive number 
$\alpha < \frac{2}{d(d-1)}$. Our aim is to show that for a good choice of the lengths 
$|I_{i_1,\ldots,i_d}|$, the maps $f_j := f_{j+1,j}$, $1 \leq j \leq d-1$, defined as in 
\S \ref{accion} using the maps from \S \ref{denjoy-pixton} instead of affine maps 
are $C^{1+\alpha}$-diffeomorphisms 
of the corresponding (non necessarily normalized) interval $I$. From now on, we will assume 
that $d \geq 3$. Although the case $d=2$ can be ruled out by a slightly modified construction, 
it is also covered by the (much simpler) construction leading to Theorem C. In all what follows, 
$M$ will denote a universal constant whose explicit value is irrelevant for our purposes.

\vsp

We begin by choosing number $p_d \!\in ]1,5/4]$, and for $1 \leq j \leq d-1$ 
we choose $p_j>0$ so that the following properties are satisfied:

\vsp\vsp

\noindent $\mathrm{(i_B)}$ \esp\esp\esp $p_1>p_2>\ldots>p_{d-1}>p_d>1$, 

\vsp\vsp

\noindent $\mathrm{(ii_B)}$ \esp\esp\esp
$\frac{1}{p_1}+\frac{1}{p_2}+\ldots+\frac{1}{p_{d-1}}+\frac{1}{p_d}<1$, 

\vsp\vsp

\noindent $\mathrm{(iii_B)}$ \esp\esp\esp $\alpha \leq \frac{p_d}{(p_d-1)p_1}$, 

\vsp\vsp

\noindent $\mathrm{(iv_B)}$ \esp\esp\esp 
$\alpha \leq \frac{p_d}{p_d-1}\left(\frac{1}{p_j}-\frac{1}{p_{j-1}}\right)$ 
\esp for all \esp $1<j<d$, 

\vsp\vsp

\noindent $\mathrm{(v_B)}$ \esp\esp\esp 
$\alpha \leq \frac{1}{p_d}-\frac{1}{p_{d-1}}$.

\vsp\vsp\vsp

A concrete choice is \esp $p_j := \frac{1}{j\alpha(1-1/p_d)}.$ \esp 
(Hence, one may take $p_d := \frac{5}{4}$ and $p_j := \frac{5}{j\alpha}$ for 
$1 \leq j \leq d-1$.) Indeed, the first property is easy 
to check. For the second one, we have 
$$\sum_{j=1}^{d}\frac{1}{p_j}=\frac{1}{p_d}+\sum_{j=1}^{d-1}j\alpha 
\left( 1-\frac{1}{p_d} \right) 
= \frac{1}{p_d}+\alpha \left( 1-\frac{1}{p_d} \right) \frac{d(d-1)}{2} 
< \frac{1}{p_d} + \left( 1 - \frac{1}{p_d} \right) = 1,$$
where the inequality comes from the hypothesis 
$\alpha < \frac{2}{d(d-1)}$. 
For the third and fourth properties, we actually have 
equalities with our choice. Finally, since $d \geq 3$, 
$$\alpha < \frac{2}{d(d-1)} < \frac{2}{3} \leq \frac{1/p_d}{2 - 1/p_d} \leq 
\frac{1/p_d}{1 + (d-1) (1 - 1/p_d)}.$$ 
Hence,
$$\alpha \left[ 1 + (d-1) \Big( 1 - \frac{1}{p_d} \Big) \right] \leq \frac{1}{p_d},$$ 
that is, 
$$\alpha \leq \frac{1}{p_d} - \alpha (d-1) \Big( 1 - \frac{1}{p_d} \Big) 
= \frac{1}{p_d} - \frac{1}{p_{d-1}},$$
which shows $\mathrm{(v_B)}$.

\vsp

It is worth mentioning that for \esp $\alpha \geq \frac{2}{d(d-1)}$, \esp 
the properties above are incompatible. Indeed, from $\mathrm{(iii_B)}$ we get 
\esp $\frac{1}{p_1} \geq \frac{\alpha (p_d - 1)}{p_d}$. \esp Using $\mathrm{(iv_B)}$ 
inductively, we obtain \esp $\frac{1}{p_j} \geq \frac{j \alpha (p_d - 1)}{p_d}$ 
\esp for $1 \leq j \leq d-1$. \esp This yields 
$$\sum_{j=1}^{d} \frac{1}{p_j} \esp\esp
\geq \esp\esp \sum_{j=1}^{d-1} \frac{j \alpha (p_d - 1)}{p_d} + \frac{1}{p_d} 
\esp\esp = \esp\esp \frac{\alpha (p_d - 1) \esp d \esp (d-1)}{2 p_d} + \frac{1}{p_d}.$$
If \esp \esp $\alpha \geq \frac{2}{d(d-1)}$, \esp \esp the right-side 
expression is greater than or equal to 1, contrary to $\mathrm{(ii_B)}$.
 
\vsp\vsp\vsp\vsp\vsp

Now fixing any choice of the $p_j$'s as above, we let 
$$|I_{i_1,\ldots,i_d}| \esp \esp := 
\esp \esp \frac{1}{|i_1|^{p_1} + \ldots + |i_d|^{p_d} + 1}.$$
According to \cite[\S 3]{KN}, property $\mathrm{(ii_B)}$ implies that the sum of the 
lengths $|I_{i_1,\ldots,i_d}|$ is finite. We next proceed to show that the induced 
maps $f_j$ are $C^{1+\alpha}$-diffeomorphisms of the corresponding interval $I$. 


\subsection{The map $f_1$ is a $C^{1+\alpha}$-diffeomorphim}

{\bf I.} First we consider $x,y$ in the same interval $I_{i_1,\ldots,i_d}$. We have
$$\frac{\vert \log Df_1(x) - \log Df_1(y) \vert}{\vert x-y\vert}\leq
\frac{M}{\vert I_{i_1,\ldots,i_d}\vert}\left\vert\frac{\vert 
I_{i_1,\ldots,i_d}\vert}{\vert I_{i_1+1,\ldots,i_d}\vert}\frac{\vert 
I_{i_1+1,\ldots,i_d-1}\vert}{\vert I_{i_1,\ldots,i_d-1}\vert}-1\right\vert.$$
Hence,
$$\frac{\vert \log Df_1(x) - \log Df_1(y) \vert}{\vert x-y\vert^{\alpha}}
\leq\ M\left\vert\frac{\vert I_{i_1,\ldots,i_d}\vert}{\vert I_{i_1+1,\ldots,i_d}\vert}
\frac{\vert I_{i_1+1,\ldots,i_d-1}\vert}{\vert I_{i_1,\ldots,i_d-1}\vert}-1\right\vert\vert 
I_{i_1,\ldots,i_d}\vert^{-\alpha}.$$
The right-side expression is bounded from above by 
$$M \left\lvert\frac{(\lvert i_1\lvert^{p_1}+\ldots+\lvert i_d-1\lvert^{p_d}+1)
(\lvert i_1+1\lvert^{p_1}+\ldots+\lvert i_d\lvert^{p_d}+1)}
{(\lvert i_1+1\lvert^{p_1}+\ldots+\lvert i_d-1\lvert^{p_d}+1)
(\lvert i_1\lvert^{p_1}+\ldots+\lvert i_d\lvert^{p_d}+1)}-1
\right\lvert(\lvert i_1\lvert^{p_1}
+\ldots+\lvert i_d\lvert^{p_d}+1)^{\alpha},$$
which equals
$$M\!\left\lvert\frac{(\lvert i_d\lvert^{p_d}-\lvert i_d-1\lvert^{p_d})
(\lvert i_1\lvert^{p_1}-\lvert i_1+1\lvert^{p_1})}{(\lvert i_1+1
\lvert^{p_1}+\lvert i_2\lvert^{p_2}+\ldots+\lvert i_d-1\lvert^{p_d}+1)
(\lvert i_1\lvert^{p_1}+\lvert i_2\lvert^{p_2}+\ldots+\lvert i_d\lvert^{p_d}+1)}\right\lvert\!
(\lvert i_1\lvert^{p_1}+\ldots+\lvert i_d\lvert^{p_d}+1)^{\alpha}\!\!\!.$$
By the Mean Value Theorem, this expression is bounded from above by
\begin{equation}\label{eq}
M\frac{(\lvert i_d\lvert+1)^{p_d-1}(\lvert i_1\lvert+1)^{p_1-1}}
{(\lvert i_1\lvert^{p_1}+\lvert i_2\lvert^{p_2}+\ldots+\lvert i_d\lvert^{p_d}+1)^{2-\alpha}}.
\end{equation}
In the case $\lvert i_1\lvert^{p_1}\leq\lvert i_d\lvert^{p_d}$, this is bounded by 
$$M\frac{(\lvert i_d\lvert+1)^{p_d-1}(\lvert i_d\lvert^{\frac{p_d}{p_1}}+1)^{p_1-1}}
{(\lvert i_d\lvert^{p_d}+1)^{2-\alpha}}.$$ 
This expression is uniformly bounded when 
$p_d-1+\frac{p_d}{p_1}(p_1-1)\leq p_d(2-\alpha)$, that is, when 
$\alpha\leq\frac{1}{p_1}+\frac{1}{p_d}$, which is ensured by the condition 
$\mathrm{(v_B)}$. \esp In the case $\lvert i_d\lvert^{p_d}\leq\lvert i_1\lvert^{p_1}$, 
we have the upper bound 
$$M\frac{(\lvert i_1\lvert^{\frac{p_1}{p_d}}+1)^{p_d-1}
(\lvert i_1\lvert+1)^{p_1-1}}{(\lvert i_1\lvert^{p_1}+1)^{2-\alpha}}.$$ 
This expression is uniformly bounded when 
$p_1-1+\frac{p_1}{p_d}(p_d-1)\leq p_1(2-\alpha)$, 
that is, when $\alpha\leq\frac{1}{p_1}+\frac{1}{p_d}$, 
which \esp --as we have already seen-- \esp is 
ensured by the condition $\mathrm{(v_B)}$.\\

\vsp\vsp

\noindent{\bf II.} Now we consider $x,y$ so that $x \in I_{i_1,\ldots,i_{d-1},i_d}$ 
and $y \in I_{i_1,\ldots,i_{d-1},i_d'}$ for some $i_d < i_d'$. To simplify, we will 
just deal with positive $i_d, i_d'$, the other cases being analogous.

If $i_d' = i_d +1$, then letting $z$ be the right endpoint of the interval 
$I_{i_1,\ldots,i_{d-1},i_d}$, we have
$$\frac{\vert \log Df_1(x) - \log Df_1(y) \vert}{|x-y|^{\alpha}} \leq 
\frac{\vert \log Df_1(x) - \log Df_1(z) \vert}{|x-z|^{\alpha}} + 
\frac{\vert \log Df_1(z) - \log Df_1(y) \vert}{|z-y|^{\alpha}},$$
and both terms of the sum above are uniformly bounded by the previous case. 

\vsp

Assume henceforth that $i_d'-i_d \geq 2$. By property (\ref{dos-tsuboi}), 
the value of \esp \esp $\vert \log Df_1(x) - \log Df_1(y) \vert$ \esp \esp 
is bounded from above by 
\begin{small}
$$\left\lvert\log\!\frac{\lvert I_{i_1+1,\ldots,i_d}\!\lvert}{\lvert I_{i_1,\ldots,i_d}\lvert} 
\!-\! \log\!\frac{\lvert I_{i_1+1,\ldots,i_d'}\!\lvert}{\lvert I_{i_1,\ldots,i_d'}\lvert}\!\right\lvert 
+\left\lvert\log\!\frac{\lvert I_{i_1+1,\ldots,i_d}\!\lvert}{\lvert I_{i_1,\ldots,i_d}\lvert} 
\!-\! \log\!\frac{\lvert I_{i_1+1,\ldots,i_d-1}\!\lvert}{\lvert I_{i_1,\ldots,i_d-1}\lvert}\!\right\lvert 
+ \left\lvert\log\!\frac{\lvert I_{i_1+1,\ldots,i_d'}\!\lvert}{\lvert I_{i_1,\ldots,i_d'}\lvert} 
\!-\! \log\!\frac{\lvert I_{i_1+1,\ldots,i_d'-1}\!\lvert}
{\lvert I_{i_1,\ldots,i_d'-1}\lvert}\!\right\lvert\!.$$
\end{small}Since $i \mapsto \frac{|I_{i_1 +1,\ldots,i_{d-1},i}|}{|I_{i_1,\ldots,i_{d-1},i}|}$ 
is a monotonous function, this expression is smaller than or equal to 
\begin{multline*}
3 \left\lvert\log \left( \frac{\lvert I_{i_1+1,\ldots,i_d'}\lvert}{\lvert I_{i_1,\ldots,i_d'}\lvert} \right)
- \log \left( 
\frac{\lvert I_{i_1+1,\ldots,i_d-1}\lvert}{\lvert I_{i_1,\ldots,i_d-1}\lvert} \right) \right\lvert =\\
= \left\lvert\log\frac{(\lvert i_1\lvert^{p_1}+\lvert i_2\lvert^{p_2}+\ldots+\lvert i_d -1 \lvert^{p_d}+1)
(\lvert i_1+1\lvert^{p_1}+\lvert i_2\lvert^{p_2}+\ldots+\lvert i_d'\lvert^{p_d}+1)}
{(\lvert i_1+1\lvert^{p_1}+\lvert i_2\lvert^{p_2}+\ldots+\lvert i_d - 1 \lvert^{p_d}+1)
(\lvert i_1\lvert^{p_1}+\lvert i_2\lvert^{p_2}+\ldots+\lvert i_d'\lvert^{p_d}+1)}\right\lvert =\\
= \left\lvert\log\left(1+\frac{(i_d'^{p_d}-(i_d-1)^{p_d})
(\lvert i_1\lvert^{p_1}-\lvert i_1+1\lvert^{p_1})}
{(\lvert i_1+1\lvert^{p_1}+\lvert i_2\lvert^{p_2}+\ldots+\lvert i_d-1 \lvert^{p_d}+1)
(\lvert i_1\lvert^{p_1}+\lvert i_2\lvert^{p_2}+\ldots+\lvert i_d'\lvert^{p_d}+1)}\right)\right\lvert.
\end{multline*}
Since the expression in brackets in the right-side term equals 
$$\frac{(\lvert i_1\lvert^{p_1}+\lvert i_2\lvert^{p_2}+\ldots+\lvert i_d-1\lvert^{p_d}+1)
(\lvert i_1+1\lvert^{p_1}+\lvert i_2\lvert^{p_2}+\ldots+\lvert i_d'\lvert^{p_d}+1)}
{(\lvert i_1+1\lvert^{p_1}+\lvert i_2\lvert^{p_2}+\ldots+\lvert i_d-1 \lvert^{p_d}+1)
(\lvert i_1\lvert^{p_1}+\lvert i_2\lvert^{p_2}+\ldots+\lvert i_d'\lvert^{p_d}+1)},$$ 
it is bounded from below by a positive number. Therefore, 
$$\big| \log Df_1(x) - \log Df_1(y) \big| \leq 
M \! \left\lvert\frac{(i_d'^{p_d}-i_d^{p_d})
(\lvert i_1\lvert^{p_1}-\lvert i_1+1\lvert^{p_1})}
{(\lvert i_1+1\lvert^{p_1}+\lvert i_2\lvert^{p_2}+\ldots+\lvert i_d\lvert^{p_d}+1)
(\lvert i_1\lvert^{p_1}+\lvert i_2\lvert^{p_2}+\ldots+\lvert i_d'\lvert^{p_d}+1)}\right\lvert\!.$$
By The Mean Value Theorem, the last expression is bounded from above by
$$M\frac{i_d'^{p_d-1}(i_d'-i_d)(\lvert i_1\lvert+1)^{p_1-1}}
{(\lvert i_1\lvert^{p_1}+\ldots+i_d^{p_d}+1)(\lvert i_1\lvert^{p_1}+\ldots+i_d'^{p_d}+1)}.$$
Thus, in order to get an upper bound for \esp\esp
$\frac{\vert \log f_1'(x) - \log f_1'(y) \vert}{\vert x-y\vert^{\alpha}},$ 
\esp\esp we need to estimate the expression
\begin{equation}\label{vvv}
\frac{i_d'^{p_d-1}(i_d'-i_d)(\lvert i_1\lvert+1)^{p_1-1}}
{(\lvert i_1\lvert^{p_1}+\ldots+i_d^{p_d}+1)
(\lvert i_1\lvert^{p_1}+\ldots+i_d'^{p_d}+1)\vert x-y\vert^{\alpha}}.
\end{equation}
We will split the general case into four ones:

\vsp\vsp

\noindent (a) $i_d'\leq 2i_d+1$,

\vsp\vsp

\noindent (b) $i_d'^{p_d}\leq\lvert i_1\lvert^{p_1}+\ldots+\lvert i_{d-1}\lvert^{p_{d-1}}$,

\vsp\vsp

\noindent (c) $i_d'\geq 2i_d+2$ \esp \esp and \esp\esp 
$i_d^{p_d}\geq\lvert i_1\lvert^{p_1}+\ldots+\lvert i_{d-1}\lvert^{p_{d-1}}$,

\vsp\vsp

\noindent (d) $i_d'\geq 2i_d+2$ \esp \esp and \esp\esp 
$i_d^{p_d} \leq \lvert i_1\lvert^{p_1}+\ldots+\lvert i_{d-1}\lvert^{p_{d-1}} \leq i_d'^{p_d}$.\\

\vsp

In case (a), the estimate \esp \esp
$\lvert x-y \lvert \hspace{0.2cm} \geq \esp\esp (i_d'-i_d-1)\lvert I_{i_1,i_2,\ldots,i_d'}\lvert$ 
\esp\esp shows that the expression (\ref{vvv}) is bounded from above by 
\begin{equation}\label{www}
M\frac{i_d'^{p_d-1}(i_d'-i_d)^{1-\alpha}(\lvert i_1\lvert+1)^{p_1-1}}
{(\lvert i_1\lvert^{p_1}+\ldots+i_d^{p_d}+1)(\lvert i_1\lvert^{p_1}+\ldots+i_d'^{p_d}+1)^{1-\alpha}}.
\end{equation}
By the condition $i_d' \!\leq\! 2i_d+1$, the latter expression is smaller than or equal to 
$$M\frac{i_d^{p_d-\alpha}(\lvert i_1\lvert+1)^{p_1-1}}{(\lvert i_1\lvert^{p_1}+i_d^{p_d}+1)^{2-\alpha}}.$$
If $\lvert i_1\lvert^{p_1}\leq i_d^{p_d}$, then $\frac{i_d^{p_d-\alpha}(\lvert i_1\lvert+1)^{p_1-1}}
{(\lvert i_1\lvert^{p_1}+i_d^{p_d}+1)^{2-\alpha}}\leq\frac{i_d^{p_d-\alpha}(i_d^{\frac{p_d}{p_1}}+1)^{p_1-1}}
{(i_d^{p_d}+1)^{2-\alpha}}$, and the last expression 
is uniformly bounded by condition $\mathrm{(iii_B)}$. 
If $i_d^{p_d} \leq \lvert i_1\lvert^{p_1}$, then 
$\frac{i_d^{p_d-\alpha}(\lvert i_1\lvert+1)^{p_1-1}}
{(\lvert i_1\lvert^{p_1}+i_d^{p_d}+1)^{2-\alpha}} \leq 
\frac{\lvert i_1\lvert^{\frac{p_1}{p_d}(p_d-\alpha)}(\lvert i_1\lvert+1)^{p_1-1}}
{(\lvert i_1\lvert^{p_1}+1)^{2-\alpha}}$, and this is uniformly bounded again 
by condition $\mathrm{(iii_B)}$.

\vsp

In case (b), the expression (\ref{vvv}) is still bounded from above by 
(\ref{www}), which in its turn is smaller than or equal to 
$$M\frac{i_d'^{p_d-\alpha}(\lvert i_1\lvert+1)^{p_1-1}}
{(\lvert i_1\lvert^{p_1}+\ldots+\lvert i_{d-1}\lvert^{p_{d-1}}+1)^{2-\alpha}}.$$
Now using the condition \esp 
\esp $i_d'^{p_d}\leq\lvert i_1\lvert^{p_1}+\ldots+\lvert i_{d-1}\lvert^{p_{d-1}}$, 
\esp \esp we see that this last expression is bounded from above by
$$M \frac{(\lvert i_1\lvert+1)^{p_1-1}}{(\lvert i_1\lvert^{p_1}
+\ldots+\lvert i_{d-1}\lvert^{p_{d-1}}+1)^{1-\alpha+\frac{\alpha}{p_d}}} \esp\esp \leq \esp \esp
\frac{(\lvert i_1\lvert+1)^{p_1-1}}{(\lvert i_1\lvert^{p_1}+1)^{1-\alpha+\frac{\alpha}{p_d}}}.$$
Finally, the right-side expression is uniformly bounded by condition $\mathrm{(iii_B)}$.\\

\vsp

In case (c), we first need to estimate the value of $\lvert x-y\lvert$:  
\begin{multline*}
\vert x-y \vert \esp\esp 
\geq \sum_{i_d< j <i_d'} \lvert I_{i_1,\ldots,i_{d-1},j} \vert 
= \sum_{i_d < j < i_d'}\frac{1}{\lvert i_1\lvert^{p_1}+\ldots+\lvert i_{d-1}\lvert^{p_{d-1}}+j^{p_d}+1} 
\geq \\ \geq 
\sum_{i_d<j<i_d'}\frac{1}{i_d^{p_d}+j^{p_d}+1}\geq\sum_{i_d<j<i_d'}\frac{1}{3j^{p_d}} 
\geq \int_{i_d+1}^{i_d'}\frac{1}{3x^{p_d}}dx 
\geq \\ \geq \frac{M}{(i_d+1)^{p_d-1}} 
\left(1-\left(\frac{i_d+1}{i_d'}\right)^{p_d-1}\right) 
\geq \\ \geq \frac{M}{(i_d+1)^{p_d-1}}\left(1-\left(\frac{1}{2}\right)^{p_d-1}\right)
\geq \frac{M}{(i_d+1)^{p_d-1}},
\end{multline*}
where in the second inequality we used the hypothesis \esp \esp  
$i_d^{p_d} \geq \lvert i_1\lvert^{p_1}+\ldots+\lvert i_{d-1}\lvert^{p_{d-1}}$. \esp 
\esp Using this, the value of (\ref{vvv}) is easily seen to be smaller than or equal to 
$$M\frac{i_d'^{p_d-1}(i_d'-i_d)(\lvert i_1\lvert+1)^{p_1-1}(i_d+1)^{(p_d-1)\alpha}}
{(\lvert i_1\lvert^{p_1}+\ldots+i_d^{p_d}+1)(\lvert i_1\lvert^{p_1}+\ldots+i_d'^{p_d}+1)}
\esp\esp\esp \leq \esp\esp\esp  
M\frac{(\lvert i_1\lvert+1)^{p_1-1}(i_d+1)^{(p_d-1)\alpha}}
{\lvert i_1\lvert^{p_1}+\ldots+i_d^{p_d}+1}.$$
Since by hypothesis we have \esp $i_d^{p_d} \geq \lvert i_1\lvert^{p_1}$, 
\esp the right-side expression above is bounded from above by
$$M\frac{(i_d^{\frac{p_d}{p_1}}+1)^{p_1-1}(i_d+1)^{(p_d-1)\alpha}}{i_d^{p_d}+1},$$
which is uniformly bounded by the condition $\mathrm{(iii_B)}$.\\

\vsp

Let us finally consider the case (d). Letting 
$$S := 1 + |i_1|^{p_1} + |i_2|^{p_2} + \ldots + |i_{d-1}|^{p_{d-1}},$$ 
we first observe that 
$$\vert x-y \vert \geq \sum_{i_d<j<i_d'} \vert I_{i_1,\ldots,i_{d-1},j}\vert 
= \sum_{i_d<j<i_d'}\frac{1}{S + j^{p_d}} 
\geq \int_{i_d+1}^{i_d'}\frac{dx}{x^{p_d} + S}
\geq \int_{i_d+1}^{i_d'}\frac{dx}{(x+ S^{1/p_d})^{p_d}}.$$
The last integral equals  
$$\frac{1}{(p_d-1)} \! \left[ \! \frac{1}{(i_d \!+\! 1 \!+\! S^{1/p_d})^{p_d - 1}} 
- \frac{1}{(i_d' \!+\! S^{1/p_d})^{p_d - 1}} \! \right] \!=\! \frac{1}{(p_d - 1)} \! 
\left[ \! \frac{(i_d' + S^{1/p_d})^{p_d - 1} \!- (i_d + 1 + S^{1/p_d})^{p_d - 1}}
{(i_d \!+\! 1 \!+\! S^{1/p_d})^{p_d - 1}(i_d' \!+\! S^{1/p_d})^{p_d - 1}} \! \right]\!\!.$$
Using the Mean Value Theorem, we conclude that
$$\vert x-y \vert \esp\geq\esp 
\frac{i_d'-i_d-1}{(i_d + 1 + S^{1/p_d})^{p_d - 1} \esp (i'_d + S^{1/p_d})}.$$ 
Using this, we conclude that (\ref{vvv}) is smaller than or equal to 
\begin{multline}\label{rrr}
\frac{i_d'^{p_d-1}(i_d'-i_d)(\lvert i_1\lvert+1)^{p_1-1}(i_d+1+S^{1/p_d})^{\alpha(p_d-1)}
(i_d'+S^{1/p_d})^{\alpha}}{(\lvert i_1\lvert^{p_1}+\ldots+i_d^{p_d}+1)
(\lvert i_1\lvert^{p_1}+\ldots+i_d'^{p_d}+1)(i_d'-i_d-1)^{\alpha}} 
\leq \\ \leq 
\frac{(\lvert i_1\lvert+1)^{p_1-1}(i_d+1+S^{1/p_d})^{\alpha(p_d-1)}
(i_d'+S^{1/p_d})^{\alpha}}{(\lvert i_1\lvert^{p_1}+\ldots+i_d^{p_d}+1)
(i_d'-i_d-1)^{\alpha}}.
\end{multline}
By hypothesis, \esp\esp $1 + i_d^{p_d} \leq S $, thus \esp $i_d \leq S^{1/p_d}$. \esp 
Since (by definition) \esp $|i_1|^{p_1} \leq S$, \esp this yields
\begin{equation}\label{a-a}
\frac{(\lvert i_1\lvert+1)^{p_1-1}(i_d+1+S^{1/p_d})^{\alpha(p_d-1)}}
{(\lvert i_1\lvert^{p_1}+\ldots+i_d^{p_d}+1)} 
\leq M S^{^{\frac{p_1 - 1}{p_1} + \frac{\alpha(p_d-1)}{p_d} - 1}}.
\end{equation}
By hypothesis, we also have \esp $S \leq 1 + i_d'^{p_d}$ \esp 
and \esp $i_d'\geq 2i_d+2$, \esp which gives 
\begin{equation}\label{b-b}
\frac{(i_d'+S^{1/p_d})^{\alpha}}{(i_d'-i_d-1)^{\alpha}} \leq M.
\end{equation}
Putting together (\ref{a-a}) and (\ref{b-b}), and using again that 
\esp \esp $S \leq 1 + i_d'^{p_d}$, \esp we conclude that the 
expression in (\ref{rrr}) is bounded from above by 
$$M S^{^{\frac{p_1 - 1}{p_1} + \frac{\alpha(p_d-1)}{p_d} - 1}},$$
which is uniformly bounded by the condition $\mathrm{(iii_B)}$.\\

\vsp\vsp

\noindent{\bf III.} Finally, we consider $x,y$ so that $x \in I_{i_1,\ldots,i_{d-1},i_d}$ and 
$y \in I_{i_1',\ldots,i_{d-1}',i_d'}$ for $(i_1,\ldots,i_{d-1}) \prec (i_1',\ldots,i_{d-1}')$, 
where $\prec$ stands for the lexicographic ordering. Letting $z$ (resp. $z'$) be the right 
endpoint (resp. left endpoint) of $\bigcup_{j \in \mathbb{Z}} I_{i_1,\ldots,i_{d-1},j}$ (resp. 
$\bigcup_{j \in \mathbb{Z}} I_{i_1',\ldots,i_{d-1}',j}$), by construction we have \esp 
$D f_1 (z) = D f_1 (z') = 1$. \esp Hence
$$\frac{\vert \log Df_1(x) - \log Df_1(y) \vert}{|x-y|^{\alpha}} \leq 
\frac{\vert \log Df_1(x) - \log Df_1(z) \vert}{|x-z|^{\alpha}} + 
\frac{\vert \log Df_1(z') - \log Df_1(y) \vert}{|z'-y|^{\alpha}},$$
and both terms of the sum above are uniformly bounded by the previous case. 

\vsp 

To conclude the proof of the regularity of $f_1$, notice that 
slightly modified arguments apply to the inverse $f_1^{-1}$, 
thus showing that $f_1$ is a $C^{1+\alpha}$-diffeomorphism.


\subsection{For $2 \leq j \leq d-1$, the map $f_j$ is a $C^{1+\alpha}$-diffeomorphism}

{\bf I.} We first consider $x,y$ in the same interval $I_{i_1,\ldots,i_d}$. We have
$$\frac{\vert \log Df_j(x) - \log Df_j(y) \vert}{\vert x-y\vert}\leq
\frac{M}{\vert I_{i_1,\ldots,i_d}\vert}\left\vert\frac{\vert I_{i_1,\ldots,i_d}\vert}
{\vert I_{i_1,\ldots,i_j+i_{j-1},\ldots,i_d}\vert}
\frac{\vert I_{i_1,\ldots,i_j+i_{j-1},\ldots,i_d-1}\vert}{\vert I_{i_1,\ldots,i_d-1}\vert}-1\right\vert.$$
Hence,
$$\frac{\vert \log Df_j(x) - \log Df_j(y) 
\vert}{\vert x-y\vert^{\alpha}}\leq\ M\left\vert\frac{\vert I_{i_1,\ldots,i_d}\vert}
{\vert I_{i_1,\ldots,i_j+i_{j-1},\ldots,i_d}\vert}
\frac{\vert I_{i_1,\ldots,i_j+i_{j-1},\ldots,i_d-1}\vert}
{\vert I_{i_1,\ldots,i_d-1}\vert}-1\right\vert\vert I_{i_1,\ldots,i_d}\vert^{-\alpha}.$$
One readily checks that the right-side expression equals
$$M\left\lvert\frac{(\lvert i_d\lvert^{p_d}-\lvert i_d-1\lvert^{p_d})
(\lvert i_j\lvert^{p_j}-\lvert i_j+i_{j-1}\lvert^{p_j})}
{(\lvert i_1\lvert^{p_1}+\ldots+\lvert i_j+i_{j-1}\lvert^{p_j}+\ldots+\lvert i_d-1\lvert^{p_d}+1)
(\lvert i_1\lvert^{p_1}+\ldots+\lvert i_j\lvert^{p_j}+\ldots+\lvert i_d\lvert^{p_d}+1)^{1-\alpha}}\right\lvert.$$
By The Mean Value Theorem, and since \esp $p_{j-1} \!>\! p_j$, 
\esp the last expression is bounded from above by  
\begin{equation}\label{roro}
M\frac{(\lvert i_d\lvert+1)^{p_d-1}(\lvert i_j\lvert+\lvert i_{j-1}\lvert)^{p_j-1}\lvert i_{j-1}\lvert}
{(\lvert i_1\lvert^{p_1}+\ldots+\lvert i_j\lvert^{p_j}+\ldots+\lvert i_d\lvert^{p_d}+1)^{2-\alpha}}.
\end{equation}
To estimate this expression, let us first assume that \esp 
$\lvert i_j\lvert^{p_j}\leq\lvert i_{j-1}\lvert^{p_{j-1}}$. \esp  
In this case, (\ref{roro}) is bounded from above by 
\begin{equation}\label{uyuy}
M\frac{(\lvert i_d\lvert+1)^{p_d-1}
(\lvert i_{j-1}\lvert^{\frac{p_{j-1}}{p_j}}+\lvert i_{j-1}\lvert)^{p_j-1}\lvert i_{j-1}\lvert}
{(\lvert i_{j-1}\lvert^{p_{j-1}}+\lvert i_d\lvert^{p_d}+1)^{2-\alpha}}.
\end{equation} 
If \esp $\lvert i_d\lvert^{p_d}\leq\lvert i_{j-1}\lvert^{p_{j-1}}$, \esp 
then this expression is smaller than or equal to
$$M\frac{(\lvert i_{j-1}\lvert^{\frac{p_{j-1}}{p_d}}+1)^{p_d-1}
(\lvert i_{j-1}\lvert^{\frac{p_{j-1}}{p_j}}+\lvert i_{j-1}\lvert)^{p_j-1}
\lvert i_{j-1}\lvert}{(\lvert i_{j-1}\lvert^{p_{j-1}}+1)^{2-\alpha}}.$$ 
Since \esp $p_{j-1} / p_j \geq 1$, \esp this is uniformly bounded if 
$$\frac{p_{j-1}}{p_d} (p_d - 1) + \frac{p_{j-1}}{p_j} (p_j - 1) 
+ 1 - p_{j-1} (2 - \alpha) \hspace{0.15cm} \leq \hspace{0.15cm} 0,$$
that is, \esp $\alpha\leq\frac{1}{p_d}+\frac{1}{p_j}-\frac{1}{p_{j-1}}$, 
\esp and this is ensured by conditions $\mathrm{(i_B)}$ and $\mathrm{(v_B)}$. 
\esp\esp If \esp $\lvert i_{j-1}\lvert^{p_{j-1}}\leq\lvert i_d\lvert^{p_d}$, 
\esp then the expression (\ref{uyuy}) is smaller than or equal to 
$$M\frac{(\lvert i_d\lvert+1)^{p_d-1}
(\lvert i_{d}\lvert^{\frac{p_{d}}{p_j}}+\lvert i_{d}
\lvert^{\frac{p_d}{p_{j-1}}})^{p_j-1}\lvert i_{d}\lvert^{\frac{p_d}{p_{j-1}}}}
{(\lvert i_d\lvert^{p_d}+1)^{2-\alpha}}.$$ 
Since \esp $p_d / p_{j-1} \leq p_d / p_j$, \esp this is uniformly bounded if 
$$p_d - 1 + \frac{p_d}{p_j} (p_j - 1) + \frac{p_d}{p_{j-1}} 
- p_d (2 - \alpha) \hspace{0.15cm} \leq \hspace{0.15cm} 0,$$
which is again ensured by conditions $\mathrm{(i_B)}$ and  $\mathrm{(v_B)}$.

Assume now that \esp\esp $\lvert i_{j-1}\lvert^{p_{j-1}}\leq\lvert i_{j}\lvert^{p_{j}}$. \esp\esp  
In this case, (\ref{roro}) is bounded from above by
$$M\frac{(\lvert i_d\lvert+1)^{p_d-1}
(\lvert i_j\lvert+\lvert i_{j}\lvert^{\frac{p_j}{p_{j-1}}})^{p_j-1}
\lvert i_{j}\lvert^{\frac{p_j}{p_{j-1}}}}{(\lvert i_j\lvert^{p_j}
+\lvert i_d\lvert^{p_d}+1)^{2-\alpha}}.$$
Proceeding as in the previous case, one readily checks that this expression is 
uniformly bounded when \esp $\alpha\leq\frac{1}{p_d}+\frac{1}{p_j}-\frac{1}{p_{j-1}}$, 
\esp which is ensured by conditions $\mathrm{(i_B)}$ and  $\mathrm{(v_B)}$.\\

\vsp\vsp

{\bf II.} Now we consider the case where $x \!\in\! I_{i_1,i_2,\ldots,i_d}$ and 
$y \!\in\! I_{i_1,i_2,\ldots,i_d'}$ for different $i_d,i_d'$. As in the case of $f_1$, 
we may restrict ourselves to the case where \esp $i_d'\!-i_d \geq 2$, \esp with $i_d,i_d'$ 
both positive.

Once again, property (\ref{dos-tsuboi}) implies that \esp 
$\vert \log Df_j(x) - \log Df_j(y) \vert$ \esp 
is smaller than or equal to the sum 
\begin{multline*}
\left\lvert\log\frac{\lvert I_{i_1,\ldots,i_j+i_{j-1},\ldots,i_d}\lvert}
{\lvert I_{i_1,\ldots,i_j,\ldots,i_d}\lvert}-
\log\frac{\lvert I_{i_1,\ldots,i_j+i_{j-1},\ldots,i_d'}\lvert}
{\lvert I_{i_1,\ldots,i_j,\ldots,i_d'}\lvert}\right\lvert
+ \\ + 
\left\lvert\log\frac{\lvert I_{i_1,\ldots,i_j+i_{j-1},\ldots,i_d}\lvert}
{\lvert I_{i_1,\ldots,i_j,\ldots,i_d}\lvert}-
\log\frac{\lvert I_{i_1,\ldots,i_j+i_{j-1},\ldots,i_d-1}\lvert}
{\lvert I_{i_1,\ldots,i_j,\ldots,i_d-1}\lvert}\right\lvert 
+ \\ +  
\left\lvert\log\frac{\lvert I_{i_1,\ldots,i_j+i_{j-1},\ldots,i_d'-1}\lvert}
{\lvert I_{i_1,\ldots,i_j,\ldots,i_d'-1}\lvert}-
\log\frac{\lvert I_{i_1,\ldots,i_j+i_{j-1},\ldots,i_d'}\lvert}
{\lvert I_{i_1,\ldots,i_j,\ldots,i_d'}\lvert}\right\lvert.
\end{multline*}As in previous estimates of similar expressions, we have 
$$\vert \log Df_j(x) - \log Df_j(y) \vert \leq 3 
\left\lvert\log \left( \frac{\lvert I_{i_1,\ldots,i_j+i_{j-1},\ldots,i_d-1}\lvert}
{\lvert I_{i_1,\ldots,i_j,\ldots,i_d-1}\lvert} \right) 
- \log \left( \frac{\lvert I_{i_1,\ldots,i_j + i_{j-1},\ldots,i_d'}\lvert}
{\lvert I_{i_1,\ldots,i_j,\ldots,i_d'}\lvert} \right) \right\lvert.$$
The last expression equals 
$$3 \left\lvert\log\frac{(\lvert i_1\lvert^{p_1}+\ldots+\lvert i_j\lvert^{p_j}
+\ldots+\lvert i_d-1\lvert^{p_d}+1)(\lvert i_1\lvert^{p_1}+\ldots+\lvert i_j+i_{j-1}\lvert^{p_j}
+\ldots+\lvert i_d'\lvert^{p_d}+1)}{(\lvert i_1\lvert^{p_1}+\ldots+\lvert i_j+i_{j-1}\lvert^{p_j}
+\ldots+\lvert i_d-1\lvert^{p_d}+1)(\lvert i_1\lvert^{p_1}+\ldots+\lvert i_j\lvert^{p_j}+\ldots
+\lvert i_d'\lvert^{p_d}+1)}\right\lvert,$$
that is,
\begin{small}\begin{equation}\label{call-1}
3 \left\lvert\log\left(1+\frac{(|i_d-1|^{p_d}-i_d'^{p_d})
(\lvert i_j+i_{j-1}\lvert^{p_j}-\lvert i_j\lvert^{p_j})}
{(\lvert i_1\lvert^{p_1}+\ldots+\lvert i_j+i_{j-1}\lvert^{p_j}+\ldots+\lvert i_d-1 \lvert^{p_d}+1)
(\lvert i_1\lvert^{p_1}+\ldots+\lvert i_j\lvert^{p_j}+\ldots+\lvert i_d'\lvert^{p_d}+1)}\right)\right\lvert.
\end{equation}
\end{small}The expression into brackets in the right-side term equals 
$$\frac{(\lvert i_1\lvert^{p_1}+\ldots+\lvert i_j\lvert^{p_j}+\ldots+\lvert i_d-1\lvert^{p_d}+1)
(\lvert i_1\lvert^{p_1}+\ldots+\lvert i_j+i_{j-1}\lvert^{p_j}+\ldots+\lvert i_d'\lvert^{p_d}+1)}
{(\lvert i_1\lvert^{p_1}+\ldots+\lvert i_j+i_{j-1}\lvert^{p_j}+\ldots+\lvert i_d-1\lvert^{p_d}+1)
(\lvert i_1\lvert^{p_1}+\ldots+\lvert i_j\lvert^{p_j}+\ldots+\lvert i_d'\lvert^{p_d}+1)},$$ 
hence it is uniformly bounded from below by a positive number. 
Therefore, the value of (\ref{call-1}) is smaller than or equal to  
$$M\left\lvert\frac{(i_d^{p_d}-i_d'^{p_d})
(\lvert i_j+i_{j-1}\lvert^{p_j}-\lvert i_j\lvert^{p_j})}
{(\lvert i_1\lvert^{p_1}+\ldots+\lvert i_j+i_{j-1}\lvert^{p_j}+\ldots+\lvert i_d\lvert^{p_d}+1)
(\lvert i_1\lvert^{p_1}+\ldots+\lvert i_j\lvert^{p_j}+\ldots+\lvert i_d'\lvert^{p_d}+1)}\right\lvert.$$
Using the Mean Value Theorem and the condition \esp $p_{j-1} \!>\! p_j$, \esp this last 
expression is easily seen to be bounded from above by 
$$M\frac{i_d'^{p_d-1}(i_d'-i_d)
(\lvert i_j\lvert+\lvert i_{j-1}\lvert)^{p_j-1}\lvert i_{j-1}\lvert}
{(\lvert i_1\lvert^{p_1}+\ldots+\lvert i_j\lvert^{p_j}+\ldots+i_d^{p_d}+1)
(\lvert i_1\lvert^{p_1}+\ldots+\lvert i_j\lvert^{p_j}+\ldots+i_d'^{p_d}+1)}.$$
Therefore, \esp\esp 
$\frac{\vert \log Df_j(x) - \log Df_j(y) \vert}{\vert x-y\vert^{\alpha}}$ 
\esp \esp is smaller than or equal to  
\begin{equation}\label{call-2}
M \frac{i_d'^{p_d-1}(i_d'-i_d)(\lvert i_j\lvert+\lvert i_{j-1}\lvert)^{p_j-1}\lvert i_{j-1}\lvert}
{(\lvert i_1\lvert^{p_1}+\ldots+\lvert i_j\lvert^{p_j}+\ldots+i_d^{p_d}+1)
(\lvert i_1\lvert^{p_1}+\ldots+\lvert i_j\lvert^{p_j}+\ldots+i_d'^{p_d}+1)\vert x-y\vert^{\alpha}}.
\end{equation}
In order to estimate this expression, we will again consider 
separately the cases (a), (b), (c) and (d) of the previous section.\\

The case (a) is \esp $i_d'\leq 2i_d+1$. \esp 
Here the estimate $\lvert x-y\lvert\geq(i_d'-i_d-1)\lvert I_{i_1,i_2,\ldots,i_d'}\lvert$ 
shows that (\ref{call-2}) is bounded from above by 
\begin{equation}\label{corro}
M\frac{i_d'^{p_d-1}(i_d'-i_d)^{1-\alpha}(\lvert i_j\lvert+\lvert i_{j-1}\lvert)^{p_j-1}\lvert i_{j-1}\lvert}
{(\lvert i_1\lvert^{p_1}+\ldots+\lvert i_j\lvert^{p_j}+\ldots+i_d^{p_d}+1)
(\lvert i_1\lvert^{p_1}+\ldots+\lvert i_j\lvert^{p_j}+\ldots+i_d'^{p_d}+1)^{1-\alpha}},
\end{equation}
which is smaller than or equal to 
\begin{equation}\label{cras}
M\frac{i_d^{p_d-\alpha}(\lvert i_j\lvert+\lvert i_{j-1}\lvert)^{p_j-1}\lvert i_{j-1}\lvert}
{(\lvert i_1\lvert^{p_1}+\ldots+\lvert i_j\lvert^{p_j}+\ldots+i_d^{p_d}+1)^{2-\alpha}}.
\end{equation}
There are three subcases:

\vsp 

\noindent -- If $\lvert i_j\lvert^{p_j}\leq\lvert i_{j-1}\lvert^{p_{j-1}}$ 
and $i_d^{p_d}\leq\lvert i_{j-1}\lvert^{p_{j-1}}$, then 
$$M \frac{i_d^{p_d-\alpha}(\lvert i_j\lvert+\lvert i_{j-1}\lvert)^{p_j-1}\lvert i_{j-1}\lvert}
{(\lvert i_1\lvert^{p_1}+\ldots+\lvert i_j\lvert^{p_j}+\ldots+i_d^{p_d}+1)^{2-\alpha}} 
\leq M \frac{\lvert i_{j-1}\lvert^{\frac{p_{j-1}}{p_d}(p_d-\alpha)}(\lvert i_{j-1}
\lvert^{\frac{p_{j-1}}{p_j}}+\lvert i_{j-1}\lvert)^{p_j-1}\lvert i_{j-1}\lvert}
{(\lvert i_{j-1}\lvert^{p_{j-1}}+1)^{2-\alpha}}.$$
The last expression is easily seen to be uniformly bounded by condition $\mathrm{(iv_B)}$.

\vsp

\noindent -- If $\lvert i_j\lvert^{p_j}\leq i_d^{p_d}$ 
and $\lvert i_{j-1}\lvert^{p_{j-1}} \leq i_d^{p_d}$, then 
$$M \frac{i_d^{p_d-\alpha}(\lvert i_j\lvert+\lvert i_{j-1}\lvert)^{p_j-1}\lvert i_{j-1}\lvert}
{(\lvert i_1\lvert^{p_1}+\ldots+\lvert i_j\lvert^{p_j}+\ldots+i_d^{p_d}+1)^{2-\alpha}} 
\leq M\frac{i_d^{p_d-\alpha}(i_d^{\frac{p_d}{p_j}}+i_d^{\frac{p_d}
{p_{j-1}}})^{p_j-1}i_d^{\frac{p_d}{p_{j-1}}}}{(i_d^{p_d}+1)^{2-\alpha}},$$
and the last expression is uniformly bounded by condition $\mathrm{(iv_B)}$.

\vsp

\noindent -- If $\lvert i_{j-1}\lvert^{p_{j-1}}\leq\lvert i_j\lvert^{p_j}$ 
and $i_d^{p_d}\leq\lvert i_j\lvert^{p_j}$, then one has 
$$M \frac{i_d^{p_d-\alpha}(\lvert i_j\lvert+\lvert i_{j-1}\lvert)^{p_j-1}\lvert i_{j-1}\lvert}
{(\lvert i_1\lvert^{p_1}+\ldots+\lvert i_j\lvert^{p_j}+\ldots+i_d^{p_d}+1)^{2-\alpha}} 
\leq M \frac{\lvert i_j\lvert^{\frac{p_j}{p_d}(p_d-\alpha)}
(\lvert i_j\lvert+\lvert i_j\lvert^{\frac{p_j}{p_{j-1}}})^{p_j-1}
\lvert i_j\lvert^{\frac{p_j}{p_{j-1}}}}{(\lvert i_j\lvert^{p_j}+1)^{2-\alpha}},$$
and the last expression is uniformly bounded by condition $\mathrm{(iv_B)}$.\\

\vsp

In case (b), we still have the upper bound (\ref{corro}) for (\ref{call-2}). Now, using the 
condition 
$$i_d'^{p_d}\leq\lvert i_1\lvert^{p_1}+\ldots+\lvert i_{d-1}\lvert^{p_{d-1}},$$ 
\esp the value of (\ref{corro}) is easily seen to be bounded from above by 
$$M\frac{i_d'^{p_d-\alpha}(\lvert i_j\lvert+
\lvert i_{j-1}\lvert)^{p_j-1}\lvert i_{j-1}\lvert}{(\lvert i_1\lvert^{p_1}+\ldots+
\lvert i_{d-1}\lvert^{p_{d-1}}+1)^{2-\alpha}} \esp\esp\esp \leq \esp\esp\esp
M\frac{(\lvert i_j\lvert+\lvert i_{j-1}\lvert)^{p_j-1}\lvert i_{j-1}\lvert}
{(\lvert i_1\lvert^{p_1}+\ldots+\lvert i_{d-1}\lvert^{p_{d-1}}+1)^{1-\alpha+\frac{\alpha}{p_d}}}.$$
To estimate the right-side expression of this inequality, we consider two subcases:

\vsp

\noindent -- If $\lvert i_{j-1}\lvert^{p_{j-1}}\leq\lvert i_j\lvert^{p_j}$, 
then 
$$M \frac{(\lvert i_j\lvert+\lvert i_{j-1}\lvert)^{p_j-1}\lvert i_{j-1}\lvert}
{(\lvert i_1\lvert^{p_1}+\ldots+\lvert i_{d-1}\lvert^{p_{d-1}}+1)^{1-\alpha+\frac{\alpha}{p_d}}} 
\leq M \frac{(\lvert i_j\lvert+\lvert i_j\lvert^{\frac{p_j}{p_{j-1}}})^{p_j-1}
\lvert i_j\lvert^{\frac{p_j}{p_{j-1}}}}{(\lvert i_j\lvert^{p_j}+1)^{1-\alpha+\frac{\alpha}{p_d}}},$$ 
and the last expression is easily seen to be uniformly bounded by condition $\mathrm{(iv_B)}$.

\vsp

\noindent -- If $\lvert i_j\lvert^{p_j}\leq\lvert i_{j-1}\lvert^{p_{j-1}}$, 
then 
$$M \frac{(\lvert i_j\lvert+\lvert i_{j-1}\lvert)^{p_j-1}\lvert i_{j-1}\lvert}
{(\lvert i_1\lvert^{p_1}+\ldots+\lvert i_{d-1}\lvert^{p_{d-1}}+1)^{1-\alpha+\frac{\alpha}{p_d}}} 
\leq M \frac{(\lvert i_{j-1}\lvert^{\frac{p_{j-1}}{p_j}}+
\lvert i_{j-1}\lvert)^{p_j-1}\lvert i_{j-1}\lvert}{(\lvert i_{j-1}
\lvert^{p_{j-1}}+1)^{1-\alpha+\frac{\alpha}{p_d}}},$$ 
and the last expression is easily seen to be uniformly bounded by condition $\mathrm{(iv_B)}$.\\

\vsp

In case (c), we had the estimate 
\begin{equation}\label{la-de-c}
\vert x-y\vert \esp\esp\esp \geq \esp\esp\esp \frac{M}{(i_d+1)^{p_d-1}},
\end{equation}
which shows that (\ref{call-2}) is bounded from above by 
$$M\frac{i_d'^{p_d-1}(i_d'-i_d)(\lvert i_j\lvert+\lvert i_{j-1}\lvert)^{p_j-1}
\lvert i_{j-1}\lvert(i_d+1)^{(p_d-1)\alpha}}{(\lvert i_1\lvert^{p_1}+\ldots+
\lvert i_j\lvert^{p_j}+\ldots+i_d^{p_d}+1)(\lvert i_1\lvert^{p_1}+\ldots+
\lvert i_j\lvert^{p_j}+\ldots+i_d'^{p_d}+1)}.$$
This is smaller than or equal to   
\begin{equation}\label{uay}
M\frac{(\lvert i_j\lvert+\lvert i_{j-1}\lvert)^{p_j-1}
\lvert i_{j-1}\lvert(i_d+1)^{(p_d-1)\alpha}}
{\lvert i_1\lvert^{p_1}+\ldots+\lvert i_j\lvert^{p_j}+\ldots+i_d^{p_d}+1}.
\end{equation}
Now from the condition \hspace{0.15cm} 
$i_d^{p_d}\geq\lvert i_1\lvert^{p_1}+\ldots+\lvert i_{d-1}\lvert^{p_{d-1}}$ \hspace{0.15cm} 
it follows that \hspace{0.15cm} $\lvert i_j\lvert^{p_j} \leq i_d^{p_d}$ \hspace{0.15cm} 
and \hspace{0.15cm} $\lvert i_{j-1}\lvert^{p_{j-1}} \leq i_d^{p_d}$. \hspace{0.15cm} 
Therefore, (\ref{uay}) is bounded from above by
$$M\frac{(i_d^{p_d/p_j}+i_d^{p_d/p_{j-1}})^{p_j-1} 
i_d^{p_d/p_{j-1}}(i_d+1)^{(p_d-1)\alpha}}{i_d^{p_d}+1},$$
and this expression is easily seen to be uniformly 
bounded by condition $\mathrm{(iv_B)}$.\\

\vsp

In case (d), we had the estimate 
\begin{equation}\label{clef}
\vert x-y \vert \esp\geq\esp 
\frac{i_d'-i_d-1}{(i_d + 1 + S^{1/p_d})^{p_d-1} \esp (i'_d + S^{1/p_d})},
\end{equation} 
where \esp\esp $S := 1 + |i_1|^{p_1} + |i_2|^{p_2} \ldots |i_{d-1}|^{p_{d-1}}.$ 
\esp\esp Thus, (\ref{call-2}) is bounded from above by
$$M \frac{i_d'^{p_d-1}(i_d'-i_d)(\lvert i_j\lvert
+\lvert i_{j-1}\lvert)^{p_j-1}\lvert i_{j-1}\lvert 
(i_d+1+S^{1/p_d})^{\alpha(p_d-1)}(i_d'+S^{1/p_d})^{\alpha}}
{(\lvert i_1\lvert^{p_1}+\ldots+\lvert i_j\lvert^{p_j}+
\ldots+i_d^{p_d}+1)(\lvert i_1\lvert^{p_1}+\ldots+\lvert i_j\lvert^{p_j}+
\ldots+i_d'^{p_d}+1)(i_d'-i_d-1)^{\alpha}},$$
hence by  
$$M \frac{(\lvert i_j\lvert
+\lvert i_{j-1}\lvert)^{p_j-1}\lvert i_{j-1}\lvert 
(i_d+1+S^{1/p_d})^{\alpha(p_d-1)}(i_d'+S^{1/p_d})^{\alpha}}
{(\lvert i_1\lvert^{p_1}+\ldots+\lvert i_j\lvert^{p_j}+
\ldots+i_d^{p_d}+1)(i_d'-i_d-1)^{\alpha}}.$$
Since the condition \esp $1 + i_d^{p_d} \leq S $ \esp yields  
$i_d \leq S^{1/p_d}$, this expression is smaller than or equal to 
$$M \esp \esp \frac{(\lvert i_j\lvert + \lvert i_{j-1}\lvert)^{p_j-1} 
\lvert i_{j-1}\lvert (i_d'+S^{1/p_d})^{\alpha}} {(i_d'-i_d-1)^{\alpha}} 
\esp\esp S^{^{\frac{\alpha(p_d-1)}{p_d} - 1}}.$$
The conditions \esp $1 + i_d^{p_d} \leq S$ \esp and \esp $p_{j-1} \geq p_j$ \esp also 
yield \esp $|i_j| \leq S^{1/p_j}$ \esp and \esp $|i_{j-1}| \leq S^{1/p_{j-1}} \leq S^{1/p_j}$, 
\esp thus showing that the last expression  is smaller than or equal to 
$$M \esp \esp \frac{(i_d'+S^{1/p_d})^{\alpha}} {(i_d'-i_d-1)^{\alpha}} 
\esp\esp S^{^{\frac{p_j-1}{p_j} + \frac{1}{p_{j-1}} + \frac{\alpha(p_d-1)}{p_d} - 1}}.$$
Using the conditions \esp $i_d'\geq 2i_d+2$ \esp and $S \leq 1 + i_d'^{p_d}$, \esp this 
last expression is easily seen to be bounded from above by 
$$M S^{\frac{p_j-1}{p_j} + \frac{1}{p_{j-1}} + \frac{\alpha(p_d-1)}{p_d} - 1},$$
which is uniformly bounded by the condition $\mathrm{(iv_B)}$.

\vsp\vsp 

\noindent{\bf III.} Finally, in the case where $x \in I_{i_1,\ldots,i_d}$ 
and $y \in I_{i_1',\ldots,i_d'}$ for different $(i_1,\ldots,i_{d-1})$ and 
$(i_1',\ldots,i_{d-1}')$, one may apply the same argument as that of $f_1$. 


\subsection{The map $f_d$ is a $C^{1+\alpha}$-diffeomorphism}

{\bf I.} First we consider $x,y$ in the same interval $I_{i_1,\ldots,i_d}$. We have
$$\frac{\vert \log Df_d(x) - \log Df_d(y) \vert}{\vert x-y\vert} 
\leq \frac{M}{\vert I_{i_1,\ldots,i_d}\vert}\left\vert\frac{\vert I_{i_1,\ldots,i_d}\vert}
{\vert I_{i_1,\ldots,i_d+i_{d-1}}\vert}\frac{\vert I_{i_1,\ldots,i_d+i_{d-1}-1}\vert}
{\vert I_{i_1,\ldots,i_d-1}\vert}-1\right\vert,$$
hence
$$\frac{\vert \log Df_d(x) - \log Df_d(y) \vert}{\vert x-y\vert^{\alpha}}
\leq\ M\left\vert\frac{\vert I_{i_1,\ldots,i_d}\vert}
{\vert I_{i_1,\ldots,i_d+i_{d-1}}\vert}\frac{\vert I_{i_1,\ldots,i_d+i_{d-1}-1}\vert}
{\vert I_{i_1,\ldots,i_d-1}\vert}-1\right\vert\vert I_{i_1,\ldots,i_d}\vert^{-\alpha}.$$
The right-side term above is smaller than or equal to 
$$M\left\lvert\frac{(\lvert i_1\lvert^{p_1}+\ldots+\lvert i_d-1\lvert^{p_d}+1)
(\lvert i_1\lvert^{p_1}+\ldots+\lvert i_d+i_{d-1}\lvert^{p_d}+1)}
{(\lvert i_1\lvert^{p_1}+\ldots+\lvert i_d+i_{d-1}-1\lvert^{p_d}+1)
(\lvert i_1\lvert^{p_1}+\ldots+\lvert i_d\lvert^{p_d}+1)}-1\right\lvert
(\lvert i_1\lvert^{p_1}+\ldots+\lvert i_d\lvert^{p_d}+1)^{\alpha},$$
which equals
$$M \! \left\lvert\frac{\sum_{k=1}^{d-1}\lvert i_k\lvert^{p_k}
(\lvert i_d+i_{d-1}\lvert^{p_d}-\lvert i_d\lvert^{p_d}) \!+\! 
\sum_{k=1}^{d-1}\lvert i_k\lvert^{p_k}(\lvert i_d-1\lvert^{p_d}-\lvert i_d+i_{d-1}-1\lvert^{p_d})+C}
{(\lvert i_1\lvert^{p_1}+\ldots+\lvert i_d+i_{d-1}-1\lvert^{p_d}+1)
(\lvert i_1\lvert^{p_1}+\ldots+\lvert i_d\lvert^{p_d}+1)}\right\lvert
\! \big( S+\lvert i_d \lvert^{p_d} \big)^{\alpha},$$
where 
$$C \esp := \esp \vert i_d-1\vert^{p_d}\lvert i_d+i_{d-1}\lvert^{p_d}-\lvert 
i_d+i_{d-1}-1\lvert^{p_d}\lvert i_d\lvert^{p_d}+
\lvert i_d-1\lvert^{p_d}-\lvert i_d+i_{d-1}-1\lvert^{p_d}+
\lvert i_d+i_{d-1}\lvert^{p_d}-\lvert i_d\lvert^{p_d}$$
and, as before, \esp 
$S := 1 + |i_1|^{p_1} + |i_2|^{p_2} \ldots |i_{d-1}|^{p_{d-1}}$. \esp  
By the Mean Value Theorem, the last expression is bounded from above by 
$$M\frac{\sum_{k=1}^{d-1}\lvert i_k\lvert^{p_k}(\lvert i_d\lvert
+\lvert i_{d-1}\lvert)^{p_d-1}\lvert 
i_{d-1}\lvert+\sum_{k=1}^{d-1}\lvert i_k\lvert^{p_k}(\lvert i_d\lvert+\lvert 
i_{d-1}\lvert+1)^{p_d-1}\lvert i_{d-1}\lvert+C'}{(\lvert i_1\lvert^{p_1}+
\ldots+\lvert i_d+i_{d-1}-1\lvert^{p_d}+1)(\lvert i_1\lvert^{p_1}+\ldots+
\lvert i_d\lvert^{p_d}+1)^{1-\alpha}},$$
where 
\begin{multline*}
C'\esp\esp\esp\esp := 
\esp\esp\esp\esp \lvert i_d+i_{d-1}\lvert^{p_d}(\lvert i_d\lvert+1)^{p_d-1}+
\lvert i_d\lvert^{p_d}(\lvert i_d\lvert+\lvert i_{d-1}\lvert+1)^{p_d-1}+\\+
(\lvert i_d\lvert+\lvert i_{d-1}\lvert+1)^{p_d-1}\lvert i_{d-1}\lvert+
(\lvert i_d\lvert+\lvert i_{d-1}\lvert)^{p_d-1}\lvert i_{d-1}\lvert.
\end{multline*}
To get an upper bound for this last expression, it is enough to do so for 
\begin{equation}\label{l-2}
\frac{\lvert i_d+i_{d-1}\lvert^{p_d}(\lvert i_d\lvert+1)^{p_d-1}}{(\lvert i_1\lvert^{p_1}+
\ldots+\lvert i_d+i_{d-1}-1\lvert^{p_d}+1)(\lvert i_1\lvert^{p_1}+\ldots+
\lvert i_d\lvert^{p_d}+1)^{1-\alpha}}
\end{equation}
and 
\begin{equation}\label{l-1}
\frac{\lvert i_k\lvert^{p_k}(\lvert i_d\lvert+\lvert i_{d-1}\lvert)^{p_d-1}\lvert i_{d-1}\lvert}
{(\lvert i_1\lvert^{p_1}+\ldots+\lvert i_d+i_{d-1}-1\lvert^{p_d}+1)(\lvert i_1\lvert^{p_1}+
\ldots+\lvert i_d\lvert^{p_d}+1)^{1-\alpha}}, 
\end{equation}
where \esp $1 \leq k \leq n.$ 

\vsp

Expression (\ref{l-2}) may be written as 
$$\frac{\lvert i_d+i_{d-1}\lvert^{p_d}}{(\lvert i_1\lvert^{p_1}+\ldots+\lvert i_d+i_{d-1}-1\lvert^{p_d}+1)} \esp 
\esp \frac{(\lvert i_d\lvert+1)^{p_d-1}}{(\lvert i_1\lvert^{p_1}+\ldots+\lvert i_d\lvert^{p_d}+1)^{1-\alpha}}.$$
The first factor is uniformly bounded, whereas the second is smaller than or equal to
$$\frac{(\lvert i_d\lvert+1)^{p_d-1}}{(\lvert i_d\lvert^{p_d}+1)^{1-\alpha}}.$$ 
This last expression is uniformly bounded provided that \esp 
$p_d - 1 - p_d (1-\alpha) \esp\esp\esp \leq \esp\esp\esp 0$, 
\esp which is a consequence of condition $\mathrm{(v_B)}$.

\vsp

Concerning expression (\ref{l-1}), notice that since $p_{d-1} \! > \! p_d$, 
it is smaller than or equal to 
$$\frac{\lvert i_k\lvert^{p_k}
(\lvert i_d\lvert+\lvert i_{d-1}\lvert)^{p_d-1}\lvert i_{d-1}\lvert}
{(\lvert i_1\lvert^{p_1}+\ldots+\lvert i_d\lvert^{p_d}+1)^{2-\alpha}} \hspace{0.2cm}
\leq \hspace{0.2cm} \frac{(\lvert i_d\lvert+\lvert i_{d-1}\lvert)^{p_d-1}\lvert i_{d-1}\lvert}
{(\lvert i_1\lvert^{p_1}+\ldots+\lvert i_d\lvert^{p_d}+1)^{1-\alpha}}.$$ 
On the one hand, if \esp $\lvert i_d\lvert^{p_d}\leq\lvert i_{d-1}\lvert^{p_{d-1}}$, 
\esp then 
$$\frac{(\lvert i_d\lvert+\lvert i_{d-1}\lvert)^{p_d-1}\lvert i_{d-1}\lvert}
{(\lvert i_1\lvert^{p_1}+\ldots+\lvert i_d\lvert^{p_d}+1)^{1-\alpha}} \esp\esp 
\leq \esp\esp \frac{(\lvert i_{d-1}\lvert^{\frac{p_{d-1}}{p_d}}+
\lvert i_{d-1}\lvert)^{p_d-1}\lvert i_{d-1}\lvert}{(\lvert i_{d-1}\lvert^{p_{d-1}}+1)^{1-\alpha}},$$ 
and the last term is uniformly bounded when \esp
$\frac{p_{d-1}}{p_d}(p_d-1)+1\leq p_{d-1}(1-\alpha)$, 
\esp which is ensured by condition $\mathrm{(v_B)}$. On the other hand, if 
\esp $\lvert i_{d-1}\lvert^{p_{d-1}}\leq\lvert i_d\lvert^{p_d}$, \esp then 
$$\frac{(\lvert i_d\lvert+\lvert i_{d-1}\lvert)^{p_d-1}\lvert i_{d-1}\lvert}
{(\lvert i_1\lvert^{p_1}+\ldots+\lvert i_d\lvert^{p_d}+1)^{1-\alpha}} \esp\esp 
\leq \esp\esp \frac{(\lvert i_d\lvert+\lvert i_d\lvert^{\frac{p_d}{p_{d-1}}})^{p_d-1}
\lvert i_d\lvert^{\frac{p_d}{p_{d-1}}}}{(\lvert i_d\lvert^{p_d}+1)^{1-\alpha}},$$ 
which is uniformly bounded when \esp $p_d-1+\frac{p_d}{p_{d-1}}\leq p_d(1-\alpha)$, 
\esp that is, when condition $\mathrm{(v_B)}$ holds.\\ 

\vsp\vsp

\noindent{\bf II.} Next we consider the case where $x \in I_{i_1,i_2,\ldots,i_d}$ and 
$y \in I_{i_1,i_2,\ldots,i_d'}$, with $i_d'-i_d \geq 2$. (For the case where $i_d' = i_d + 1$, 
we apply a similar argument to that of the previous maps.) Once again, we will only deal with 
positive $i_d,i_d'$. As in previous cases, \esp $\vert \log Df_d (x) - \log Df_d (y) \vert$ \esp 
is bounded from above by 
\begin{small}
$$3 \left\lvert\log\frac{\lvert I_{i_1,\ldots,i_d+i_{d-1}-1}\lvert}
{\lvert I_{i_1,\ldots,i_d-1}\lvert}\frac{\lvert I_{i_1,\ldots,i_d'}\lvert}
{\lvert I_{i_1,\ldots,i_d'+i_{d-1}}\lvert}\right\lvert 
\leq M \left\lvert\log\frac{(\lvert i_1\lvert^{p_1}+\ldots+i_d^{p_d}+1)
(\lvert i_1\lvert^{p_1}+\ldots+\lvert i_d'+i_{d-1}\lvert^{p_d}+1)}
{(\lvert i_1\lvert^{p_1}+\ldots+\lvert i_d+i_{d-1}\lvert^{p_d}+1)
(\lvert i_1\lvert^{p_1}+\ldots+i_d'^{p_d}+1)}\right\lvert\!.$$
\end{small}Note that the right-side term may be rewritten as 
\begin{equation}\label{a-r}
M \esp \left\lvert\log\left(1+\frac{\sum_{k=1}^{d-1}\lvert i_k\lvert^{p_k}
(\lvert i_d' +i_{d-1}\lvert^{p_d}-i_d'^{p_d})+\sum_{k=1}^{d-1}
\lvert i_k\lvert^{p_k}(i_d^{p_d}-\lvert i_d+i_{d-1}\lvert^{p_d}) + \bar{C}}
{(\lvert i_1 \lvert^{p_1}+\ldots+\lvert i_d+i_{d-1}\lvert^{p_d}+1)
(\lvert i_1\lvert^{p_1}+\ldots+i_d'^{p_d}+1)}\right) \right\lvert,
\end{equation}
where
\begin{eqnarray*}
\bar{C} 
\!\!\!&:=&\!\!\!  
i_d^{p_d}\lvert i_d' + i_{d-1}\lvert^{p_d} - \lvert i_d+i_{d-1} \lvert^{p_d} i_d'^{p_d}
+i_d^{p_d}- i_d'^{p_d} +\lvert i_d'+i_{d-1}\lvert^{p_d} - \lvert i_d+i_{d-1}\lvert^{p_d}\\
\!\!\!&=&\!\!\! i_d^{p_d} \big[ \lvert i_d' + i_{d-1}\lvert^{p_d} - i_d'^{p_d} \big] 
\!\!+\! i_d'^{p_d} \big[ i_d^{p_d} \!-\! \lvert i_d+i_{d-1} \lvert^{p_d} \big]
\!+\! \big[ i_d^{p_d} \!-\! \lvert i_d+i_{d-1}\lvert^{p_d} \big]  
\!+\! \big[ \lvert i_d'+i_{d-1}\lvert^{p_d} - i_d'^{p_d} \big]\!.
\end{eqnarray*}
Since the expression 
$$\frac{(\lvert i_1\lvert^{p_1}+\ldots+i_d^{p_d}+1)
(\lvert i_1\lvert^{p_1}+\ldots+\lvert i_d'+i_{d-1}\lvert^{p_d}+1)}
{(\lvert i_1\lvert^{p_1}+\ldots+\lvert i_d+i_{d-1}\lvert^{p_d}+1)
(\lvert i_1\lvert^{p_1}+\ldots+i_d'^{p_d}+1)}$$ 
is uniformly bounded from below by a positive number, 
(\ref{a-r}) is bounded from above by
$$M\left\lvert\frac{\sum_{k=1}^{d-1}\lvert i_k\lvert^{p_k}(\lvert i_d'+i_{d-1}\lvert^{p_d}-i_d'^{p_d})+
\sum_{k=1}^{d-1}\lvert i_k\lvert^{p_k}(i_d^{p_d}-\lvert i_d+i_{d-1}\lvert^{p_d})+\bar{C}}
{(\lvert i_1\lvert^{p_1}+\ldots+\lvert i_d+i_{d-1}\lvert^{p_d}+1)
(\lvert i_1\lvert^{p_1}+\ldots+i_d'^{p_d}+1)}\right\lvert.$$
By The Mean Value Theorem, and since $p_{d-1}\!>\!p_d$, this last expression is smaller than or equal to 
$$M\frac{\sum_{k=1}^{d-1}\lvert i_k\lvert^{p_k}(i_d'+\lvert i_{d-1}\lvert)^{p_d-1} \lvert i_{d-1}\lvert 
+\sum_{k=1}^{d-1}\lvert i_k\lvert^{p_k}(i_d+\lvert i_{d-1}\lvert)^{p_d-1}\lvert i_{d-1}\lvert+\bar{C}'}
{(\lvert i_1\lvert^{p_1}+\ldots+i_d^{p_d}+1)(\lvert i_1\lvert^{p_1}+\ldots+i_d'^{p_d}+1)},$$
where $\bar{C}'$ equals 
$$i_d^{p_d}(i_d'+\lvert i_{d-1}\lvert)^{p_d-1}\lvert i_{d-1}\lvert+i_d'^{p_d}
(i_d+\lvert i_{d-1}\lvert)^{p_d-1}\lvert i_{d-1}\lvert+(i_d+\lvert i_{d-1}\lvert)^{p_d-1}
\lvert i_{d-1}\lvert+(i_d'+\lvert i_{d-1}\lvert)^{p_d-1}\lvert i_{d-1}\lvert.$$
Therefore, in order to get an upper bound for the value of \esp \esp \esp 
$\frac{\vert \log Df_d(x) - \log Df_d(y) \vert}{\vert x-y\vert^{\alpha}}$, 
\esp \esp \esp we only need to do so with 
\begin{equation}\label{a-voir}
\frac{\lvert i_k\lvert^{p_k}(i_d'+\lvert i_{d-1}\lvert)^{p_d-1}\lvert 
i_{d-1}\lvert}{(\lvert i_1\lvert^{p_1}+\ldots+i_d^{p_d}+1)(\lvert i_1\lvert^{p_1}
+\ldots+i_d'^{p_d}+1)\vert x-y\vert^{\alpha}}, 
\quad \mbox{ where }1 \leq k\leq n, 
\end{equation}
and
\begin{equation}\label{a-manger}
\frac{i_d'^{p_d}(i_d+\lvert i_{d-1}\lvert)^{p_d-1}\lvert i_{d-1}\lvert}{(\lvert i_1\lvert^{p_1}
+\ldots+i_d^{p_d}+1)(\lvert i_1\lvert^{p_1}+\ldots+i_d'^{p_d}+1)\vert x-y\vert^{\alpha}}.
\end{equation}

\vsp

Expression (\ref{a-voir}) is easy to deal with. Indeed, since  
\begin{equation}\label{simple}
\vert x-y\vert \esp 
\geq \esp (i_d'-i_d-1)\vert I_{i_1,i_2,\ldots,i_d'}\vert \esp
= \esp \left(\frac{i_d'-i_d-1}{\lvert i_1\lvert^{p_1}+\ldots+i_d'^{p_d}+1}\right),
\end{equation}
we have
\begin{multline*}
\frac{\lvert i_k\lvert^{p_k}(i_d'+\lvert i_{d-1}\lvert)^{p_d-1}\lvert i_{d-1}\lvert}
{(\lvert i_1\lvert^{p_1}+\ldots+i_d^{p_d}+1)(\lvert i_1\lvert^{p_1}+\ldots+i_d'^{p_d}+1)\vert x-y\vert^{\alpha}}
\leq \\ \leq \frac{\lvert i_k\lvert^{p_k}(i_d'+\lvert i_{d-1}\lvert)^{p_d-1}\lvert i_{d-1}\lvert}
{(\lvert i_1\lvert^{p_1}+\ldots+i_d^{p_d}+1)(\lvert i_1\lvert^{p_1}+\ldots+i_d'^{p_d}+1)^{1-\alpha}} 
\leq \\ \leq \frac{(i_d'+\lvert i_{d-1}\lvert)^{p_d-1}\lvert i_{d-1}\lvert}
{(\lvert i_1\lvert^{p_1}+\ldots+i_d'^{p_d}+1)^{1-\alpha}}.
\end{multline*}
To estimate the right-side expression, we consider two cases. If, on the one hand, 
we have \esp $i_d'^{p_d}\leq\lvert i_{d-1}\lvert^{p_{d-1}}$, \esp then 
$$\frac{(i_d'+\lvert i_{d-1}\lvert)^{p_d-1}\lvert i_{d-1}\lvert}
{(\lvert i_1\lvert^{p_1}+\ldots+i_d'^{p_d}+1)^{1-\alpha}}
\hspace{0.1cm}\leq\hspace{0.1cm}
\frac{(\lvert i_{d-1}\lvert^{\frac{p_{d-1}}{p_d}}+\lvert i_{d-1}\lvert)^{p_d-1}\lvert i_{d-1}\lvert}
{(\lvert i_{d-1}\lvert^{p_{d-1}}+1)^{1-\alpha}}.$$ 
This is uniformly bounded when \esp $\frac{p_{d-1}}{p_d}(p_d-1)+1\leq p_{d-1}(1-\alpha)$, 
\esp which is equivalent to condition $\mathrm{(v_B)}$. On the other hand, 
if \esp $\lvert i_{d-1}\lvert^{p_{d-1}} \leq i_d'^{p_d}$, then \esp 
$$\frac{(i_d' + \lvert i_{d-1}\lvert)^{p_d-1}\lvert i_{d-1}\lvert}
{(\lvert i_1\lvert^{p_1}+\ldots+i_d'^{p_d}+1)^{1-\alpha}}
\hspace{0.1cm} \leq \hspace{0.1cm}
\frac{(i_d'+(i_d')^{\frac{p_d}{p_{d-1}}})^{p_d-1}(i_d')^{\frac{p_d}{p_{d-1}}}}
{(i_d'^{p_d}+1)^{1-\alpha}},$$ 
and the right-side term is uniformly bounded 
provided that condition $\mathrm{(v_B)}$ holds.\\

To obtain an upper bound for (\ref{a-manger}), we will consider separately 
the cases (a), (b), (c) and (d) of the previous two sections.\\

In case (a) we have \esp $i_d \leq i_d'\leq 2i_d+1$. \esp Hence, the upper bound already 
obtained for (\ref{a-voir}) with $k\!=\!d$ is an upper bound for (\ref{a-manger}).\\ 

In case (b), we have \esp 
$i_d'^{p_d}\leq\lvert i_1\lvert^{p_1}+\ldots+\lvert i_{d-1}\lvert^{p_{d-1}}$. \esp 
Hence, (\ref{a-manger}) is smaller than or equal to 
$$\sum_{k=1}^{d-1} \frac{\lvert i_k\lvert^{p_k}(i_d'+\lvert i_{d-1}\lvert)^{p_d-1}\lvert i_{d-1}\lvert}
{(\lvert i_1\lvert^{p_1}+\ldots+i_d^{p_d}+1)(\lvert i_1\lvert^{p_1}+\ldots+i_d'^{p_d}+1)
\vert x-y\vert^{\alpha}},$$
and we have already seen that each term of this sum is uniformly bounded.\\

\vsp

In case (c), we use (\ref{la-de-c}) to obtain
$$\frac{i_d'^{p_d}(i_d+\lvert i_{d-1}\lvert)^{p_d-1}\lvert i_{d-1}\lvert}
{(\lvert i_1\lvert^{p_1}+\ldots+i_d^{p_d}+1)(\lvert i_1\lvert^{p_1}+\ldots+i_d'^{p_d}+1)
\vert x-y\vert^{\alpha}} 
\leq 
M \esp \frac{(i_d+\lvert i_{d-1}\lvert)^{p_d-1}\lvert i_{d-1}\lvert }
{(\lvert i_1\lvert^{p_1}+\ldots+i_d^{p_d}+1)^{1-\alpha}} \esp 
\frac{(i_d+1)^{(p_d-1)\alpha}}{(\lvert i_1\lvert^{p_1}+\ldots+i_d^{p_d}+1)^{\alpha}}.$$
In the right-side expression, the second factor \esp 
$\frac{(i_d+1)^{(p_d-1)\alpha}}{(i_d^{p_d} + 1)^{\alpha}}$ 
is uniformly bounded. To show that the same holds with the first factor, 
one may proceed as at the end of the estimates for (\ref{a-voir}) 
just changing $i_d'$ by $i_d$.\\

\vsp

Finally, in case (d), the estimate (\ref{clef}) 
shows that (\ref{a-manger}) is smaller than or equal to 
$$\frac{(i_d+\lvert i_{d-1}\lvert)^{p_d-1}\lvert i_{d-1}\lvert 
(i_d + 1 + S^{1/p_d})^{\alpha(p_d-1)} \esp (i'_d + S^{1/p_d})^{\alpha}}
{(\lvert i_1\lvert^{p_1}+\ldots+i_d^{p_d}+1)(i_d' - i_d - 1)^{\alpha}}.$$
Since the condition \esp $1 + i_d^{p_d} \leq S $ \esp yields  
$i_d \leq S^{1/p_d}$, this expression is smaller than or equal to 
$$M \esp \esp \frac{(i_d + \lvert i_{d-1}\lvert)^{p_d-1} 
\lvert i_{d-1}\lvert (i_d'+S^{1/p_d})^{\alpha}} {(i_d'-i_d-1)^{\alpha}} 
\esp\esp S^{^{\frac{\alpha(p_d-1)}{p_d} - 1}}.$$
Moreover, by the definition of $S$, we have 
\esp $|i_{d-1}| \leq S^{1/p_{d-1}} \leq S^{1/p_d}$, \esp which 
shows that the last expression  is smaller than or equal to 
$$M \esp \esp \frac{(i_d'+S^{1/p_d})^{\alpha}} {(i_d'-i_d-1)^{\alpha}} 
\esp\esp S^{^{\frac{p_d-1}{p_d} + \frac{1}{p_{d-1}} + \frac{\alpha(p_d-1)}{p_d} - 1}}.$$
Because of the conditions \esp $i_d'\geq 2i_d+2$ \esp and \esp $S \leq 1 + i_d'^{p_d}$, 
\esp this last expression is bounded from above by 
$$MS^{\frac{p_d-1}{p_d} + \frac{1}{p_{d-1}} + \frac{\alpha(p_d-1)}{p_d} - 1},$$
which is uniformly bounded by the condition $\mathrm{(v_B)}$.

\vsp\vsp 

\noindent{\bf III.} Finally, in the case where $x \in I_{i_1,\ldots,i_d}$ 
and $y \in I_{i_1',\ldots,i_d'}$ for different $(i_1,\ldots,i_{d-1})$ and 
$(i_1',\ldots,i_{d-1}')$, one may apply the same argument of the previous maps. 


\section{On a family of metabelian subgroups of $\dac$}
\label{last}

\hspace{0.45cm} For each couple of integers $(i,j)$, let $I_{i,j}$ be an interval of length 
$|I_{i,j}|$ so that the sum $\sum_{i,j} \vert I_{i,j}\vert$ is finite. Joining these intervals 
lexicographically, we obtain a closed interval $I$. Following \cite[\S 2.3]{FF}, we will deal 
with a particular family of nilpotent groups $N_d$ acting on $I$. Each $N_d$ has 
nilpotence degree $d+1$ and is metabelian. Moreover, $N_1$ coincides with the 
Heisenberg group $N_2$. 

The group $N_d$ has a presentation 
$$\big\langle f,g_0,g_1,\ldots,g_d \! : [g_i,g_j] = id, \esp \esp [f,g_0] = id, 
\esp \esp [f,g_i] = g_{i-1} \mbox{ for all } i \geq 1 \big\rangle.$$ 
As maps, the generators send each interval $I_{i,j}$ 
into a certain $I_{i',j'}$ and coincide with the diffeomorphism 
$\varphi_{I_{i,j},I_{i,j-1}}^{I_{i',j'},I_{i',j'-1}}$ therein. The map $f$ sends 
$I_{i,j}$ into $I_{i+1,j}$. The maps $g_0$ and $g_1$ send $I_{i,j}$ into $I_{i,j+1}$ 
or $I_{i,j+i},$ respectively. To describe the dynamics of $g_2,\ldots,g_d$, for each 
$0 < k \leq d$ and each $i \in \mathbb{Z}$, we let 
$$r_{k}(i) = \frac{i(i+1)(i+2)\ldots (i+k-1)}{k!},$$ 
and we define $r_0 (i) = 1$ for all $i.$ 
(Note that $\vert r_{k}(i)\vert \leq \vert i\vert^{k}$ for $k > 0$.) 
Then the element $g_k$ sends the interval $I_{i,j}$ into $I_{i,j+r_k(i)}$. 

Now fix a positive number $\alpha < 1$. To carry out the preceding construction so that the 
resulting maps are $C^{1+\alpha}$-diffeomorphisms of $I$, we need to make a careful choice of 
the lengths $|I_{i,j}|$. We let $q > 1$ be such that the following conditions are satisfied:

\vsp

\noindent (i$_{\mathrm{C}}$) \esp \esp $1<q<2$,

\vsp

\noindent (ii$_{\mathrm{C}}$) \esp \esp $\alpha < 2-q$,

\vsp

\noindent (iii$_{\mathrm{C}}$) \esp \esp $\alpha < \frac{q}{2q-1}$,

\vsp

\noindent (iv$_{\mathrm{C}}$) \esp \esp $\alpha < \frac{1}{q}$.

\vsp\vsp

Note that since $\alpha < 1$ and the preceding right-side expressions go to 1 or 
to infinity as $q$ tends to 1 from above, we may choose $q$ very near to $1$ 
in such a way that these conditions are fulfilled.

\vsp

Now let \esp $p := \frac{2q-1}{q-1}.$ \esp Clearly, we may also impose the following 
supplementary conditions:

\vsp

\noindent (v$_{\mathrm{C}}$) \esp \esp $p > dq$,

\vsp

\noindent (vi$_{\mathrm{C}}$) \esp \esp $\alpha \leq \frac{1}{q}-\frac{d}{p}$,

\vsp

\noindent (vii$_{\mathrm{C}}$) \esp \esp $\alpha < \frac{1}{q-1}-\frac{dq}{2q-1}$. 
 
\vsp\vsp

We then define 
$$\vert I_{i,j}\vert:=\frac{1}{\vert i \vert^{p}+\vert j\vert^{q}+1}.$$ 
Since \hspace{0.04cm} $1/p + 1/ q < 1$, \hspace{0.04cm} it follows from \cite[\S 3]{KN} 
that the sum $\sum_{i,j} \vert I_{i,j}\vert$ is finite. We claim that the group 
$N_d$ obtained by using the maps from \S \ref{denjoy-pixton} is formed 
by $C^{1+\alpha}$-diffeomorphisms of $I$. 


\subsection{The map $f$ is a $C^{1+\alpha}$-diffeomorphism}
 
\hspace{0.45cm} For simplicity, we only deal with points in the intervals 
$I_{i,j}$ with $i \geq 0$ and $j \geq 0$ (the other cases are analogous).

\vsp

First we consider $x,y$ in the same interval $I_{i,j}$. We have
$$\frac{\vert \log Df(x) - \log Df(y) \vert}{\vert x-y\vert}\leq
\frac{M}{\vert I_{i,j}\vert}\left\vert\frac{\vert I_{i,j}\vert}
{\vert I_{i+1,j}\vert}\frac{\vert I_{i+1,j-1}\vert}{\vert I_{i,j-1}\vert}-1\right\vert.$$
Hence, 
$$\frac{\vert \log Df(x) - \log Df(y) \vert}{\vert x-y\vert^{\alpha}}
\leq\frac{M}{\vert I_{i,j}\vert}\left\vert\frac{\vert I_{i,j}\vert}{\vert I_{i+1,j}\vert}
\frac{\vert I_{i+1,j-1}\vert}{\vert I_{i,j-1}\vert}-1\right\vert\vert I_{i,j}\vert^{1-\alpha}
= M\left\vert\frac{\vert I_{i,j}\vert}{\vert I_{i+1,j}\vert}\frac{\vert I_{i+1,j-1}\vert}
{\vert I_{i,j-1}\vert}-1\right\vert\vert I_{i,j}\vert^{-\alpha}.$$
This yields 
$$\frac{\vert \log Df(x) - \log Df(y) \vert}{\vert x-y\vert^{\alpha}}\leq 
M\left\vert\ \frac{(i+1)^{p}+j^{q}+1}{i^{p}+j^{q}+1}\frac{i^{p}+(j-1)^{q}+1}
{(i+1)^{p}+(j-1)^{q}+1}-1 \right\vert(i^{p}+j^{q}+1)^{\alpha}.$$
Therefore, the value of $\frac{\vert \log Df(x) - \log Df(y) \vert}{\vert x-y\vert^{\alpha}}$ 
is bounded from above by 
$$M\left\vert\ \frac{((i+1)^{p}+j^{q}+1)(i^{p}+(j-1)^{q}+1)-(i^{p}+j^{q}+1)
((i+1)^{p}+(j-1)^{q}+1)}{(i^{p}+j^{q}+1)^{1-\alpha}((i+1)^{p}+(j-1)^{q}+1)}\right\vert,$$
which equals 
$$M \frac{\big( (i+1)^p - i^p \big) \big( j^q - (j-1)^q \big)}
{(i^{p}+j^{q}+1)^{1-\alpha}((i+1)^{p}+(j-1)^{q}+1)}.$$
By the Mean Value Theorem, this expression is bounded from above by 
$$M \frac{i^{p-1}j^{q-1}}{(i^{p}+j^{q}+1)^{1-\alpha}((i+1)^{p}+(j-1)^{q}+1)}.$$
Thus 
$$\frac{\vert \log Df(x) - \log Df(y) \vert}{\vert x-y\vert^{\alpha}} 
\esp \leq \esp M\frac{i^{p-1}}{(i+1)^{p}}\frac{j^{q-1}}{j^{q(1-\alpha)}}.$$
Now notice that the last expression is uniformly bounded when 
$q-1\leq q(1-\alpha)$, which is satisfied by condition (iv$_{\mathrm{C}}$). 

\vsp

Next we consider $x \in I_{i,j}$ and $y \in I_{i,j'}$, with $j < j'$. 
The definition of $f$ and property (\ref{uno-tsuboi}) yield 
$$\log Df(x) = \log D\phi^{a',a}_{b',b} (x-w) = 
\log \frac{b}{a} + \log D\psi_{\log(b'a/a'b)} \left(\frac{x-w}{a}\right),$$
where $I_{i,j}=[w,w+a]$, $I_{i,j-1}=[w+a',w]$, $I_{i+1,j}=[w',w'+b]$, and $I_{i+1,j-1}=[w'+b',w']$. 
Analogously, 
$$\log Df(y) = \log D\phi^{c',c}_{d',d}(y-u) = 
\log \frac{d}{c}+\log D\psi_{\log(d'c/c'd)}\left(\frac{y-u}{c}\right),$$ 
where $I_{i,j'}=[u,u+c]$, $I_{i,j'-1}=[u+c',u]$, $I_{i+1,j'}=[u',u'+d]$, 
and $I_{i+1,j'-1}=[u'+d',u']$. By property (\ref{dos-tsuboi}),  
\begin{eqnarray*}
\big\vert \log Df(x) - \log Df(y) \big\vert 
&\leq& \left\lvert\log \frac{b}{a}-\log \frac{d}{c}\right\lvert + 
\left\lvert\log D\psi_{\log(b'a/a'b)}\left(\frac{x-w}{a}\right)\right\lvert + 
\left\lvert\log D\psi_{\log(d'c/c'd)}\left(\frac{y-u}{c}\right)\right\lvert \\
&\leq& \left\lvert\log \frac{b}{a}-\log \frac{d}{c}\right\lvert + 
\left\lvert\log \frac{b'}{a'}-\log \frac{b}{a}\right\lvert + 
\left\lvert\log \frac{d'}{c'}-\log \frac{d}{c}\right\lvert.
\end{eqnarray*}
Note that the last expression corresponds to 
$$\left \vert \log \left( \frac{|I_{i+1,j}|}{|I_{i,j}|} \right)
 - \log \left( \frac{|I_{i+1,j'}|}{|I_{i,j'}|}\right) 
\right\vert+\left \vert \log \left( \frac{|I_{i+1,j-1}|}{|I_{i,j-1}|} \right)
 - \log \left( \frac{|I_{i+1,j}|}{|I_{i,j}|}\right) 
\right\vert+\left \vert \log \left( \frac{|I_{i+1,j'-1}|}{|I_{i,j'-1}|} \right)
 - \log \left( \frac{|I_{i+1,j'}|}{|I_{i,j'}|}\right) \right\vert.$$
Since the function $j \mapsto \frac{|I_{i+1,j}|}{|I_{i,j}|}$ is non-decreasing, 
the preceding inequality yields 
\begin{eqnarray*}
\big\vert \log Df(x) - \log Df(y) \big\vert 
&\leq& 
3 \left \vert \log \left( \frac{|I_{i+1,j'}|}{|I_{i,j'}|} \right)
 - \log \left( \frac{|I_{i+1,j-1}|}{|I_{i,j-1}|}\right) \right\vert\\ 
&=& 3 \left\vert \log \frac{\vert I_{i+1,j'} \vert}{\vert I_{i,j'} \vert}
\frac{\vert I_{i,j-1} \vert}{\vert I_{i+1,j-1} \vert} \right\vert\\
&=& 
3 \left\vert \log \frac{\big( i^{p}+j'^{q}+1\big)}{\big((i+1)^{p}+j'^{q}+1\big)}
\frac{\big((i+1)^{p}+(j-1)^{q}+1\big)}{\big(i^{p}+(j-1)^{q}+1\big)}\right\vert.
\end{eqnarray*}
Hence, the value of 
$\big\vert \log Df(x) - \log Df(y) \big\vert $ 
is bounded from above by
\begin{equation}
M \left\vert \log \left( 1+\frac{(i^{p}+j'^{q}+1)
((i+1)^{p}+(j-1)^{q}+1)-((i+1)^{p}+j'^{q}+1)
(i^{p}+(j-1)^{q}+1)}{((i+1)^{p}+j'^{q}+1)
(i^{p}+(j-1)^{q}+1)} \right) \right\vert.
\label{expre}
\end{equation}
Since $\frac{ i^{p}+j'^{q}+1}{(i+1)^{p}+j'^{q}+1} \frac{(i+1)^{p}+(j-1)^{q}+1}{i^{p}+(j-1)^{q}+1}$ is uniformly 
bounded from below, namely 
$$\frac{ i^{p}+j'^{q}+1}{(i+1)^{p}+j'^{q}+1} \frac{(i+1)^{p}+(j-1)^{q}+1}{i^{p}+(j-1)^{q}+1} 
\geq \frac{ i^{p}+j'^{q}+1}{2^{p}i^{p}+j'^{q}+1} 
\geq \frac{ i^{p}+j'^{q}+1}{2^{p}i^{p}+2^{p}j'^{q}+2^{p}} 
= \frac{1}{2^{p}},$$
the expression in (\ref{expre}) is bounded from above by 
$$M \bigg\vert\frac{(i^{p}+j'^{q}+1)((i+1)^{p}+(j-1)^{q}+1)-((i+1)^{p}+j'^{q}+1)
(i^{p}+(j-1)^{q}+1)}{((i+1)^{p}+j'^{q}+1)(i^{p}+(j-1)^{q}+1)}\bigg\vert,$$
which equals 
$$M \bigg\vert\frac{(j'^{q}-(j-1)^{q})((i+1)^{p}-i^{p})}{((i+1)^{p}+j'^{q}+1)(i^{p}+(j-1)^{q}+1)}\bigg\vert.$$
By the Mean Value Theorem, this expression is bounded from above by 
$$M \frac{i^{p-1}j'^{q-1}(j'-j+1)}{((i+1)^{p}+j'^{q}+1)(i^{p}+(j-1)^{q}+1)}.$$
Therefore,
\begin{equation}
\big\vert \log Df(x) - \log Df(y) \big\vert \leq 
M \frac{i^{p-1}j'^{q-1}(j'-j)}{(i^{p}+j'^{q})(i^{p}+j^{q})}.
\label{expresion-1}
\end{equation}
We will split the general case into the following four cases:

\vspace{0.15cm}

\noindent (a) \esp \esp $j'\leq 2j+1$,

\vspace{0.15cm}

\noindent (b) \esp \esp $j'^{q}\leq i^{p}$,

\vspace{0.15cm}

\noindent (c) \esp \esp $j'>2j+1 \, , j'^{q}>i^{p} \, , j^{q} \geq i^{p}$,

\vspace{0.15cm}

\noindent (d) \esp \esp $j'>2j+1 \, , j'^{q}>i^{p} \, , j^{q} < i^{p}$. 

\vspace{0.15cm}

In cases (a) and (b), notice that from 
$$\vert x-y \vert \esp \geq \esp (j'-j-1)I_{i,j'}$$ 
it follows that
$$\vert x-y \vert^{\alpha} \esp 
\geq \esp \left( \frac{j'-j-1}{i^{p}+j'^{q}+1} \right)^{\alpha}.$$
Hence by (\ref{expresion-1}),
$$\frac{\vert \log Df(x) - \log Df(y) \vert}{\vert x-y \vert^{\alpha}} \esp \leq \esp  
M \frac{i^{p-1}j'^{q-1}(j'-j)(i^{p}+j'^{q}+1)^{\alpha}}{(i^{p}+j'^{q})(i^{p}+j^{q})(j'-j-1)^{\alpha}},$$
that is,
\begin{equation}
\frac{\vert \log Df(x) - \log Df(y) \vert}{\vert x-y \vert^{\alpha}} \esp \leq \esp 
M \frac{i^{p-1}j'^{q-1}(j'-j)^{1-\alpha}}{(i^{p}+j'^{q})^{1-\alpha}(i^{p}+j^{q})}.
\label{primera}
\end{equation}

\vsp

In case (a), we have $j' \leq 2j+1$, and hence the right-side of (\ref{primera}) 
is bounded from above by 
$$M \frac{i^{p-1}j^{q-1}j^{1-\alpha}}{(i^{p}+j^{q})^{1-\alpha}(i^{p}+j^{q})} =
M \frac{i^{p-1}j^{q-\alpha}}{(i^{p}+j^{q})^{2-\alpha}}.$$
On the one hand, if $i \leq j^{\frac{1}{p-1}}$, then this expression is bounded by 
$\frac{j^{q-\alpha+1}}{j^{q(2-\alpha)}} = j^{(\alpha-1)(q-1)}$. 
Since $\alpha < 1$, the expression is uniformly bounded.
On the other hand, if $j \leq i^{p-1}$, then we have the upper bound 
$$\frac{i^{p-1+(p-1)(q-\alpha)}}{i^{p(2-\alpha)}} = i^{\alpha - 1 - p - q + pq},$$ 
and this expression is uniformly bounded by the condition (ii$_{\mathrm{C}}$).\\

In case (b), we have $j'^{q}\leq i^{p}$, and hence the right-side expression in (\ref{primera}) 
is bounded from above by 
$$\frac{i^{p-1}i^{\frac{p}{q}(q-1)}i^{\frac{p}{q}(1-\alpha)}}{i^{p+p(1-\alpha)}} 
= i^{\alpha (p - \frac{p}{q}) - 1},$$
which is uniformly bounded by condition (iii$_{\mathrm{C}}$).\\

In case (c), we have
$$\vert x-y \vert \geq \sum_{j<n<j'} I_{i,n} = 
\sum_{j<n<j'} \frac{1}{i^{p}+n^{q}+1}\geq \sum_{j<n<j'} \frac{1}{j^{q}+n^{q}+1} 
\geq \sum_{j<n<j'} \frac{1}{3n^{q}}\geq \int_{j+1}^{j'}\frac{dx}{3x^{q}},$$
and hence
$$\vert x-y \vert \geq M \frac{1}{(j+1)^{q-1}} 
\Big( 1- \big( \frac{j+1}{j'} \big)^{q-1} \Big) \geq \frac{1}{(j+1)^{q-1}} 
\Big( 1- \big( \frac{1}{2} \big)^{q-1} \Big).$$
Therefore, 
\begin{equation}
\vert x-y \vert^{\alpha} \geq M \frac{1}{(j+1)^{(q-1)\alpha}},
\label{uf}
\end{equation}
and this yields
$$\frac{\vert \log Df(x) - \log Df(y) \vert}{\vert x-y \vert^{\alpha}} \leq 
M \frac{i^{p-1}j'^{q-1}(j'-j)j^{(q-1)\alpha}}{(i^{p}+j'^{q})(i^{p}+j^{q})} = 
M \Big( \frac{j'^{q-1}(j'-j)}{i^{p}+j'^{q}} \Big) 
\Big( \frac{i^{p-1}j^{(q-1)\alpha}}{i^{p}+j^{q}} \Big).$$
In the last expression, the first factor is uniformly bounded, 
while the second one is bounded by 
$$M \hspace{0.05cm} \frac{j^{\frac{q}{p}(p-1)}j^{(q-1)\alpha}}{j^{q}}.$$ 
This last expression is uniformly bounded when $\frac{q}{p}(p-1)+(q-1)\alpha \leq q$, 
which is ensured by the condition (iii$_{\mathrm{C}}$).\\

The last case (d) is
$$j'>2j+1 \, , j'^{q}>i^{p} \, , j^{q} < i^{p}.$$
For the distance between $x$ and $y$ we now have the estimate 
$$\vert x-y \vert \geq \sum_{j<n<j'} I_{i,n}=\sum_{j<n<j'} \frac{1}{i^{p}+n^{q}+1} \geq 
\int_{j+1}^{j'} \frac{1}{i^{p}+x^{q}+1}dx \geq \int_{j+1}^{j'} \frac{1}{(i^{\frac{p}{q}}+x+1)^{q}}dx.$$
The last integral is essentially
$$M \hspace{0.05cm} 
\frac{(i^{\frac{p}{q}}+j'+1)^{q-1}-(i^{\frac{p}{q}}+j+2)^{q-1}}{(i^{\frac{p}{q}}+j+2)^{q-1}(i^{\frac{p}{q}}+j'+1)^{q-1}},$$
and by the Mean Value Theorem, this is larger than
$$M \hspace{0.05cm} \frac{j'-j-1}{(i^{\frac{p}{q}}+j'+1)^{2-q}(i^{\frac{p}{q}}+j+2)^{q-1}(i^{\frac{p}{q}}+j'+1)^{q-1}}.$$
Therefore, 
\begin{equation}
\vert x-y \vert ^{\alpha} \esp \esp \geq \esp \esp 
M \frac{(j'-j-1)^{\alpha}}{(i^{\frac{p}{q}}+j'+1)^{\alpha}(i^{\frac{p}{q}}+j+2)^{(q-1)\alpha}}.
\label{vamos}
\end{equation}
This yields 
\begin{eqnarray*}
\frac{\vert \log Df(x) - \log Df(y) \vert}{\vert x-y \vert^{\alpha}} 
&\leq& M \frac{i^{p-1}j'^{q-1}(j'-j)(i^{\frac{p}{q}}+j'+1)^{\alpha}
(i^{\frac{p}{q}}+j+2)^{(q-1)\alpha}}{(i^{p}+j'^{q})(i^{p}+j^{q})(j'-j-1)^{\alpha}}\\
&\leq& M \frac{i^{p-1}(i^{\frac{p}{q}}+j'+1)^{\alpha}
(i^{\frac{p}{q}}+j+2)^{(q-1)\alpha}}{(i^{p}+j^{q})(j'-j-1)^{\alpha}}\\ 
&\leq& M \frac{i^{p-1}(2j'+1)^{\alpha}
(2i^{\frac{p}{q}}+2)^{(q-1)\alpha}}{(i^{p}+j^{q})(\frac{j'}{2})^{\alpha}}\\
&\leq& M \frac{i^{p-1}(i^{\frac{p}{q}}+1)^{(q-1)\alpha}}{(i^{p}+j^{q})},
\end{eqnarray*}
and by condition (iii$_{\mathrm{C}}$) the last expression is uniformly bounded.

This completes the proof of the $C^{1+\alpha}$ regularity 
of $f$. Similar arguments apply to its inverse $f^{-1}$, thus 
showing that $f$ is a $C^{1 + \alpha}$-diffeomorphism of $I$.


\subsection{Each map $g_k$ is a $C^{1+\alpha}$-diffeomorphism}

\hspace{0.45cm} Again, we will only consider the case of positive 
$i,j$. First, we take $x,y$ in the same interval $I_{i,j}$. We have
$$\frac{\vert \log Dg_k(x) - \log Dg_k(y)\vert}{\vert x-y\vert} \esp \leq \esp  
\frac{M}{\vert I_{i,j}\vert}\left\vert\frac{\vert I_{i,j}\vert}{\vert I_{i,j+r_k(i)}\vert} 
\frac{\vert I_{i,j+r_k(i)-1}\vert}{\vert I_{i,j-1}\vert}-1\right\vert.$$
Hence 
\begin{eqnarray*}
\frac{\vert \log Dg_k(x) - \log Dg_k(y)\vert}{\vert x-y\vert^{\alpha}} 
&\leq& 
\frac{M}{\vert I_{i,j}\vert}\left\vert\frac{\vert I_{i,j}\vert}{\vert I_{i,j+r_k(i)}\vert} 
\frac{\vert I_{i,j+r_k(i)-1}\vert}{\vert I_{i,j-1}\vert}-1\right\vert\vert I_{i,j}\vert^{1-\alpha}\\ 
&=& 
M \esp \left\vert\frac{\vert I_{i,j}\vert}{\vert I_{i,j+r_k(i)}\vert} 
\frac{\vert I_{i,j+r_k(i)-1}\vert}{\vert I_{i,j-1}\vert}-1\right\vert\vert I_{i,j}\vert^{-\alpha}\\
&\leq& M \esp \left\vert\ \frac{i^{p}+(j+r_k(i))^{q}+1}{i^{p}+j^{q}+1} 
\frac{i^{p}+(j-1)^{q}+1}{i^{p}+(j+r_k(i)-1)^{q}+1}-1 \right\vert(i^{p}+j^{q}+1)^{\alpha}.
\end{eqnarray*}
The last expression may be rewritten as 
$$M \left \vert \frac{(i^{p}+(j+r_k(i))^{q}+1)(i^{p}+(j-1)^{q}+1)-(i^{p}+j^{q}+1)(i^{p}+(j+r_k(i)-1)^{q}+1)}
{(i^{p}+j^{q}+1)^{1-\alpha}(i^{p}+(j+r_k(i)-1)^{q}+1)}\right \vert.$$
By the Mean Value Theorem, the value of this expression is bounded from above by 
$$M \esp \frac{i^{p}j^{q-1}+i^{p}(j+r_k(i))^{q-1}+j^{q-1}+(j+r_k(i))^{q-1} + 
(j+r_k(i))^{q}j^{q-1}+(j+r_k(i))^{q-1}j^{q}}{(i^{p}+j^{q})^{1-\alpha}(i^{p}+(j+r_k(i)-1)^{q})},$$
and hence by 
$$M \esp \frac{i^{p}(j+r_k(i))^{q-1}+(j+r_k(i))^{q}j^{q-1}+(j+r_k(i))^{q-1}j^{q}}{(i^{p}+j^{q})^{2-\alpha}} 
\leq M \esp \frac{i^{p}(j+i^{k})^{q-1}+(j+i^{k})^{q}j^{q-1}+(j+i^{k})^{q-1}j^{q}}{(i^{p}+j^{q})^{2-\alpha}}.$$
We claim that the preceding right-expression is uniformly bounded. Indeed, if $i^{p}\leq j^{q}$, then 
it is smaller than or equal to 
$$M\frac{j^{q}(j+j^{\frac{qk}{p}})^{q-1}+(j+j^{\frac{qk}{p}})^{q}j^{q-1} + 
(j+j^{\frac{qk}{p}})^{q-1}j^{q}}{j^{q(2-\alpha)}},$$ 
which is uniformly bounded by the conditions (iv$_{\mathrm{C}}$) and (v$_{\mathrm{C}}$). 
If $j^{q}\leq i^{p}$, then it is smaller than or equal to 
$$M\frac{i^{p}(i^{\frac{p}{q}}+i^{k})^{q-1} + 
(i^{\frac{p}{q}}+i^{k})^{q}i^{\frac{p}{q}(q-1)}+(i^{\frac{p}{q}}+i^{k})^{q-1}i^{p}}{i^{p(2-\alpha)}},$$ 
which is again uniformly bounded by the conditions (iv$_{\mathrm{C}}$) and (v$_{\mathrm{C}}$).\\

Next we consider the case where $x \in I_{i,j}$ and $y \in I_{i,j'}$, with $j \leq j'$. In 
this case, $\vert \log \, Dg_k(x) - \log \, Dg_k(y) \vert$ is smaller than or equal to 
$$\left \vert \log \frac{\vert I_{i,j+r_k(i)}\vert}{\vert I_{i,j}\vert} - 
\log \frac{\vert I_{i,j'+r_k(i)} \vert}{\vert I_{i,j'}\vert}\right \vert + 
\left \vert \log \frac{\vert I_{i,j+r_k(i)-1}\vert}{\vert I_{i,j-1}\vert} - 
\log \frac{\vert I_{i,j+r_k(i)} \vert}{\vert I_{i,j}\vert}\right \vert + 
\left \vert \log \frac{\vert I_{i,j'+r_k(i)-1}\vert}{\vert I_{i,j'-1}\vert} - 
\log \frac{\vert I_{i,j'+r_k(i)} \vert}{\vert I_{i,j'}\vert}\right \vert.$$
The estimates for the last two terms are similar to those above, 
and we leave the computations to the reader. The first term equals
$$\left \vert \log \frac{\vert I_{i,j}\vert}{\vert I_{i,j+r_k(i)}\vert} 
\frac{\vert I_{i,j'+r_k(i)} \vert}{\vert I_{i,j'}\vert}\right \vert 
= \left \vert \log \frac{i^{p}+j'^{q}+1}{i^{p} + (j'+r_k(i))^{q}+1} 
\frac{i^{p}+(j+r_k(i))^{q}+1}{i^{p}+j^{q}+1} \right \vert,$$
that is, 
$$\left \vert \log \left (1+\frac{(i^{p}+j'^{q}+1)(i ^{p}+(j+r_k(i))^{q}+1) - 
(i^{p}+(j'+r_k(i))^{q}+1)(i^{p}+j^{q}+1)}{(i^{p}+(j'+r_k(i))^{q}+1)
(i^{p}+j^{q}+1)} \right)\right\vert.$$
We claim that the expression $\frac{i^{p}+j'^{q}+1}{i^{p} + (j'+r_k(i))^{q}+1} 
\frac{i^{p} + (j+r_k(i))^{q}+1}{i^{p}+j^{q}+1}$ is bounded from below by a 
positive number. Indeed, the first factor is uniformly bounded because: 

\vsp 

\noindent -- if $j'^{q} \leq i^{p}$, then $\frac{i^{p}+j'^{q}+1}{i^{p}+(j'+r_k(i))^{q}+1} \geq 
\frac{i^{p}+j'^{q}+1}{i^{p}+(j' + r_k(i))^{q}+1} \geq \frac{i^{p}+1}{i^{p}+(i^{\frac{p}{q}}+i^{k})^{q}+1}$, 
and the last expression is uniformly bounded from below by a positive number;

\vsp

\noindent -- if $i^{p}\leq j'^{q}$, then $\frac{i^{p}+j'^{q}+1}{i^{p}+(j'+r_k(i))^{q}+1} \geq 
\frac{i^{p}+j'^{q}+1}{i^{p}+(j'+i^{k})^{q}+1} \geq \frac{j'^{q}+1}{j'^{q}+(j'+j'^{\frac{qk}{p}})^{q}+1}$, 
which is uniformly bounded from below by a positive number. 

\vsp \vsp

\noindent The second factor is uniformly bounded as well because: 

\noindent -- if $i^{p} \leq j^{q}$, then $0\leq j-j^{\frac{qk}{p}}\leq j+r_{k}(i)$, 
thus $\frac{i^{p}+(j+r_k(i))^{q}+1}{i^{p}+j^{q}+1} \geq \frac{i^{p}+(j-j^{\frac{qk}{p}})^{q}+1}{i^{p}+j^{q}+1} 
\geq \frac{(j-j^{\frac{qk}{p}})^{q}+1}{2j^{q}+1}$, which is uniformly bounded from below by a positive number;

\vsp

\noindent -- if $j^{q} \leq i^{p}$, then $\frac{i^{p} + (j+r_k(i))^{q} + 1}{i^{p}+j^{q}+1} \geq 
\frac{i^{p}+1}{2 i^{p}+1}$, which is uniformly bounded from below by a positive number. 

\vsp\vsp

From the preceding, we deduce the estimate 
\begin{small}
$$\left \vert \log \frac{\vert I_{i,j+r_k(i)}\vert}{\vert I_{i,j}\vert} -\log \frac{\vert I_{i,j'+r_k(i)} \vert}
{\vert I_{i,j'}\vert}\right \vert 
\leq M\left \vert\frac{(i^{p}+j'^{q}+1)(i^{p} + (j+r_k(i))^{q} + 1) - 
(i^{p} + (j'+r_k(i))^{q}+1)(i^{p}+j^{q}+1)}{(i^{p} + (j'+r_k(i))^{q}+1) (i^{p}+j^{q}+1)} \right\vert.$$
\end{small}Using the Mean Value Theorem and the inequalities $j < j'$ and $r_k(i)\leq i^{k}$, 
the right-side term above is easily seen to be smaller than or equal to 
\begin{equation}
M\frac{i^{p+k}(j'+i^{k})^{q-1}+j^{q}(j'+i^{k})^{q-1}i^{k}+j'^{q}(j+i^{k})^{q-1}i^{k}}{(i^{p}+j'^{q})(i^{p}+j^{q})}.
\label{ultima}
\end{equation}
To get an upper bound for this expression, we separately consider the cases (a), (b), (c), and (d), from 
the previous section.\\

The first case (a) is $j'\leq 2j+1$. We have \esp 
$\vert x-y \vert \geq \frac{j'-j-1}{i^{p}+j'^{q}+1}$, and hence 
\begin{eqnarray*}
\frac{\vert \log Dg_k(x) - \log Dg_k(y) \vert}{\vert x-y \vert^{\alpha}} 
\!&\leq&\! M \esp \frac{i^{p+k}(j'+i^{k})^{q-1}+j^{q}(j'+i^{k})^{q-1}i^{k}+j'^{q}(j+i^{k})^{q-1}i^{k}}
{(i^{p}+j'^{q})(i^{p}+j^{q})}\frac{(i^{p}+j'^{q}+1)^{\alpha}}{(j'-j-1)^{\alpha}} + M\\
&\leq&\! M \esp \frac{i^{p+k}(j'+i^{k})^{q-1}+j^{q}(j'+i^{k})^{q-1}i^{k}+j'^{q}(j+i^{k})^{q-1}i^{k}}
{(i^{p}+j'^{q})^{1-\alpha}(i^{p}+j^{q})} + M\\
&\leq&\! M \esp \frac{i^{p+k}(j+i^{k})^{q-1} + 2j^{q}(j+i^{k})^{q-1}i^{k}}
{(i^{p}+j'^{q})^{1-\alpha}(i^{p}+j^{q})} + M. 
\end{eqnarray*}
We will deal with the expressions $\frac{i^{p+k}(j+i^{k})^{q-1}}{(i^{p}+j'^{q})^{1-\alpha}(i^{p}+j^{q})}$ 
and $\frac{j^{q}(j+i^{k})^{q-1}i^{k}}{(i^{p}+j'^{q})^{1-\alpha}(i^{p}+j^{q})}$ separately. For the first 
we have
$$\frac{i^{p+k}(j+i^{k})^{q-1}}{(i^{p}+j'^{q})^{1-\alpha}(i^{p}+j^{q})} \leq 
\frac{i^{p+k}(j+i^{k})^{q-1}}{(i^{p}+j^{q})^{2-\alpha}}.$$
Now notice that 

\vsp

\noindent -- if $i^{p} \leq j^{q}$, then $\frac{i^{p+k}(j+i^{k})^{q-1}}{(i^{p}+j^{q})^{2-\alpha}} 
\leq \frac{j^{\frac{q}{p}(p+k)}(j+j^{\frac{qk}{p}})^{q-1}}{j^{q(2-\alpha)}}$, and 
this is bounded when $\frac{q}{p}(p+k)+q-1 \leq q(2-\alpha)$, that is, when 
$\alpha \leq \frac{1}{q}-\frac{k}{p}$, which is our condition (vi$_{\mathrm{C}}$);

\vsp

\noindent -- if $j^{q} \leq i^{p}$, then $\frac{i^{p+k}(j+i^{k})^{q-1}}{(i^{p}+j^{q})^{2-\alpha}} 
\leq \frac{i^{p+k}(i^{\frac{p}{q}}+i^{k})^{q-1}}{i^{p(2-\alpha)}}$, and this is bounded when 
$p+k+\frac{p}{q}(q-1) \leq p(2-\alpha)$, that is, when $\alpha \leq \frac{1}{q}-\frac{k}{p}$.

\vsp

\noindent For the second expression we have 
$$\frac{j^{q}(j+i^{k})^{q-1}i^{k}}{(i^{p}+j'^{q})^{1-\alpha}(i^{p}+j^{q})} 
\leq \frac{j^{q}(j+i^{k})^{q-1}i^{k}}{(i^{p}+j^{q})^{2-\alpha}}.$$
Again, notice that 

\vsp

\noindent -- if $i^{p}\leq j^{q}$, then $\frac{j^{q}(j+i^{k})^{q-1}i^{k}}{(i^{p}+j^{q})^{2-\alpha}} 
\leq \frac{j^{q}(j+j^{\frac{qk}{p}})^{q-1}j^{\frac{qk}{p}}}{j^{q(2-\alpha)}}$, and as before, 
this is bounded when $\alpha\leq\frac{1}{q}-\frac{k}{p}$;

\vsp

\noindent -- if $j^{q}\leq i^{p}$, then $\frac{j^{q}(j+i^{k})^{q-1}i^{k}}{(i^{p}+j^{q})^{2-\alpha}} 
\leq \frac{i^{p}(i^{\frac{p}{q}}+i^{k})^{q-1}i^{k}}{i^{p(2-\alpha)}}$, and as before, 
this is bounded when $\alpha \leq \frac{1}{q}-\frac{k}{p}$.\\

The second case (b) is $j'^{q}\leq i^{p}$. The inequality 
$\vert x-y \vert \geq \frac{j'-j-1}{i^{p}+j'^{q}+1}$ yields 
\begin{eqnarray*}
\frac{\vert \log Dg_k(x) - \log Dg_k(y) \vert}{\vert x-y \vert^{\alpha}} 
\!&\leq&\! M \esp \frac{i^{p+k}(j'+i^{k})^{q-1}+j^{q}(j'+i^{k})^{q-1}i^{k}+j'^{q}(j+i^{k})^{q-1}i^{k}}
{(i^{p}+j'^{q})(i^{p}+j^{q})}\frac{(i^{p}+j'^{q}+1)^{\alpha}}{(j'-j-1)^{\alpha}} + M\\ 
&\leq&\! M \esp \frac{i^{p+k}(j'+i^{k})^{q-1}+j^{q}(j'+i^{k})^{q-1}i^{k}+j'^{q}(j+i^{k})^{q-1}i^{k}}
{(i^{p}+j^{q})^{2-\alpha}} + M\\
&\leq&\! M \esp \frac{i^{p+k}(i^{\frac{p}{q}}+i^{k})^{q-1}+j^{q}
(i^{\frac{p}{q}}+i^{k})^{q-1}i^{k}+i^{p}(j+i^{k})^{q-1}i^{k}}{(i^{p}+j^{q})^{2-\alpha}} + M.
\end{eqnarray*}
To estimate the last expression, we will bound the following three expressions: 
$$\frac{i^{p+k}(i^{\frac{p}{q}}+i^{k})^{q-1}}{(i^{p}+j^{q})^{2-\alpha}}, \quad 
\frac{j^{q}(i^{\frac{p}{q}}+i^{k})^{q-1}i^{k}}{(i^{p}+j^{q})^{2-\alpha}}, \hspace{0.5cm} 
\mbox{ and } \hspace{0.5cm} \frac{i^{p}(j+i^{k})^{q-1}i^{k}}{(i^{p}+j^{q})^{2-\alpha}}.$$ 
For the first we have  
$$\frac{i^{p+k}(i^{\frac{p}{q}}+i^{k})^{q-1}}{(i^{p}+j^{q})^{2-\alpha}} 
\esp \leq \esp \frac{i^{p+k}(i^{\frac{p}{q}}+i^{k})^{q-1}}{i^{p(2-\alpha)}},$$ 
and the right-side member is bounded provided that $p+k+\frac{p}{q}(q-1) \leq p(2-\alpha)$, 
that is, $\alpha\leq \frac{1}{q}-\frac{k}{p}$, which is our condition (vi$_{\mathrm{C}}$). 
For the second expression, notice that 

\vsp

\noindent -- if $i^{p}\leq j^{q}$, then
$\frac{j^{q}(i^{\frac{p}{q}}+i^{k})^{q-1}i^{k}}{(i^{p}+j^{q})^{2-\alpha}} 
\leq \frac{j^{q}(j+j^{\frac{qk}{p}})^{q-1}j^{\frac{qk}{p}}}{j^{q(2-\alpha)}}$, and this 
is bounded when $q+\frac{qk}{p}+q-1\leq q(2-\alpha)$, that is, $\alpha\leq \frac{1}{q}-\frac{k}{p}$; 

\vsp

\noindent -- if $j^{q}\leq i^{p}$, then $\frac{j^{q}(i^{\frac{p}{q}}+i^{k})^{q-1}i^{k}}{(i^{p}+j^{q})^{2-\alpha}} 
\leq \frac{i^{p}(i^{\frac{p}{q}}+i^{k})^{q-1}i^{k}}{i^{p(2-\alpha)}}$, and the last term is bounded 
because $\alpha \leq \frac{1}{q} - \frac{k}{p}$.

\vsp

\noindent For the third expression, we have that  

\vsp

\noindent -- if $i^{p}\leq j^{q}$, then $\frac{i^{p}(j+i^{k})^{q-1}i^{k}}{(i^{p}+j^{q})^{2-\alpha}} 
\leq \frac{j^{q}(j+j^{\frac{qk}{p}})^{q-1}j^{\frac{qk}{p}}}{j^{q(2-\alpha)}},$ which is bounded 
when $\alpha \leq \frac{1}{q} - \frac{k}{p}$;

\vsp

\noindent -- if $j^{q}\leq i^{p}$, then $\frac{i^{p}(j+i^{k})^{q-1}i^{k}}{(i^{p}+j^{q})^{2-\alpha}} 
\leq \frac{i^{p}(i^{\frac{p}{q}}+i^{k})^{q-1}i^{k}}{i^{p(2-\alpha)}}$, which is also bounded when 
$\alpha \leq \frac{1}{q}-\frac{k}{p}$.\\ 

\vsp

The third case (c) is $j'>2j+1 \, , \esp j'^{q}>i^{p} \, , \esp j^{q} \geq i^{p}$. Using (\ref{uf}) 
we obtain 
$$\frac{\vert \log \, Dg_k(x) - \log \, Dg_k(y) \vert}{\vert x-y \vert^{\alpha}} \leq 
M \esp \frac{i^{p+k}(j'+i^{k})^{q-1}+j^{q}(j'+i^{k})^{q-1}i^{k}+j'^{q}(j+i^{k})^{q-1}i^{k}}
{(i^{p}+j'^{q})(i^{p}+j^{q})}j^{(q-1)\alpha} + M.$$
To estimate the preceding right-side expression, we deal separately with \esp 
$$\frac{i^{p+k}(j'+i^{k})^{q-1}}{(i^{p}+j'^{q})(i^{p}+j^{q})}j^{(q-1)\alpha}, \quad 
\frac{j^{q}(j'+i^{k})^{q-1}i^{k}}{(i^{p}+j'^{q})(i^{p}+j^{q})}j^{(q-1)\alpha}, \hspace{0.3cm} 
\mbox{ and } \quad \frac{j'^{q}(j+i^{k})^{q-1}i^{k}}{(i^{p}+j'^{q})(i^{p}+j^{q})}j^{(q-1)\alpha}.$$ 
For the first of these expressions one has 
$$\frac{i^{p+k}(j'+i^{k})^{q-1}}{(i^{p}+j'^{q})(i^{p}+j^{q})}j^{(q-1)\alpha} 
\leq \frac{i^{k}(j'+i^{k})^{q-1}j'^{(q-1)\alpha}}{i^{p}+j'^{q}} 
\leq \frac{j'^{\frac{qk}{p}}(j'+j'^{\frac{qk}{p}})^{q-1}j'^{(q-1)\alpha}}{j'^{q}},$$
and the right-side term is bounded by condition (vii$_{\mathrm{C}}$). 
For the second expression one has
$$\frac{j^{q}(j'+i^{k})^{q-1}i^{k}}{(i^{p}+j'^{q})(i^{p}+j^{q})}j^{(q-1)\alpha} 
\leq \frac{(j'+i^{k})^{q-1}i^{k}j'^{(q-1)\alpha}}{i^{p}+j'^{q}} \leq 
\frac{(j'+j'^{\frac{qk}{p}})^{q-1}j'^{\frac{qk}{p}}j'^{(q-1)\alpha}}{j'^{q}},$$
and the right-side term is uniformly bounded by the condition (vii$_{\mathrm{C}}$). 
Finally, for the third expression one has
$$\frac{j'^{q}(j+i^{k})^{q-1}i^{k}}{(i^{p}+j'^{q})(i^{p}+j^{q})}j^{(q-1)\alpha} 
\leq \frac{(j+i^{k})^{q-1}i^{k}j^{(q-1)\alpha}}{i^{p}+j^{q}} 
\leq \frac{(j+j^{\frac{qk}{p}})^{q-1}j^{\frac{qk}{p}}j^{(q-1)\alpha}}{j^{q}},$$
and the right-side term is also bounded by condition (vii$_{\mathrm{C}}$).\\ 

The last case (d) is $j'>2j+1 \, , j'^{q}>i^{p} \, , j^{q} < i^{p}$. Inequality (\ref{vamos}) shows 
that \esp \esp $\frac{\vert \log \, Dg_k(x) - \log \, Dg_k(y) \vert}{\vert x-y \vert^{\alpha}}$ 
\esp \esp is smaller than or equal to
$$M \frac{i^{p+k}(j'+i^{k})^{q-1}+j^{q}(j'+i^{k})^{q-1}i^{k}+j'^{q}(j+i^{k})^{q-1}i^{k}}{(i^{p}+j'^{q})(i^{p}+j^{q})} 
\frac{(i^{\frac{p}{q}}+j'+1)^{\alpha}(i^{\frac{p}{q}}+j+2)^{(q-1)\alpha}}{(j'-j-1)^{\alpha}} + M.$$
In this expression, the term $\frac{(i^{\frac{p}{q}}+j'+1)^{\alpha}}{(j'-j-1)^{\alpha}}$ is bounded by 
$\frac{(2j'+1)^{\alpha}}{(\frac{j'}{2})^{\alpha}}$, and hence it is uniformly bounded. Therefore, 
$$\frac{\vert \log \, Dg_k(x) - \log \, Dg_k(y) \vert}{\vert x-y \vert^{\alpha}} \leq
M\frac{i^{p+k}(j'+i^{k})^{q-1}+j^{q}(j'+i^{k})^{q-1}i^{k}+j'^{q}(j+i^{k})^{q-1}i^{k}}
{(i^{p}+j'^{q})(i^{p}+j^{q})}(i^{\frac{p}{q}}+j)^{(q-1)\alpha}.$$
To estimate the right-side expression, we will deal separately with 
$$\frac{i^{p+k}(j'+i^{k})^{q-1}} {(i^{p}+j'^{q})(i^{p}+j^{q})}(i^{\frac{p}{q}}+j)^{(q-1)\alpha}, 
\quad \frac{j^{q}(j'+i^{k})^{q-1}i^{k}} {(i^{p}+j'^{q})(i^{p}+j^{q})}(i^{\frac{p}{q}}+j)^{(q-1)\alpha}, 
\hspace{0.1cm} \mbox{ and } \hspace{0.1cm} \frac{j'^{q}(j+i^{k})^{q-1}i^{k}}
{(i^{p}+j'^{q})(i^{p}+j^{q})}(i^{\frac{p}{q}}+j)^{(q-1)\alpha}.$$ 
For the first of these expressions one has
$$\frac{i^{p+k}(j'+i^{k})^{q-1}}
{(i^{p}+j'^{q})(i^{p}+j^{q})}(i^{\frac{p}{q}}+j)^{(q-1)\alpha}\leq\frac{i^{k}(j'+i^{k})^{q-1}}
{(i^{p}+j'^{q})}(i^{\frac{p}{q}}+j)^{(q-1)\alpha}\leq\frac{j'^{\frac{qk}{p}}(j'+j'^{\frac{qk}{p}})^{q-1}}
{j'^{q}}(j'+j')^{(q-1)\alpha},$$
and, as before, we know that the right-side term is uniformly bounded by the condition 
(vii$_{\mathrm{C}}$). For the second expression one has
$$\frac{j^{q}(j'+i^{k})^{q-1}i^{k}}
{(i^{p}+j'^{q})(i^{p}+j^{q})}(i^{\frac{p}{q}}+j)^{(q-1)\alpha}\leq\frac{(j'+i^{k})^{q-1}i^{k}}
{(i^{p}+j'^{q})}(i^{\frac{p}{q}}+j)^{(q-1)\alpha}\leq\frac{(j'+j'^{\frac{qk}{p}})^{q-1}j'^{\frac{qk}{p}}}
{j'^{q}}(j'+j')^{(q-1)\alpha},$$
and the right-side term is uniformly bounded by the condition (vii$_{\mathrm{C}}$). Finally, 
for the third expression one has
$$\frac{j'^{q}(j+i^{k})^{q-1}i^{k}}
{(i^{p}+j'^{q})(i^{p}+j^{q})}(i^{\frac{p}{q}}+j)^{(q-1)\alpha}\leq\frac{(j+i^{k})^{q-1}i^{k}}
{(i^{p}+j^{q})}(i^{\frac{p}{q}}+j)^{(q-1)\alpha}\leq\frac{(i^{\frac{p}{q}}+i^{k})^{q-1}i^{k}}
{i^{p}}(i^{\frac{p}{q}}+i^{\frac{p}{q}})^{(q-1)\alpha}.$$
Here the right-side term is bounded when \esp $\frac{p}{q}(q-1)+k+\frac{p}{q}(q-1)\alpha \leq p$, 
\esp which is true by the condition (vii$_{\mathrm{C}}$).


\begin{footnotesize}

\vspace{0.2cm}

\noindent Gonzalo Castro, Andr\'es Navas\\

\noindent Dpto de Matem\'atica y C.C., Fac. de Ciencia, USACH\\

\noindent Alameda 3363, Estaci\'on Central, Santiago, Chile\\

\noindent gonzalo.castro@usach.cl, andres.navas@usach.cl

\vspace{0.2cm}

\noindent Eduardo Jorquera\\

\noindent Inst. de Matem\'aticas, Pontificia Universidad Cat\'olica de Valpara\'{\i}so\\

\noindent Blanco Viel 596, Carro Bar\'on, Valpara\'{\i}so, Chile\\

\noindent eduardo.jorquera@ucv.cl

\end{footnotesize}

\end{document}